\documentclass{amsart}
\usepackage{amssymb,latexsym,graphics}
\usepackage[mathscr]{eucal}

\numberwithin{equation}{section}

\newtheorem{thm}{Theorem}[section]

\newtheorem{lem}[thm]{Lemma}
\newtheorem{prop}[thm]{Proposition}

\theoremstyle{remark}
\newtheorem{rem}[thm]{Remark}

\theoremstyle{definition}

\newcommand{\RR}{ {\mathbb{R}^3} }

\newcommand{\lp}[2]{\Vert \, #1 \, \Vert_{#2}}
\newcommand{\llp}[2]{|||  \, #1 \, |||_{#2}}
\newcommand{\cBox}{\Box^{\mathbb{C}}}
\newcommand{\bL}{\underline{L}}
\newcommand{\bu}{\underline{u}}
\newcommand{\Div}{(\hbox{div})}
\newcommand{\tr}{(\hbox{trace})}
\newcommand{\sym}{(\hbox{symm})}
\newcommand{\sF}{ {{}^\star \! F} }
\newcommand{\gi}{ {{}^I \widetilde{ g }}         }
\newcommand{\gii}{   {{}^{I\! I} \widetilde{g  }}   }
\newcommand{\cqi}{   {   {}^I\widetilde{Q}  }   }
\newcommand{\cqii}{ {   {}^{I\! I} \widetilde{Q}   }   }
\newcommand{\nablai}{ {{}^I\widetilde{\nabla}}  }
\newcommand{\nablaii}{ {{}^{I\! I}\widetilde{\nabla}}  }

\newcommand{\sD}{ {\, {\slash\!\!\!\! D  }}}
\newcommand{\sA}{ {\, {\slash\!\!\!\! A  }}}

\newcommand{\bv}{ \underline{v}}

\newcommand{\tp}{\tau_+}
\newcommand{\tm}{\tau_-}
\newcommand{\tz}{\tau_0}
\newcommand{\balpha}{ \underline{\alpha}}

\newcommand{\tF}{\widetilde{F}}
\newcommand{\oF}{\overline{F}}
\newcommand{\sJ}{ {\, {\slash\!\!\! J  }} }
\newcommand{\td}{\widetilde}
\newcommand{\snabla}{ {\slash\!\!\! \nabla} }
\newcommand{\TK}{\overline{K}_0^s}
\newcommand{\LC}{ \mathcal{L}^\mathbb{C} }
\newcommand{\cdiv}{ \hbox{(div)}^\mathbb{C} }

\newcommand{\ret}{\vspace{.3cm}}

\begin{document}

\title[Stability for CSF on Minkowski Space]
{Global Stability for Charge--Scalar Fields on Minkowski Space}
\author{Hans Lindblad and Jacob Sterbenz\\ \\ University of California at San Diego}
\email{lindblad@math.ucsd.edu} \email{jsterben@math.ucsd.edu}
\date{}
\dedicatory{}
\commby{}

%%% ----------------------------------------------------------------------

\begin{abstract}
We prove that the charge-scalar field equations are globally stable
on $(3+1)$ dimensional Minkowski space for small initial data in
certain gauge covariant weighted Sobolev spaces. These spaces can be
chosen to be almost scale invariant with respect to the
homogeneity of the equations, and our result is valid for initial
data with non-zero charge that is also non-stationary at space-like
infinity. The method of proof is a tensor-geometric approach which
is based on a sharp family of weighted bilinear $L^2$ space-time
estimates.
\end{abstract}

%%% ----------------------------------------------------------------------
\maketitle
%%% ----------------------------------------------------------------------
%%%%%%%%%%%%%%%%%%%%%%%%%%%%%%%%%%%%%%%%%%%%%%%%%%%%%%%%%%%%%%%%%%%%%%%%%%
%%%-----------------------------------------------------------------------

%-------------------------------------------------------------------------
%-------------------------------------------------------------------------
%-------------------------------------------------------------------------
%-------------------------------------------------------------------------
%-------------------------------------------------------------------------
%%%%%%%%%%%%%%%%%%%%%%%%%%%%%%%%%%%%%%%%%%%%%%%%%%%%%%%%%%%%%%%%%%%%%%%%%%
%-------------------------------------------------------------------------

\section{Introduction}

In this paper, we study the global in time behavior of small
amplitude classical solutions to the massless  Charge--Scalar--Field
equations (CSF)  equations on Minkowski space. These are also
sometimes referred to in the literature as Maxwell--Higgs or the
Maxwell--Klein--Gordon equations. They are constructed as follows: \
We let $\mathcal{M}$ denote the $(3+1)$ dimensional Minkowski space
with metric $g = (-1,1,1,1)$ and compatible Levi-Civita connection
which we refer to here as $\nabla$. Then let:
\begin{equation}
    V \ = \ \mathcal{M}\times\mathbb{C} \ , \notag
\end{equation}
denote a complex line bundle over $\mathcal{M}$ with
hermitian inner product $\langle \cdot , \cdot \rangle_V $ and
compatible connection $D$. That is, one has the formula:
\begin{equation}
    X \langle \psi, \phi\rangle_V \ = \
    \langle D_X \psi, \phi\rangle_V +
    \langle \psi, D_X \phi\rangle_V \ , \label{bundle_metric}
\end{equation}
for  vector-fields $X$ on $\mathcal{M}$. The connection $D$ can be
extended in a natural way to sections of the tensor bundle $V\otimes
T^k_l\mathcal M$ of complex valued tensors of type $(k,l)$ on
Minkowski space such that one has the Leibnitz rule:
\begin{equation}
        D_X( \phi\otimes T) \ = \ D_X(\phi)\otimes T + \phi\otimes
        \nabla_X \phi \ . \notag
\end{equation}
In this setup, we call a section $\phi$ of $V$  a \emph{complex
scalar field} if it satisfies the equation:
\begin{equation}
    \Box^{\mathbb{C}}\phi \ = \
    D^\alpha D_\alpha \phi \ = \ 0 \ , \label{complex_field}
\end{equation}
where $D^2\phi$ is the covariant Hessian of $\phi$, which is
computed regarding $D\phi$ as a complex valued one-form on
$\mathcal{M}$. Since $V$ is globally trivial, we can let
$\textbf{1}_V$ denote a unit normalized global section. That is one
has $\langle \textbf{1}_V , \textbf{1}_V\rangle \equiv 1 $. Such a
choice is referred to as a \emph{gauge}. In terms of this, any
section $\phi$ to $V$ can be identified with a complex valued
function on $\mathcal{M}$ where we abusively write $\phi = \phi
\textbf{1}_V$. This allows us to define complex conjugation
$\overline{\phi}$ as one would for ordinary complex numbers. We will
use this operation to define a symmetric and symplectic inner
product on $V$ as follows:
\begin{align}
        \Re \langle\phi,\overline{\psi}\rangle_V \ &= \
        \frac{1}{2}(\phi\overline{\psi} + \overline{\phi}\psi) \ ,
        &\Im \langle\phi,\overline{\psi}\rangle_V \ &= \
        \frac{1}{2}(\phi\overline{\psi} - \overline{\phi}\psi) \ .
        \notag
\end{align}
Also, with respect to this frame for $V$, and for any frame
$\{e_\alpha\}$ on $\mathcal{M}$, we define the connection one-form
$\{A_\alpha\}$ via the relations:
\begin{equation}
    D_\alpha \textbf{1}_V \ = \ \sqrt{-1} A_\alpha \ . \label{A_con_def}
\end{equation}
Notice that the compatibility condition \eqref{bundle_metric}
immediately implies that the $\{A_\alpha\}$ are real. Using this
notation, the covariant derivative of any complex valued tensor $T$
with respect to a vector-field $X$ in $\mathcal{M}$ can be written
as:
\begin{equation}
    D_X T \ = \ \nabla_X(T) + \sqrt{-1}
    A(X) \cdot T \ . \label{cov_der_A}
\end{equation}
Of course, the choice of $\textbf{1}_V$ and hence representation of
$D_X$ is somewhat
arbitrary because one can always perform a local (or global)
unitary transformation of it. If we let
 $\textbf{1}_V \rightsquigarrow \widetilde{\textbf{1}}_V$
denote the transformation given by
$\widetilde{\textbf{1}}_V = e^{-i\chi} \textbf{1}_V$,
where $\chi$ is some real valued function on $\mathcal{M}$,
then \eqref{cov_der_A} shows that the potentials $\{A_\alpha\}$
transform as:
\begin{equation}
    \widetilde{A}_\alpha \ = \ A_\alpha - e_\alpha(\chi) \ . \notag
\end{equation}
Now it turns out that this ambiguity which is inherent to the
equation \eqref{cov_der_A} does \emph{not} need to be resolved in
order to proceed with a detailed analysis of the scalar field
equation \eqref{complex_field}. This is because the connection
\eqref{bundle_metric} has a basic geometric invariant on which it is
possible to base literally all analysis of \eqref{complex_field} so
long as one is content  to define all analytic objects of interest
in terms of geometric quantities. This is the approach we shall take
here. The basic geometric invariant is, of course, the curvature of
the connection $D$ which  arises from the operation of commuting
covariant differentiation. Specifically, there exists a (real)
two-form $F$ on $\mathcal{M}$, such that for any two vector-fields
$X,Y$ we have the relation:
\begin{equation}
    D_XD_Y \phi - D_YD_X\phi - D_{[X,Y]}\phi \ =\
    \sqrt{-1}F(X,Y)\cdot\phi \ . \label{curvature_def}
\end{equation}
In terms of the frame $\{e_\alpha\}$ we can simply write:
\begin{equation}
    D_\alpha D_\beta \phi - D_\beta D_\alpha \phi \ = \
    \sqrt{-1}F_{\alpha\beta}\cdot\phi \ . \label{index_version}
\end{equation}
Here $D_\alpha D_\beta \phi$ denotes the $(\alpha,\beta)$
component of the Hessian $D^2\phi$ and should not be
confused with the repeated directional covariant
differentiation $D_{e_\alpha}D_{e_\beta}\phi$. To see the difference,
a short computation
of \eqref{curvature_def} in the frame $\{e_\alpha\}$ shows that
one has the identity:
\begin{equation}
    D_{e_\alpha}D_{e_\beta}\phi -
    D_{e_\beta}D_{e_\alpha}\phi - D_{[e_\alpha,e_\beta]}\phi
    \ = \ \sqrt{-1}\big( e_\alpha(A_\beta) - e_\beta(A_\alpha)
     - c_{\alpha\beta}^\gamma A_\gamma \big)\cdot \phi \ , \notag
\end{equation}
where $[e_\alpha,e_\beta] = c^\gamma_{\alpha\beta} e_\gamma$ are the
structure ``constants'' of the frame $\{e_\alpha\}$. In terms of the
indexed version \eqref{index_version}, using \eqref{cov_der_A} this
last line reads:
\begin{equation}
        D_\alpha D_\beta \phi - D_\beta D_\alpha \phi \ = \
        \sqrt{-1}(\nabla_\alpha A_\beta - \nabla_\beta
        A_\alpha)\cdot\phi \ , \notag
\end{equation}
Either way, this gives the well known ``Bianchi identity'':
\begin{equation}
    F \ = \ dA \ . \label{bianchi}
\end{equation}\\

To define the full CSF system, we now couple the complex scalar
field $\phi$ to the curvature of the connection $D$ in such a way
that one ends up with a Lagrangian field theory. For an intuitive
approach to how this is done, we consider $F$ as a solution to
Maxwell's equations. Because $F$ satisfies the Bianchi identity
\eqref{bianchi}, we have that:
\begin{subequations}\label{basic_Maxwell}
\begin{align}
    \nabla^\beta F_{\alpha\beta} \ = \ J_\alpha \ ,
    \label{basic_Maxwell1}\\
    \nabla^\beta \, \sF_{\alpha\beta} \ = \ 0 \ ,
    \label{basic_Maxwell2}
\end{align}
\end{subequations}
for some one-form $J$ which obeys the rule $\nabla^\alpha J_\alpha
=0$. In the above formulae $\sF$ denotes the Hodge dual of $F$ which
is given by the expression $\sF_{\alpha\beta} = \frac{1}{2}
\in_{\alpha\beta}^{\ \ \ \gamma\delta}F_{\gamma\delta}$. Here
$\in_{\alpha\beta\gamma\delta}$ denotes the volume form on Minkowski
space. As usual, we use this duality operation and contraction along
the time-like vector-field $\partial_t$ to define the electric and
magnetic field strengths:
\begin{align}
        E_i \ &= \ F_{0i} \ ,
        &H_i \ &= \ {}^*F_{0i} \ . \label{em_decomp}
\end{align}
Using these quantities, the system \eqref{basic_Maxwell} takes the
following familiar form which will be of use in the sequel:
\begin{subequations}\label{EH_Maxwell}
\begin{align}
    \partial_t E \ - \ \nabla_x\times H \ &= \ \underline{J} \ ,
    &\Div E \ &= \ J_0 \ , \label{EH_Maxwell1}\\
    \partial_t H \ + \ \nabla_x\times E \ &= \ 0 \ ,
    &\Div H \ &= \ 0 \ . \label{EH_Maxwell2}
\end{align}
\end{subequations}
Here $\nabla_x\times = \in_{i}^{\ jk}\partial_j$ is the usual curl
operation, and $\underline{J}$ denotes the spatial part of the
current vector $J$. We note here that this latter form of the
equations will be particularly useful when dealing with questions
concerning the initial data of Lie derivatives of $F_{\alpha\beta}$.
In general it is the more natural language with which to discuss the
initial value problem for the system
\eqref{basic_Maxwell}.\\

For the Maxwell field $F$, there is an energy-momentum tensor
$Q[F]$ which  is given by:
\begin{equation}
    Q_{\alpha\beta}[F] \ = \ \frac{1}{2}\left(
    F_{\alpha\gamma}F_{\beta}^{\ \gamma}
    + \sF_{\alpha\gamma}\, \sF_{\beta}^{\ \gamma}\right) \ = \
    F_{\alpha\gamma}F_{\beta}^{\ \gamma} - \frac{1}{4}g_{\alpha\beta}
    \, F_{\gamma\delta} F^{\gamma\delta} \ . \label{F_em}
\end{equation}
A brief calculation using the identities
\eqref{basic_Maxwell1}--\eqref{basic_Maxwell2}
shows that one has the divergence identity:
\begin{equation}
    \nabla^\alpha Q_{\alpha\beta}[F] \  = \ - \
    F_{\beta\gamma} J^\gamma \ . \label{F_consv}
\end{equation}
Furthermore, the complex scalar field \eqref{complex_field} also has
an energy-momentum tensor which is analogous to the usual
energy-momentum tensor of the (non-covariant) D'Lambertian. This is
given by:
\begin{equation}
    Q_{\alpha\beta}[\phi] \ = \
    \Re ( D_\alpha \phi \overline {D_\beta \phi}) -
    \frac{1}{2}g_{\alpha\beta} \, D^\gamma\phi \overline{
    D_\gamma\phi} \ . \label{phi_em}
\end{equation}
A direct application of the field equation \eqref{complex_field},
the compatibility condition \eqref{bundle_metric},
and the commutator identity \eqref{curvature_def} shows that one has the divergence
law:
\begin{equation}
    \nabla^\alpha Q_{\alpha\beta}[\phi] \ = \
    F_{\beta\gamma} \Im (\phi \overline{D^\gamma\phi})
    \ . \label{phi_consv}
\end{equation}
Notice that this formula does not contain any explicit reference to
the connection coefficients \eqref{A_con_def}. This is because it is a
purely tensorial identity, and in particular does not depend on a
choice of coordinates or frame.\\

The natural choice of coupling $F$ to $\phi$ now comes if we simply
stipulate that the current vector on the right hand side of
\eqref{basic_Maxwell1} be given by $J_\alpha = \Im (\phi
\overline{D_\alpha\phi})$. With this extra condition satisfied, the
total energy momentum tensor $Q = Q[F] + Q[\phi]$ becomes divergence
free:
\begin{equation}
    \nabla^\alpha Q_{\alpha\beta} \ = \ 0 \ . \label{total_consv}
\end{equation}
Notice that this coupling is legitimate because as an immediate
consequence of the equation \eqref{complex_field} one has that this
choice of current vector $J_\alpha$ satisfies the continuity
equation:
\begin{equation}
    \nabla^\alpha J_\alpha \ = \ 0 \ , \label{J_cont}
\end{equation}
which is the only prerequisite for a vector $J$ to show up as the
right hand side of the Maxwell equation \eqref{basic_Maxwell1}. We
write the resulting system of equations, which we henceforth refer
to as the CSF equations, together as:
\begin{subequations}\label{basic_MKG}
\begin{align}
        \nabla^\beta F_{\alpha\beta} \ &= \ \Im (\phi \overline
    {D_\alpha\phi}) \ , \label{basic_EM}\\
    \cBox \phi \ &= \ 0 . \label{basic_CW}
\end{align}
\end{subequations}
We note here that implicit in the system \eqref{basic_MKG} is the
Bianchi identity \eqref{bianchi} or \eqref{basic_Maxwell2}. This is
because we are assuming a-priori that $F$ is the curvature of the
connection
which gives rise to $\cBox$.\\

In this work we study the global in time Cauchy problem for the
system \eqref{basic_MKG}. Since we have stated the equations in such
a way as to deemphasize a choice of gauge, this deserves some
explanation. Because we are interested in an evolution problem, we are
reduced to a discussing the notion of gauge covariant initial data.
To specify these, we first let $\underline{D}$ denote a connection
on the initial time slice bundle $\{0\}\times\RR\times\mathbb{C}$.
Since this is embedded as a hyper-surface in the original bundle
$\mathcal{M}\times \mathbb{C} = V$, we will also specify the initial
normal derivative $D_t$. The initial data for the system
\eqref{basic_MKG} can then be written in the form:
\begin{subequations}\label{initial_data}
\begin{align}
        F_{0i}(0)  \ &= \ E_i  \ ,
    &{}^*F_{0i}(0)  \ &= \  H_i  \ , \label{initial_data1}\\
    \phi(0) \ &= \ \phi_0 \ ,
    &D_t \phi\, (0) \ &= \ \dot{\phi}_0 \ . \label{initial_data2}
\end{align}
\end{subequations}
In the above notation, we have used the labels $(E,H)$ to denote
quantities which only depend on $x\in\RR$. These should not be
confused with the space-time $(E,H)$ used in the formulas \eqref{EH_Maxwell}.
Now, from the form of the system \eqref{basic_MKG} it is easy to see that
this initial data cannot be specified freely. It must also satisfy the compatibility
conditions:
\begin{align}
        \nabla^i E_i \ &= \ \Im ( \phi_0 \, \overline{\dot{\phi}_0
    }) \ ,
    &\nabla^i H_i \ &= \ 0 .
    \label{compat_conds}
\end{align}
We will call a data set $(E,H, \phi_0 , \dot{\phi}_0 )$ which
satisfies \eqref{compat_conds} \emph{admissible}. The question we
will be concerned with here is to describe in as detailed a way as possible
the global in time behavior of solutions to the CSF equations whose
initial data satisfies certain natural smallness assumptions. As
usual, these will be stated in terms of regularity. In general, we
define the covariant and gauge covariant weighted Sobolev spaces:
\begin{align}
        \lp{T}{H^{k,s_0}(\RR)}^2 \ &= \ \sum_{|I|\leqslant k}\
    \int_{\RR}\ (1 + r^2)^{s_0 + |I|} \ |\nabla_x^I\, T|^2 \
    dx \ , \label{tensor_initial_sob}\\
    \lp{\psi}{H^{k,s_0}(\RR)}^2 \ &= \ \sum_{|I|\leqslant k}\
    \int_{\RR}\ (1 + r^2)^{s_0 + |I|}\ |\underline{D}^I \psi |^2
    dx \ . \label{scalar_initial_sob}
\end{align}
Here $I$ denotes the usual multiindex notation, while $\nabla_x$
denotes the restriction of the Levi-Civita connection $\nabla$ to
the time slice $\{0\}\times \RR$. $T$ denotes an arbitrary tensor.
We will consider the Cauchy problem for admissible initial data sets
which are
in these function spaces for various values of $s_0$ and $k$.\\

Now, a basic regularity theorem for the system \eqref{basic_MKG}
says that if the norms
\eqref{tensor_initial_sob}--\eqref{scalar_initial_sob} applied to
(admissible) $(E,H,\phi_0,\dot{\phi}_0)$ are finite for certain
values of $s_0$ and $k$, then a global solution to this system
exists. In fact, there is no need to impose any of the weights
$(1+r^2)$, and almost no smoothness is needed to make the argument
work. This is the content of the fundamental regularity result of
Eardly-Moncrief \cite{EM_YM1}--\cite{EM_YM2}, and its later
significant refinement due to Klainerman-Machedon \cite{KM_MKG}.
While these results in some sense give the strongest possible global
existence result one could hope to ask for, they contain surprising
little information as to the nature of the global solution they
obtain. In fact, the only information they provide is that the
unweighted Sobolev $H^s$ norm of the solution remains bounded for
all time. No practical bounds for this quantity are obtained for $1
< s$. Furthermore, these results provide no information on the
profile of the solution, even assuming that the initial data is
localized at time $t=0$ in a way that is consistent with the norms
\eqref{tensor_initial_sob}--\eqref{scalar_initial_sob}. Finally,
these results do not lend themselves to a discussion of the
asymptotic behavior of the scalar field $\phi$ in the sense of $L^2$
scattering. For example, there is no doubt that this involves long
range corrections\footnote{We will not discuss these corrections
here, although they are implicit in the gauge-covariant approach we
use. We believe the analysis we present goes a long way towards
constructing and estimating the effect of long range corrections in
the system \eqref{basic_MKG}, which should take the from of phase
corrections which depend only in the charge. Also, we note here that
the existence of long range corrections to the CSF equations is
consistent with what is known for the Maxwell-Dirac system (see
\cite{SFT_MD}) which is similar.} due to the asymptotic behavior of
the electro-magnetic field $F_{\alpha\beta}$, but the large data
results just mentioned provide
no information as to what form these corrections should take.\\

The result we present here is a first attempt to fill in some of the
gap in understanding the asymptotic behavior of solutions to the
system \eqref{basic_MKG}, at least for initial data with small
\eqref{tensor_initial_sob}--\eqref{scalar_initial_sob} norms where
the specific values of $s_0$ and $k$ are to be specified. The method
we employ is a new variant of the tensorial-geometric approach of
Christodoulou-Klainerman \cite{CK_Fields}--\cite{CK_Ein} which also
uses in a crucial way the space-time energy estimate approach
appearing in the recent work of Lindblad-Rodnianski \cite{LR_Ein}.
This will be applied to both the electro-magnetic field
\eqref{basic_Maxwell} and to the scalar field \eqref{complex_field}
through a novel geometric adaptation of a set of fundamental
estimates which goes back to work of Morawetz \cite{Mor}. The
approach we follow also borrows ideas from previous works on the low
regularity behavior of systems of semilinear wave equations, notably
\cite{KM_MKG}, in that we make crucial use of a certain family of
sharp weighted bilinear $L^2$ space-time estimates for tensorial
contractions to control error terms which come up in our analysis.
Some of these estimates are proved through the use of a weighted
$L^2(L^\infty)$ ``Strichartz'' type estimate which holds for certain
null-components in a tensorial setting.\\

There has been previous work in general\footnote{``General'' meaning
the system \eqref{basic_MKG} with or without charge. There is also
the previous work of Christodoulou--Choquet-Bruhat (\cite{CC_gauge})
which uses conformal compactification. However, this type of
procedure requires  more decay of the initial data and does not
allow for either charge or dipole moments.} on the global asymptotic
behavior of the equations \eqref{basic_MKG}. This is the paper of
Shu \cite{Shu_MKG}. That program, while similar in spirit to ours,
was not carried out in detail and was only aimed at resolving the
case of solutions which are trivial outside of a fixed forward
cone.\footnote{While it is usually the case with scalar wave
equations that once one can prove global existence using weighted
energies for compactly supported initial data, one can move to
non-compactly supported initial data by simply adding enough weights
(at space-like infinity), for tensorial equations of the type
\eqref{basic_MKG} it is \emph{not} possible to directly pass from
the assumption of compactly supported initial data to the case of
initial data with infinite support. This problem persists
\emph{even} if one assumes that the initial data decays at a very
fast polynomial rate. The reason is not entirely obvious at first,
and is due to subtle interaction between the scaling and tensorial
properties of the equations \eqref{basic_MKG}. This will be
explained in a moment.} It should be noted that we have found no way
of closing our argument without making crucial use of space-time
$L^2(L^\infty)$  estimates which serve to eliminate certain
logarithmic divergences one would encounter when trying to pursue an
$L^\infty$ approach. As far as we can see, these types of estimates
are necessary even of one assumes compact support for the initial
data (modulo charge). There is no discussion of space-time norms in
the work \cite{Shu_MKG}, or any other mechanism which would allow
one to circumvent divergences arising in the purely $L^\infty$ and energy
estimates  approach taken there.
However, we would like to call the readers attention to
the fact that the work \cite{Shu_MKG} represents in many ways the
approach we will take here. In particular, some of the more
interesting features of the CSF nonlinearity were first discovered there,
notably the beautiful cancelation which we use in
(a slightly different form on) line
\eqref{special_canc} below.\\

Finally, we also mention the work of Psarrelli \cite{Ps_MKG} on the
asymptotic behavior for the equations \eqref{basic_MKG} with a
non-zero mass. This work again assumes compact support of the
initial data, and is ultimately very different form what we do here
because the version of \eqref{complex_field} with non-zero mass
enjoys a quite different decay estimate than its massless
counterpart. It would be an interesting and non-trivial problem to
extend the work of \cite{Ps_MKG} to the case of fields which are
non-stationary at space-like infinity, and it is not clear if our
method is a step in this direction because of our persistent use of
estimates derived from conformal invariance.\\

Before moving on the statement of our main result, it is worth
mentioning at a heuristic level two of the most difficult features
of attempting a purely physical-space analysis of the equations
\eqref{basic_MKG}. We hope this  gives the reader a bit of insight
into why certain things are necessary. The first problem has to do
with attempting to make use of the so called ``null-condition''
inherent in the system \eqref{basic_MKG}. Unlike other semilinear
geometric wave equations, for example the wave-maps \cite{Sid_WM},
there is not a direct correspondence between the semilinear model
equations studied in \cite{K_semi} and \cite{C_semi} which satisfy
the null condition directly, and the equations of gauge field
theory.\\

There are basically two known ways the null condition for the
equations \eqref{basic_MKG} can be uncovered and utilized.
The first is through
the use of an elliptic gauge, as in the works on the low regularity
properties of this and other related systems
\cite{KM_MKG}--\cite{KM_YM}. This type of procedure leads to
non-local versions of the null-forms studied in  \cite{K_semi} and
\cite{C_semi} which also have entirely different scaling properties than
their non-local versions. To analyze the system of equations which
arise in this way by using the inhomogeneous algebra
\eqref{inh_lor_alg} is quite awkward, and it is not even clear that
this can even be done correctly. For example, the commutator of the
weighted derivatives in \eqref{inh_lor_alg} with Riesz potentials
seems to cause enough of a problem that it effectively negates the
savings one gets from the standard null-forms. Furthermore, any such
analysis would have to somehow take into account the long range
corrections the system \eqref{basic_MKG} must undergo due to the
lack of decay of the electro-magnetic field \eqref{basic_Maxwell}.
These corrections are clearly a guage dependent phenomena, and it
does not seem that the usual elliptic gauge is most convenient for
writing them down.\\

The second way the null condition for the CSF system can uncovered
is directly through its tensorial structure. That is, the null
condition makes itself evident through the contraction structure of
error terms which arise after commuting the equations with
geometrically defined operations of differentiation. To be useful,
the geometric differentiation must take the form of Lie
differentiation of the Maxwell field and gauge covariant
differentiation of the scalar field. Because  the resulting error
terms then involve contractions over all possible indices, they can
be expanded in a frame which takes into account the null geometry of
Minkowski space. In this way, one immediately sees that two ``bad
components'' can never interact with each other. However, since the
system \eqref{basic_MKG} involves the interaction of two different
types of quantities, namely a two form and a scalar field, it is not
entirely clear at first that this type of tensorial structure will
cause enough ``good interactions'' to take place, or even what the
``components'' of the scalar field should be. Fortunately for us,
the CSF equations do contain a deep underlying structure which
allows one to treat \emph{both} quantities, $F_{\alpha\beta}$ and
$\phi$, as a single instance of a master ``tensorial'' object.
What's more, one can produce a single family of bilinear estimates
for this ``tensor'' which covers all possible error terms that can
generated by differentiating the system \eqref{basic_MKG}. This type
of underlying unity goes well beyond the scalar null structure of
\cite{K_semi} and \cite{C_semi}, and is akin to what one sees for
the pure Yang-Mills and the Einstein equations (see
\cite{Shu_YM}--\cite{CK_Ein}).\\

The second main difficulty we need to overcome in our analysis is a
consequence of an interesting interplay between the scaling and
decay properties of solutions to the system \eqref{basic_MKG}. While
in $(3+1)$ dimensions these equation are \emph{subcritical} with
respect to the translation invariant conserved energy \footnote{The
conserved energy is a simple consequence of the tensorial
conservation law \eqref{total_consv} and the fact that $\partial_t$
is a Killing field on Minkowski space. This conserved quantity is at
the level of the $\dot{H}^1$ norm of $\phi$.}, which ultimately is
responsible for large data regularity as mentioned previously, these
equations are in general \emph{critical}\footnote{This should really
be called ``charge critical''. The conservation of charge for these
equations, which can be written as $q = \int_{\mathbb{R}^3} \Im(\phi
\overline{D_t \phi})\equiv const.$, takes place at exactly the
critical regularity which is the $\dot{H}^\frac{1}{2}$ norm of
$\phi$.} with respect to decay at space-like infinity. To better
understand this, consider a linear wave equation of the form:
\begin{equation}
        \Box \phi \ = \ -2\sqrt{-1}\, A^\alpha \partial_\alpha \phi \
        , \label{scalar_model}
\end{equation}
where $\Box= \partial^\alpha \partial_\alpha$ denotes the usual
D'Lambertian on Minkowski space. This can be seen as a model for the
equation \eqref{complex_field} where we have fixed a gauge (which
gauge does not matter for the sake of this discussion). Now, assume
that the initial data for \eqref{scalar_model} decays like $r^{-m}$
as $r\to\infty$. The question we then ask is how
much decay should one then expect from the
potentials $A_\alpha$ in order to guarantee that the solution keeps
this decay in the far exterior region $2t < r$ (which should be the
simplest to control!). By homogeneity, one expects that if
$A_\alpha$ decays like $r^{-l}$ for $2t < r$, then the solution
$\phi$ will decay like $r^{-m -l +1}$ (integrate the second order
equation twice).\footnote{For smooth $A_\alpha$ this heuristic can
be made precise through an easy use of weighted exterior energy
estimates.}\\

One immediately sees from this simple analysis that to control
things $A_\alpha$ must decay at least as well as $r^{-1}$ in the far
exterior, and that for this critical rate of decay some work is
required to avoid logarithmic divergences in the ``error'' term on
the right hand side of \eqref{scalar_model}. This decay rate $r^{-1}$
is not a coincidence, but rather a consequence of the fact that the
equation \eqref{scalar_model} scales like $\dot{H}^\frac{1}{2}$
at the level of initial data, for which one sees that
$r^{-1}$ is precisely the
critical rate of decay. What's more, the entire system
\eqref{basic_MKG} scales like $\dot{H}^\frac{1}{2}$ at the level of
the potentials $A_\alpha$, and because of the tensorial nature of
the CSF equations at least some of the potentials \emph{must} in
general decay no better that $r^{-1}$. This is a simple effect of
the elliptic constraint equation:
\begin{equation}
        \nabla^i E_i \ = \ \Im (\phi\, \overline{D_t \phi}) \ . \label{es_eq1}
\end{equation}
Introducing a potential function $\varphi$ for the curl-free part $E^{cf}$, we
see that unless one has:
\begin{equation}
        \int_{\RR}\ \Im (\phi\, \overline{D_t \phi}) \ \equiv \ 0 \ , \notag
\end{equation}
it will be the case that $\varphi \sim \frac{1}{r}$. Since this
potential function must show up somewhere in field potentials
$A_\alpha$ regardless of the gauge, we see that in general an
equation of the form   \eqref{scalar_model}, and hence
\eqref{complex_field}, will be critical for decay at space-like
infinity.\\

To combat this critical rate of decay, one must prove sharp
space-time estimates in order to have an effective strategy. There
will in general be no extra ``convergence factors'' with which one
can use to integrate over large time intervals. In practice, this
means that all of our space-time estimates need to be proved in
a-priori (divergence) form. This kind of difficulty should be
understood in contrast to other field equations, such as the
Einstein equations. In $(3+1)$ dimensions the long range effect of
the mass for this latter system leaves plenty of room with respect
to the scaling properties of the equations (see for example
\cite{LR_Ein}). Therefore, from the point of view of decay, the
Einstein equations are somewhat more forgiving and do not require
estimates which are as precise as what we need to make things work
here (see for example \cite{LR_Ein}).\\

\noindent We are now ready to state our main result:\\

\begin{thm}[Global Stability of CSF Equations]\label{main_th}
Let $2\leqslant k$ with $k\in\mathbb{N}$, and let $s_0 = s+\gamma$
be given such that $s_0 < \frac{3}{2}$, $0<\gamma$, and
$\frac{1}{2}< s$. Let $(E,H,\phi_0,\dot{\phi}_0)$ be an admissible
initial data set, and define the \emph{charge} to be the value:
\begin{equation}
        q \ =\ \int_{\RR}\ \Im ( \phi_0\, \overline{\dot{\phi}_0}) \ . \label{the_charge}
\end{equation}
Then there exists a universal constant $\mathcal{E}_{k,s,\gamma}$,
which depends only on the parameters $k,s,\gamma$, such that if
$(E,H,\phi_0,\dot{\phi}_0)$ is an admissible initial data set which
satisfies the smallness condition:
\begin{equation}
        \lp{E^{df}}{H^{k,s_0}(\RR)} + \lp{H}{H^{k,s_0}(\RR)} +
    \lp{\underline{D}\phi_0 }{H^{k,s_0}(\RR)} +
    \lp{\dot{\phi}_0}{H^{k,s_0}(\RR)} \ \leqslant \ \mathcal{E}_{k,s,\gamma}
    \ , \label{initial_smallness}
\end{equation}
where $E = E^{df} + E^{cf}$ is the Hodge decomposition of $E$ into
its divergence free and curl free components (resp.), then there
exists a (unique) global solution to the system of equations
\eqref{basic_MKG} with this initial data set such that if
$\{L,\bL,e_A\}$ denotes a standard spherical null frame (see
\eqref{L_bL} and \eqref{ang_frame} below), then the following
point-wise properties of this solution holds:
\begin{subequations}\label{F_peeling}
\begin{align}
        |\alpha| \ &\lesssim \ \mathcal{E}_{k,s,\gamma}\cdot \tp^{-s
        - \frac{3}{2}}\cdot (w)^{-\frac{1}{2}}_\gamma \ ,
    \ \ \ \ \ \ \ \ \ \ \ \ \ \ \
    |\balpha| \ \lesssim \ \mathcal{E}_{k,s,\gamma }\cdot \tp^{-1}
        \tm^{-s-\frac{1}{2}} \cdot (w)^{-\frac{1}{2}}_\gamma \ , \label{F_peeling1}\\
    |\rho| \ &\lesssim \ q\cdot  r^{-2} \,
         \chi_{ 1 < t < r} \ + \ \mathcal{E}_{k,s,\gamma }\cdot
         \tp^{-1-s}\tm^{-\frac{1}{2}}\cdot (w)^{-\frac{1}{2}}_\gamma  \ , \label{F_peeling2}\\
     |\sigma| \ &\lesssim \ \mathcal{E}_{k,s,\gamma }\cdot
     \tp^{-1-s}\tm^{-\frac{1}{2}}\cdot (w)^{-\frac{1}{2}}_\gamma \ , \label{F_peeling3}
\end{align}
\end{subequations}
and:
\begin{subequations}\label{phi_peeling}
\begin{align}
         |\widetilde{D}_L\phi| \ &\lesssim \
    \mathcal{E}_{k,s,\gamma }\cdot \tp^{-s -
         \frac{3}{2}}\cdot (w)^{-\frac{1}{2}}_\gamma \ ,
    &|D_{\bL}\phi| \ &\lesssim \  \mathcal{E}_{k,s,\gamma }\cdot
    \tp^{-1}
        \tm^{-s-\frac{1}{2}} \cdot (w)^{-\frac{1}{2}}_\gamma \ , \label{phi_peeling1}\\
    |\sD\phi| \ &\lesssim \ \mathcal{E}_{k,s,\gamma }\cdot
     \tp^{-1-s}\tm^{-\frac{1}{2}}\cdot (w)^{-\frac{1}{2}}_\gamma \ ,
     &|\phi| \ &\lesssim \ \mathcal{E}_{k,s,\gamma }\cdot \tp^{-1}\tm^{-s+\frac{1}{2}}\cdot
     (w)^{-\frac{1}{2}}_\gamma \
     . \label{phi_peeling2}
\end{align}
\end{subequations}
Here we have set:
\begin{align}
        |\widetilde{D}_L\phi|^2 \ &= \
        |\frac{1}{r}D_L(r\phi)|^2\, \chi_{1 < t < 2r} +
    |D_L\phi|^2 \, \chi_{r < \frac{1}{2}t} \ ,
    &|\sD\phi|^2 \ &= \  \delta^{AB} D_A\phi \overline{D_B\phi} \
        , \notag
\end{align}
and $\underline{D}$ denotes the spatial part of the connection $D$.
Also, $(\alpha,\balpha,\rho,\sigma)$ denotes the components of the
null decomposition \eqref{null_decomp} of $F_{\alpha\beta}$.
Finally, the weight functions $\tau_\pm$ and $w_\gamma$ are defined
via the formulas:
\begin{align}
        \tp^2 \ &= \ 1 + (t+r)^2 \ ,
    &\tm^2 \ &= \ 1 + (t-r)^2 \ ,
    &w_\gamma \ &= \ \tm^{2\gamma}\, \chi_{t < r} \ . \notag
\end{align}
\end{thm}\ret

\begin{rem}
In the course of our proof of Theorem \ref{main_th}, we will also
show decay (peeling) properties similar to
\eqref{F_peeling}--\eqref{phi_peeling} for the higher derivatives of
$(F,\phi)$ assuming that $2 < k$. We have not stated this formally
for the sake of brevity.
\end{rem}\ret

\begin{rem}
Note the first term in the decay asymptotic for $\rho$ on line
\eqref{F_peeling2} above. As $r\to\infty$ on any given time slice
$t=const.$ this portion of the electro-magnetic field decays only
like $r^{-2}$. This is the long range effect of the electro-static
equation \eqref{es_eq1}. What is somewhat surprising is that this
long range effect \emph{only} propagates in the exterior of the
forward cone $t=r$. That is, after rescaling  at infinity, the
effect of the charge is discontinuous across this cone. This fact
goes a long way towards explaining why it is so much easier to deal
with the case of compactly supported initial data. Also, notice that
in the case where the value \eqref{the_charge} vanishes, this long
range effect is eliminated. This, of course, is the main reason why
the assumption of zero charge is amenable to other techniques, e.g.
conformal compactification.
\end{rem}\ret

\begin{rem}
Note that we can allow any value $s_0<\frac{3}{2}$ of the weight factor in the
norms \eqref{tensor_initial_sob}--\eqref{scalar_initial_sob} applied
to the initial data so long as it is bigger than the scale invariant
factor of $\frac{1}{2}$. At this value, the norms become homogeneous
(except for the extra term of $1$ in the physical space weight) with
respect to the scaling properties of the equations
\eqref{basic_MKG}. It is also at this point that essentially all of
our estimates break down due to logarithmic divergences. Also, for
values of $s_0 < 1$, note that our condition on the initial data
\eqref{initial_smallness} requires much less decay than previous
works based on weighted energies, e.g. the classical results in
\cite{K_semi} for semilinear equations.\footnote{For instance, we
take the gradient $\underline{D}\phi$ \emph{first} before applying
the weighted $L^2$ norm.} Thus, our main theorem falls just short of
a scale invariant result, and can be considered a ``low regularity''
theorem with respect to physical space decay. Finally, we should
mention that below the level of the scaling (i.e. $s_0 <
\frac{1}{2}$), there is no reason for the type of analysis we do
here to make sense. In fact, even the simple notion of the charge
\eqref{the_charge} cannot be defined in this case because this level
of decay is consistent with the asymptotics $\phi_0\sim r^{-1 +
\delta}$ and $\dot{\phi}_0 \sim r^{-2 + \delta}$ for some $0 <
\delta$, which lead to a divergence in  \eqref{the_charge}.
\end{rem}\ret

%---------------------------------------------------------------------------

\subsection{A brief outline of the work}

We give here a very quick overview of the paper. In the next section,
we make a list of some more or less standard geometric formulas
relating to tensors and complex line bundles, as well as to the geometric
structure of Minkowski space. This section is probably best left as a
reference, and can be avoided by the reader with a passing acquaintance of
this material.\\

In the third section, we build a series of weighted $L^2$ estimates for
electro-magnetic fields. This material is for the most part a
straight-forward adaptation of the standard material found in
\cite{CK_Fields}, with some new technical devices
which are crucial for the approach we take in this paper. These are: \
A fractionally weighted version of the usual Morawetz type estimate
for null-decomposed electro-magnetic fields, the decomposition of
an electro-magnetic field
into its pure charge and charge free components, and fixed time,
characteristic, and space-time estimates which respect this
decomposition. These estimates are proved with the help of a certain
weighted elliptic estimate which we demonstrate in the appendix.\\

In the fourth section, we set about proving analogs of the estimates
of the third section for complex scalar fields. This involves an
interesting new proof of the classical Morawetz energy decay estimates
for scalar fields which is directly based on the conformal geometry
of Minkowski space. With some work, this leads to a set of estimates
which are virtually identical with what is available for
electro-magnetic fields. This device is important for what we do here
because it allows us to treat the curvature and the (conjugated
gradient) scalar field on the same footing.\\

In the fifth and sixth sections, we prove $L^\infty$ type estimates
for the electro-magnetic and complex scalar fields. Interestingly
enough, these involve in a crucial way the space-time, fixed time, as
well as characteristic energy estimates developed in the preceding
two sections.\\

In the seventh section, we recast all of the estimates we prove in
this paper in a uniform form. First we introduce an auxiliary
``tensor'' which has the behavior of the
electro-magnetic-complex-scalar field which we are studying. We then
prove a series of weighted $L^2$ bilinear estimates for interactions
of various components of this ``tensor''. It turns out that all of
the estimates we need to control the error terms of the sequel can
be put in a single form, which we call the abstract parity estimate.
This is proved at the end of this section through a simple case
analysis.\\

In the final two sections of the paper, we prove error estimates for
the commutators of the field equations \eqref{basic_MKG} with certain
Lie derivatives. With the work developed in the previous sections, this
turns out to be quite easy because everything is just a special case
of the abstract parity estimate of section seven.\\

Finally, in an appendix we prove certain Sobolev type estimates
which are used in the main work. This material is included for the
sake of completeness and can simply be referred to by the reader who
wants more detailed account of certain (standard) calculations
we use.\\ \\

\noindent
{\em{Acknowledgments:}} The authors would like to thank Sergiu
Klainerman and Igor Rodnianski who played a large role in initiating
and supporting this project. This investigation originated from a conversation
with those two as to whether one could extend the low regularity
approach of \cite{KM_MKG} to discuss the scattering behavior of the
CSF system \eqref{basic_MKG}. It was S. Klainerman's insistence that one should
first ``do it classically'' in order to get a feel for the bulk behavior
of the field quantities, and to see if anything interesting turns up, that
lead to our study.
The possibility that one could use space-time estimates
to eliminate logarithmic divergences coming from
the charge was first pointed out to us by I. Rodnianski.
This plays a crucial role in our work. We would also
like to thank him for mentioning the connection with conformal geometry
that a previous proof of the divergence identities
\eqref{conf_mor_div_iden1}--\eqref{conf_mor_div_iden2} had. We have used
this elegant description in Section \ref{phi_L2_sect}.

Part of this work was done while we were Members of the
Institute for  Advanced Study, Princeton. H.L. was supported by the NSF
grant DMS-0111298 to the Institute and partially supported by the 
the NSF Grant DMS-0200226. 
J.S. was supported by a Veblin Research Instructorship and an NSF postoctoral fellowship.

\ret
%-------------------------------------------------------------------------
%%%%%%%%%%%%%%%%%%%%%%%%%%%%%%%%%%%%%%%%%%%%%%%%%%%%%%%%%%%%%%%%%%%%%%%%%%
%-------------------------------------------------------------------------

\section{Some Geometric Preliminaries}

In this section, we introduce the basic geometric concepts and
identities which we will use in our analysis of the field equations
\eqref{basic_MKG}. All of this material is completely standard, and
we review it here in detail solely for the convenience of the
reader. Experienced readers will find it worthwhile to skip this
section and simply refer to it at places in the sequel where a
specific calculation requires one of the identities we list here.\\

\subsection{Lie Derivatives}

We begin with some well known formulas involving Lie derivatives and
divergences of vector-fields and differential forms. For any
sections $X,Y$ to $T\mathcal{M}$ we have the usual formula
$\mathcal{L}_X Y = [X,Y]$. We also denote by $\mathcal{L}_X \omega$
and $\mathcal{L}_X F$ the Lie derivative of a one-form and two-form
respectively. In the frame $\{e_\alpha\}$ these can be computed as:
\begin{subequations}\label{otf_lie_der}
\begin{align}
        (\mathcal{L}_X \omega)_\alpha \ &= \ X(\omega_\alpha) -
        \omega([X,e_\alpha]) \ , \label{otf_lie_der1}\\
    (\mathcal{L}_X F)_{\alpha\beta} \ &= \
    X(F_{\alpha\beta}) - F([X,e_\alpha], e_\beta) - F(e_\alpha,
    [X,e_\beta]) \ . \label{otf_lie_der2}
\end{align}
\end{subequations}
For vector-fields $X$, we form the Lorentzian divergence:
\begin{equation}
        \Div X \ = \ \nabla_\alpha X^\alpha \ = \
    \tr \nabla X \ . \label{div_def}
\end{equation}
Also, for each vector-field $X$ on $\mathcal{M}$, we measure the
effect of its flow on the Minkowski metric by forming its Lorentzian
deformation tensor:
\begin{equation}
        \mathcal{L}_{X} g \ = \ {}^{(X)}\pi \ .
        \label{lie_def_tensor}
\end{equation}
By Lie differentiating the contraction $g_{\alpha\beta} =
g^{\gamma\delta} g_{\gamma\alpha} g_{\delta_\beta}$ we have that:
\begin{equation}
    (\mathcal{L}_X g^\dagger)^{\alpha\beta} \ = \ -\,
    {}^{(X)}\pi^{\alpha\beta} \ , \label{dagger_form}
\end{equation}
where we have set $g^\dagger = g^{\alpha\beta}$. Using this formula,
one can easily compute the commutator of the Lie derivative
$\mathcal{L}_X$ and the duality operator ${}^*$ applied to a
two-form:
\begin{align}
        [\mathcal{L}_X , {}^*\ ]\,  F_{\alpha\beta} \ &= \
    \frac{1}{2}\, \mathcal{L}_X (\in_{\alpha\beta\gamma\delta}g^{\gamma\mu}g^{\delta\nu}
    )F_{\mu\nu} \ , \notag\\
    &= \ -\in_{\alpha\beta\gamma}^{\ \ \ \ \ \delta}
    {}^{(X)}\pi^{\gamma\sigma} F_{\sigma\delta} \ + \
    \big( \Div X \big) \, {}^* F_{\alpha\beta}
    \ . \label{star_lie_comm}
\end{align}
Now, a well known calculation shows that we have:
\begin{equation}
       {}^{(X)}\pi \ = \ 2\sym \nabla X \ , \label{def_ten_id}
\end{equation}
where $\sym$ denotes the symmetric part of the tensor
$\nabla X$. Since the trace of an antisymmetric tensor is
always zero, using \eqref{div_def} this gives the following
formula:
\begin{equation}
        \Div X \ = \ \frac{1}{2} \tr\, {}^{(X)}\pi \ .
        \label{div_trace_form}
\end{equation}
A direct use of this gives the calculation:
\begin{align}
        X\big( \Div Y\big) - Y\big( \Div X\big) \ &= \
    \frac{1}{2} \mathcal{L}_X \big( g^{\alpha\beta}
    (\mathcal{L}_{Y}g)_{\alpha\beta}\big) -
    \frac{1}{2} \mathcal{L}_Y \big( g^{\alpha\beta}
    (\mathcal{L}_{X}g)_{\alpha\beta}\big) \ , \notag \\
    &= \ \frac{1}{2} g^{\alpha\beta} (\mathcal{L}_{ [ X,Y] }g)_{\alpha\beta}
    \ , \notag \\
    &= \ \Div ( \mathcal{L}_X Y ) \ . \label{div_comm_calc}
\end{align}
In particular, for any vector-field  $X$ such that $\Div X =
const.$, we have the following  useful formula:
\begin{equation}
        X\big( \Div Y\big) \ = \ \Div(\mathcal{L}_{X} Y) \ . \notag
\end{equation}\ret
This discussion goes through almost verbatim for one-forms $\omega$.
In this case we define:
\begin{equation}
        \Div \omega \ = \ \nabla_\alpha \omega^\alpha \ . \notag
\end{equation}
Then using the identification $\omega^\alpha =
g^{\alpha\beta}\omega_\beta$, we have from \eqref{div_comm_calc} and
\eqref{dagger_form} that in any frame $\{e_\alpha\}$ that:
\begin{equation}
        X\big( \Div \omega\big) - \omega^\alpha
    e_\alpha \big( \Div X\big) \ = \ -\, \nabla_\alpha \,
    {}^{(X)}\pi^{\alpha\beta} \omega_\beta +
    \Div (\mathcal{L}_X \omega) \ , \label{div_omega_calc}
\end{equation}
Thus, again assuming  that $\Div X = const.$ one can drop the second
term on the left hand side of \eqref{div_omega_calc} above.\\

To extend the formulas \eqref{div_comm_calc} and
\eqref{div_omega_calc} to higher order tensors, it is useful to
define the Lie derivative of the Levi-Civita connection $\nabla$ via
the formula:
\begin{equation}
        \mathcal{L}_X(\nabla) \ = \ \mathcal{L}_X\nabla - \nabla
        \mathcal{L}_X \ . \label{nabla_lie}
\end{equation}
It is easy to see that \eqref{nabla_lie} in fact defines the tensor,
which in our zero curvature setting is given by:
\begin{equation}
        \mathcal{L}_X(\nabla) \ = \ \nabla\nabla X \ . \notag
\end{equation}
This object acts in a multilinear fashion on sections to
$T^k_l\mathcal{M}$. For example, with respect to
a contravariant two-tensor the action reads:
\begin{equation}
        \left(\mathcal{L}_X(\nabla) T\right)_{\alpha\beta\gamma}
        \ = \ (\nabla_\alpha\nabla_\beta X^\delta) T_{\delta\gamma}
        + (\nabla_\alpha\nabla_\gamma X^\delta) T_{\beta\delta} \ .
        \notag
\end{equation}
The above formula immediately applies to two-forms. In this case, if
${}^{(X)}\pi$ is a constant multiple of the metric
$g_{\alpha\beta}$, which also implies that $\Div X = const.$ by
formula \eqref{div_trace_form}, then antisymmetry and some simple
manipulations yields the following important commutator formula:
\begin{equation}
        \mathcal{L}_X \nabla^\beta F_{\alpha\beta} \ = \
        -\, {}^{(X)}\pi^{\beta\gamma}\nabla_\gamma F_{\alpha\beta}
        +  \nabla^\beta (\mathcal{L}_X F)_{\alpha\beta} \ . \label{F_div_iden}
\end{equation}
This will be used in to differentiate the Maxwell equation
\eqref{basic_Maxwell} with respect to various vector-fields which
are Killing and conformal Killing, and will be the basis of our
analysis for that portion of the system \eqref{basic_MKG}.\\

We would now like to set up analogs of some of the above formulas
for vector-fields and one-forms which are naturally associated with
the complex line bundle $V$. We consider complexified versions of
the tangent and cotangent bundles:
\begin{align}
        T^\mathbb{C}\mathcal{M} \ &= \ V\otimes T\mathcal{M} \ ,
        &{T^*}^\mathbb{C}\mathcal{M} \ &= \ V\otimes T^*\mathcal{M} \
        . \label{complex_bundles}
\end{align}
In local coordinates, sections to these bundles can simply be
identified with complex valued vector-fields and one-forms
respectively. Of special significance is the full covariant
derivative of a section to $V$, which we denote by $D\phi$. In
accordance with the above notation, this can also be seen as the
complex exterior derivative of the scalar $\phi$ which is defined
via the natural formula:
\begin{equation}
        D\phi(X) \ = \ D_X\phi \ . \notag
\end{equation}
This should be understood in analogy with the natural formula for
the exterior derivative of a real valued scalar $f$ on $\mathcal{M}$
which is given by $df(X) = X(f)$. For real scalars, one has the Lie
derivative formula: \ $\mathcal{L}_X df = d X(f)$, and we seek an
analog of this relation in the complex setting. Therefore, we define
a covariant Lie derivative, called $\mathcal{L}^\mathbb{C}$ on the
bundles \eqref{complex_bundles} via the following formulas which
hold for arbitrary tensors (in analogy with \eqref{cov_der_A}):
\begin{equation}
        \LC_X (T) \ = \ \mathcal{L}_X T + \sqrt{-1}  \
        A(X) T \ ,  \ \label{implicit_cdef}
\end{equation}
where $\mathcal{L}_X T$ is somewhat abusive notation in that we
consider $T$ as a complex valued tensor, whose coefficients are
computed in the trivialization of $\{A_\alpha\}$. It should be noted
that by the Leibnitz rule, the formula \eqref{implicit_cdef} is
completely gauge independent. In the sequel, we will not make
explicit use of \eqref{implicit_cdef} which is instead provided as a
concrete way of visualizing $\LC$. We note here that alternatively,
one could define $\LC$ abstractly via the Leibnitz rule which would be
sufficient for our purposes.\\

Now, using the commutator formula \eqref{curvature_def}, it is seen
that the complex Lie derivative satisfies the following natural
formula with respect to $D$:
\begin{equation}
        \LC_X D\phi \ = \ D (D_X\phi) \ + \ \sqrt{-1}\,  i_XF \cdot \phi \ .
        \label{D_complex_lie}
\end{equation}
Here $(i_X F)_\alpha \ = \ X^\gamma F_{\gamma\alpha}$ is the
interior product of $X$ and $F$.\\

Next, using formula \eqref{implicit_cdef} we have the following
useful expression involving complex scalars $\phi$ and complex
valued one-forms $\eta$:
\begin{equation}
        \mathcal{L}_X \ \Im \left(\phi \cdot \overline{\eta}\right)
    \ = \ \Im\left(D_X\phi\cdot \overline{\eta}\right)
    + \Im\left(\phi\cdot \overline{\LC_X \eta} \right)
    \ . \label{basic_J_Lie}
\end{equation}
Furthermore, notice that from \eqref{implicit_cdef} the operation of
complex Lie differentiation is well behaved with respect to
contractions. For example, for a real valued two-from and real and
complex vector-fields $X,Y$ and $Z$ (resp.) we have the formula:
\begin{equation}
        D_Y F(X,Z) \ = \ \mathcal{L}_Y F (X, Z) + F( [Y,X] , Z ) +
        F(X, \LC_Y Z) \ . \label{complex_lie_cont_id}
\end{equation}
We will use this formula \emph{instead} of the corresponding formula
involving the covariant derivative $\nabla_Y$. This is necessary
because when we differentiate $F$ it will always be with respect to
$\mathcal{L}$ and \emph{not} $\nabla$. \\

Finally, using a more abstract form of the formula
\eqref{implicit_cdef}, we compute the covariant divergence of a
complex valued one-form via the implicit formula:
\begin{equation}
        \cdiv (\phi\otimes \omega) \ = \ \omega^\alpha D_\alpha\phi
        + \Div \omega\cdot \phi \ , \notag
\end{equation}
for any complex scalar field $\phi$ and real one-form $\omega$. Using this
last line in conjunction with the identity \eqref{div_omega_calc}
above, we have the following commutator formula for the complex
divergence and complex Lie derivative for complex valued one-forms
$\eta$:
\begin{equation}
        D_X (\cdiv \eta) - \eta^\alpha e_\alpha (\Div X) \ = \
        -\, D_\alpha {}^{(X)}\pi^{\alpha\beta} \eta_\beta + \cdiv
        (\LC_X\eta) + \sqrt{-1} X^\alpha F_{\alpha\beta} \eta^\beta \
        . \notag
\end{equation}
In particular, if ${}^{(X)}\pi$ is a constant multiple of the metric
$g_{\alpha\beta}$, we have the following formula for the commutator
of the covariant wave equation and the covariant derivative (or Lie
derivative if you like) $D_X$:
\begin{align}
        &D^\alpha D_\alpha D_X \phi \ , \label{cbox_comm_formula}\\
        = \ &\cdiv D (D_X\phi) \ ,
        \notag\\
        = \ &\cdiv \left( \LC_X D\phi - \sqrt{-1}\, i_X F \cdot\phi
        \right) \ , \notag\\
        = \ &D_X D^\alpha D_\alpha \phi +
        {}^{(X)}\pi^{\alpha\beta}D_\alpha D_\beta \phi +
        \sqrt{-1}\left( 2 X^\alpha F_{\alpha\beta}D^\beta\phi -
        \nabla^\alpha (X^\beta F_{\alpha\beta})\cdot \phi\right) \ .
        \notag
\end{align}
The important thing to notice here is that the above formula depends
only on the curvature $F$ and not the choice of gauge potentials
$\{A_\alpha\}$.\ret

%-------------------------------------------------------------------------

\subsection{The Null Frame, Lorentz Algebra, and Associated
Covariant Derivatives}

Of central importance to the analysis we do here is the freedom to
perform calculations on the equations \eqref{basic_MKG} in an
arbitrary frame. In particular, the ability to use frames which do
not arise from any system of coordinates. Since we are assuming that
our initial data is both smooth and well localized around the origin
in $\mathbb{R}^3$, there is a canonical choice of frame for our
problem. This is the so called \emph{standard spherical null frame}.
The first two members of this are the null generators of forward and
backward light cones which we define respectively as:
\begin{align}
        L \ &= \ \partial_t + \partial_r \ , &\bL \ &= \
    \partial_t - \partial_r
    \ . \label{L_bL}
\end{align}
To complete things, we need only define derivatives in the angular directions.
This can be done in an identical fashion on each time slice $\{t = const.\}$
so we only need to define things on $\mathbb{R}^3$. If we let
$\{e^0_A\}_{A=1,2}$ denote a local  orthonormal frame for the unit sphere
in $\mathbb{R}^3$, then for each value of the radial variable
we can by extension define:
\begin{equation}
        e_A \ = \ \frac{1}{r}\, e^0_A \ . \label{ang_frame}
\end{equation}
Thus the collection $\{e_A\}_{A=1,2}$ forms an orthonormal basis on each
sphere $\{r=const.\}$ on each fixed time slice. We can write the usual
translation invariant frame $\{\partial_i\}$
in terms of $\partial_r$ and this basis as follows:
\begin{equation}
        \partial_i \ = \ \omega_i \partial_r + \omega^A_i e_A \ ,
    \notag
\end{equation}
where $\omega^A_i = e_A(x^i)$, which follows at once from the formula:
\begin{equation}
        \omega^A_i \ = \ \langle \partial_i , e_A \rangle
    \ = \ dx^i(e_A) \ = \ e_A(x^i) \ . \label{hans_comp}
\end{equation}
In line above we have used in order orthonormality, duality,  and then
simply the definition of exterior derivative.
Also, we note here that the $\{\omega_i^A\}$ are part of a rotation
matrix at each point, which implies the following useful formulas:
\begin{align}
        \omega^i \omega_i^A \ &= \ 0 \ ,
    &\omega^i_A\omega_i^{B} \ &= \ \delta_A^{\ B} \ ,
    &\omega_i^A  \omega_A^j \ &= \ \delta_i^{\ j} - \omega_i\omega^j
    \ . \label{omega_relations}
\end{align}\ret

With respect to the full frame $\{\bL,L,e_A\}$,
the Minkowski metric reads:
\begin{equation}
        g \ = \ -2 \theta^{\bL}\otimes \theta^{L} -2
    \theta^{L}\otimes \theta^{\bL} + \delta_{AB} \theta^A\otimes\theta^B
    \ , \notag
\end{equation}
where $\{\theta^{\bL},\theta^{L},\theta^A  \}$ is the
corresponding dual frame. Associated with this null frame are
the standard optical functions:
\begin{align}
        \bu \ &= \ t+r \ , &u \ &= \  t-r \ . \label{optical_funs}
\end{align}
Note that one has the identities:
\begin{align}
        L (\bu) \ &= \ 2 \ , &\bL(\bu) \ &= \ 0 \ , \notag\\
    L (u) \ &= \ 0 \ , &\bL(u) \ &= \ 2 \ . \notag
\end{align}
In particular both $u$ and $\bu$ solve the eikonal equation
$\nabla^\alpha h \nabla_\alpha h = 0$. In terms of the frame
\eqref{L_bL}--\eqref{ang_frame} we have the standard null
decomposition of the electro-magnetic field tensor $F$, as well
as the gradient $D\phi$. As usual we define:
\begin{subequations}\label{null_decomp}
\begin{align}
        \alpha_A \ &= \ F_{LA} \ , &\balpha_A \ &= \ F_{\bL A}\ ,
        \label{null_decomp1}\\
    \rho \ &= \ \frac{1}{2}F_{\bL L} \ ,
    &\sigma \ &= \ \frac{1}{2} \in^{AB} F_{AB}
    \ . \label{null_decomp2}
\end{align}
\end{subequations}
Using the duality operator ${}^*$, we have the useful formulas for
the null decomposition of ${}^*F$, where $\in_{AB}=\frac{1}{2}
\in_{\bL L AB}$ is the volume element
on the spheres $r=const.$:
\begin{subequations}\label{dual_null_decomp}
\begin{align}
        {}^*\alpha_A \ &= \ - \in_A^{\ B} \alpha_B \ ,
    &{}^*\balpha_A \ &= \ \in_A^{\ B} \balpha_B \ ,
    \label{dual_null_decomp1}\\
    {}^* \rho \ &= \ \sigma \ ,
    &{}^*\sigma \ &= \ -\rho \ . \label{dual_null_decomp2}
\end{align}
\end{subequations}
One can expand out the energy-momentum tensors \eqref{phi_em} and
\eqref{F_em} in the null directions \eqref{L_bL}. Using the
identifications \eqref{null_decomp}, for the electro-magnetic field
these read:
\begin{subequations}\label{F_em_decomp}
\begin{align}
        Q(L,L)[F] \ &= \  \delta^{AB} \alpha_A \alpha_B \ = \ |\alpha|^2
    \ , \label{F_em_decomp1}\\
    Q(\bL,\bL)[F] \ &= \  \delta^{AB} \balpha_A \balpha_B \ = \ |\balpha|^2
    \ , \label{F_em_decomp2}\\
    Q(\bL,L)[F] \ &= \  \rho^2 + \sigma^2 \ . \label{F_em_decomp3}
\end{align}
\end{subequations}
Furthermore, for the energy-momentum tensor of the scalar field
\eqref{phi_em}, one has that:
\begin{subequations}\label{phi_em_decomp}
\begin{align}
        Q(L,L)[\phi] \ &= \ |D_L \phi|^2 \ , \label{phi_em_decomp1}\\
    Q(\bL,\bL)[\phi] \ &= \ |D_{\bL} \phi|^2 \ , \label{phi_em_decomp2}\\
    Q(\bL,L)[\phi] \ &= \ \delta^{AB} D_A\phi \overline{D_B\phi} \
    = \ |\sD\phi|^2 \ . \label{phi_em_decomp3}
\end{align}
\end{subequations}\\

Next, we record here the following standard formulas for the Levi-Civita
connection $\nabla$ with respect to this frame. This will be
used many times in the  sequel:
\begin{subequations}\label{cov_ders}
\begin{align}
        \nabla_{L} L \ &= \ 0 \ , &\nabla_{\bL} L \ &= \ 0 \ ,\label{cd1}\\
    \nabla_{\bL} \bL \ &= \ 0 \ , &\nabla_{L} \bL \ &= \ 0
    \ ,\label{cd2}\\
    \nabla_{L} e_A \ &= \ 0 \ ,  &\nabla_{\bL} e_A \ &= \ 0
    \ , \label{cd3}\\
    \nabla_{e_A}L \ &= \ \frac{1}{r} e_A \ ,
    &\nabla_{e_A}\bL \ &= \ - \frac{1}{r} e_A \ , \label{cd4}\\
    \nabla_{e_A} e_B \ &= \ \overline{\nabla}_{e_A} e_B
    + \frac{1}{2r}\delta_{AB}(\bL - L) \ , \label{cd5}\\
    &= \ \overline{\Gamma}^D_{AB} e_D +
    \frac{1}{2r}\delta_{AB}(\bL - L) \ . \label{cd6}
\end{align}
\end{subequations}
Here, $\overline{\nabla}$ denotes the intrinsic covariant
differentiation on spheres, and $\overline{\Gamma}$ its associated
frame-Christoffel symbols. The key property of the $\overline{\Gamma}$
we will use in the sequel is that they are homogeneous functions
of degree $-1$ with respect to the radial variable $r$.\\

Our result is based in part on the usual process of obtaining weighted
energy inequalities, which in turn give $L^\infty$ estimates of the
type \eqref{F_peeling}--\eqref{phi_peeling}.
The weights in these kind of estimates are closely
related to the Killing and conformal-Killing structure of Minkowski
space, as was first realized to full effect in the seminal work
\cite{K_semi}. Accordingly, we introduce the inhomogeneous Lorentz
algebra:
\begin{equation}
        \mathbb{L} \ = \ \{\partial_\alpha , \Omega_{\alpha\beta} , S \}
    \ , \label{inh_lor_alg}
\end{equation}
where we use $x^0 = t$ and $x_0 = -t$ to express:
\begin{align}
        \Omega_{\alpha\beta} \ &=  \ x_\alpha\partial_\beta - x_\beta
    \partial_\alpha \ , \label{lor_rot} \\
    S \ &= \ x^\alpha\partial_\alpha \ . \label{scaling}
\end{align}
This collection of vector-fields satisfy the well known
algebraic relations:
\begin{subequations}\label{lor_bracket_rel}
\begin{align}
    [ \partial_\alpha , \Omega_{\beta \gamma}] \ &= \
    \delta_{\alpha[\beta} \partial_{\gamma]} \ , \label{lor_bracket_rel1}\\
    [ \Omega_{\alpha\beta}, \Omega_{\gamma\sigma}] \ &= \
    - \delta_{(\alpha\gamma}\Omega_{\beta\sigma)}
    \ , \label{lor_bracket_rel2}\\
    [\partial_\alpha , S] \ &= \ \partial_\alpha \ , \label{lor_bracket_rel3}\\
    [\Omega_{\alpha\beta} , S] \ &=  0 \ . \label{lor_bracket_rel4}
\end{align}
\end{subequations}
Expanding the formulas \eqref{lor_rot}--\eqref{scaling} out in the
frame \eqref{L_bL}--\eqref{ang_frame}, and using the symbol
$\omega^A_i = \theta^A(\partial_i)$ we have the null decompositions:
\begin{subequations}\label{lor_null_decomp}
\begin{align}
    \partial_0 \ &= \ \frac{1}{2}(L + \bL) \ , \label{lor_null_decomp1}\\
    \partial_i \ &= \  \frac{\omega_i}{2}( L - \bL)
    + \omega^A_i e_A \ , \label{lor_null_decomp2}\\
        \Omega_{ij} \ &= \ \Omega_{ij}^A e_A \ = \
    \big( x_i \omega^A_j -
    x_j \omega^A_i\big)e_A \ , \label{lor_null_decomp3}\\
    \Omega_{i0} \ &= \ \frac{\omega_i}{2}( \bu L - u \bL)
    + t \omega^A_i e_A \ , \label{lor_null_decomp4}\\
    S \ &= \ \frac{1}{2}\left(\bu L + u \bL\right)
    \ . \label{lor_null_decomp5}
\end{align}
\end{subequations}
It will also be necessary for us to have a precise account of
the algebraic relations between the frame
\eqref{L_bL}--\eqref{ang_frame} and the fields \eqref{inh_lor_alg}.
These can easily be computed using the formulas
\eqref{lor_null_decomp}. We record this computation here as:
\begin{subequations}\label{lor_null_brack}
\begin{align}
    [L , \partial_0 ] \ &= \ [\bL , \partial_0 ] \ = \
    [e_A , \partial_0 ] \ = \ 0 \ , \label{lor_null_brack1}\\
    [ L ,\partial_i ] \ &= \ \omega^A_i [ L , e_A] \ = \
    - \frac{\omega^A_i}{r} e_A \ , \label{lor_null_brack2}\\
    [ \bL ,\partial_i ] \ &= \ \omega^A_i [ \bL , e_A] \ = \
    \frac{\omega^A_i}{r} e_A \ , \label{lor_null_brack3}\\
    [L , \Omega_{ij} ] \ &= \ -[\bL,\Omega_{ij}] \ = \
    0 \ , \label{lor_null_brack4}\\
    [e_B , \Omega_{ij}] \ &= \ e_B(\Omega_{ij}^A) e_A
    + \Omega_{ij}^A \textbf{c}_{BA}^D e_D \ , \label{lor_null_brack5}\\
    [L , \Omega_{i0}] \ &= \ \omega_i L - \frac{u}{r}\omega^A_i
    e_A \ , \label{lor_null_brack6}\\
    [\bL,\Omega_{i0}] \ &= \ - \omega_i \bL + \frac{\bu}{r}
    \omega^A_i e_A \ , \label{lor_null_brack7}\\
    [e_B , \Omega_{i0}] \ &= \ \frac{\omega_{Bi}}{2r}(
    \bu L - u \bL) -  t \omega_i^A \overline{\Gamma}_{AB}^D
    e_D \ , \label{lor_null_brack8}\\
    [L, S] \ &= \ L \ , \label{lor_null_brack9}\\
    [\bL,S]\ &= \ \bL \ , \label{lor_null_brack10}\\
    [e_A , S ] \ &= \ e_A . \label{lor_null_brack11}
\end{align}
\end{subequations}
In formula \eqref{lor_null_brack5} above we have set
$\textbf{c}_{AB}^D = [e_A,e_B]^D$, which is again a smooth function
homogeneous of degree $-1$ in the radial variable $r$. Also, on line
\eqref{lor_null_brack8} we have used the identities:
\begin{align}
       e_B(\omega_i^A) \ &= \
       \nabla_{e_B}\theta^A(\partial_i) \ =\
       -\delta_B^A \frac{\omega_i}{r} -
       \overline{\Gamma}_{BD}^A\omega_i^D \ ,
       &\textbf{c}_{BA}^D \ &= \ \overline{\Gamma}_{BA}^D
       - \overline{\Gamma}_{AB}^D \ . \notag
\end{align}\ret

Finally, to end this subsection, we list here the various formula
for the action of the Levi-Civita connection $\nabla$ on the algebra
\eqref{inh_lor_alg}. All of these formulas are computed in a
straightforward way using the identities \eqref{cov_ders}
and \eqref{lor_null_decomp}:
\begin{subequations}\label{cov_ders_vect_fields}
\begin{align}
        \nabla_L \partial_\alpha  \ &= \
    \nabla_{\bL} \partial_\alpha \ = \ \nabla_{e_A}\partial_\alpha
    \ = \ 0 \ , \label{cov_ders_vect_fields1}\\
    \nabla_L \Omega_{ij} \ &= \ -\nabla_{\bL}\Omega_{ij} \ = \
    \frac{1}{r}\Omega_{ij} \ , \label{cov_ders_vect_fields2}\\
    \nabla_{e_B} \Omega_{ij} \ &= \ e_B(\Omega_{ij}^A)e_A
    + \Omega^A_{ij} \overline{\Gamma}_{BA}^D e_D +
    \frac{1}{2r}\Omega^B_{ij} (\bL - L) \ , \label{cov_ders_vect_fields3}\\
    \nabla_{L}\Omega_{i0} \ &= \  \omega_i L + \omega^A_i e_A
    \ , \label{cov_ders_vect_fields4}\\
    \nabla_{\bL}\Omega_{i0} \ &= \  - \omega_i \bL + \omega^A_i e_A
    \ , \label{cov_ders_vect_fields5}\\
    \nabla_{e_B} \Omega_{i0} \ &= \ \frac{\omega_{Bi}}{2}(L + \bL)
    \ , \label{cov_ders_vect_fields6}\\
    \nabla_L S \ &= \ L  \ , \label{cov_ders_vect_fields7}\\
    \nabla_{\bL}S \ &= \ \bL \ , \label{cov_ders_vect_fields8}\\
    \nabla_{e_A} S\ &= \ e_A \ . \label{cov_ders_vect_fields9}
\end{align}
\end{subequations}\ret

\ret
%-------------------------------------------------------------------------
%%%%%%%%%%%%%%%%%%%%%%%%%%%%%%%%%%%%%%%%%%%%%%%%%%%%%%%%%%%%%%%%%%%%%%%%%%
%-------------------------------------------------------------------------

\section{Fixed Time and Space-time Energy Estimates for the
Curvature}\label{F_L2_section}

In this section, we begin to build the estimates which lie at the
center of our approach. All of these will be produced in some way
through the tensorial conservation laws \eqref{F_consv} and
\eqref{phi_consv}. In this section we will deal with weighted $L^2$
type estimates which involve the curvature tensor $F_{\alpha\beta}$.
We assume that this satisfies the Maxwell equation
\eqref{basic_Maxwell} for an unspecified current density
$J_\alpha$.\\

On Minkowski space, the two most basic energy estimates (really
identities) for $F_{\alpha\beta}$ are based respectively on time
translation invariance and the conformal structure of the Minkowski
metric. This is utilized by contracting the energy-momentum tensor
\eqref{F_em} with the following two vector fields (respectively):
\begin{align}
        T \ &= \ \partial_0 \ = \ \frac{1}{2}(L+ \bL)\ , \label{dt_field}\\
        K_0 \ &= \ (t^2+|x|^2)\partial_0 + 2t x^i\partial_i \ = \
        \frac{1}{2}( \bu^2 L + u^2\bL ) \ . \label{mor_field}
\end{align}
One readily computes their deformation tensors to be:
\begin{align}
        {}^{(T)}\pi \ &= \ 0 \ , &{}^{(K_0)}\pi \ = \ 4t \, g \
        . \notag
\end{align}
Due to the trace-free nature of the
energy-momentum tensor \eqref{F_em}, the one-form resulting from
such a contraction is seen to be divergence free. Therefore,
defining the weights:
\begin{align}
        \tp^2 \ &= \ 1 + \bu^2 \ , &\tm^2 \ &= \ 1 + u^2
\end{align}
and the hybrid vector-field:
\begin{equation}
        \overline{K}_0 \ = \ T + K_0 \ , \notag
\end{equation}
we arrive at the following fundamental estimate through the
application of the geometric Stokes theorem\footnote{We shall
discuss this type of procedure in more detail in the sequel. See
also \cite{CK_Fields} for a more thorough discussion and a
derivation of this particular estimate.} over domains of the type
$\{u\leqslant u_0\}\cap\{0\leqslant t \leqslant t_0\}$ and expanding
out the resulting quantities using the identities
\eqref{F_em_decomp}:
\begin{multline}\label{basic_mor_maxwell}
        \int_{\RR \cap \{t=t_0\}} \ \tp^2 |\alpha|^2 + \tm^2|\balpha|^2
    + (\tp^2+ \tm^2)(\rho^2 + \sigma^2)\ dx \ \\
    + \ \sup_{u}\ \int_{C(u) \cap\{0\leqslant t \leqslant t_0\}}
    \ \tp^2 |\alpha|^2 + \tm^2(\rho^2 + \sigma^2)\
    dV_{C(u)} \\
    \lesssim \ \int_0^{t_0}\int_{\RR}\ |(\overline{K}_0)^\alpha
    F_{\alpha\beta} J^\beta|\ dxdt \ + \
    \int_{\RR \cap \{t=0\}} \ (1+r^2)( |E|^2 + |H|^2 )\  dx
    \ .
\end{multline}
Here $dV_{C(u)}$ is the Euclidean volume element on the cone
$\{t-r = u\}$.\\

As it stands, estimate \eqref{basic_mor_maxwell} is not terribly
useful except in certain restricted situations.  This is because in
general it is not possible to assume that the right hand side of
\eqref{basic_mor_maxwell} above is finite, even if one assumes that
the ``dynamic'' portion of the initial data is compactly supported.
This is the effect of the constraint equation on line \eqref{EH_Maxwell1}.
Indeed, taking a Hodge decomposition $E = E^{df} + E^{cf}$ of $E$ into its
divergence-free and curl-free components (respectively), and
introducing the potential function:
\begin{equation}
        E^{cf} \ = \ \nabla \varphi \ , \notag
\end{equation}
we see that this constraint is equivalent to the elliptic equation:
\begin{equation}
        \Delta \varphi \ = \ J_0 \ . \label{charge_eq}
\end{equation}
In particular, it is not possible for $E$ to decay better that $E\sim
\frac{1}{r^2}$ unless the following quantity known as the
\emph{charge} vanishes:
\begin{equation}
        q(F)(t) \ = \ \int_{\RR} J_0(t) \ dx \ . \label{charge_def}
\end{equation}
Notice that this quantity is a constant of motion thanks to the
continuity equation \eqref{J_cont}. Therefore, we shall henceforth
refer to it as $q(F)$ and drop the dependence on $t$. Because of the
weights involved, there is no hope of obtaining a
finite value for the energy \eqref{basic_mor_maxwell} without
the extra assumption that $q(F)=0$.\\

In order to circumvent this problem it is necessary to either modify
the energy \eqref{basic_mor_maxwell} so that less decay of the
initial data is required, or to modify the field equations
\eqref{basic_Maxwell} themselves so that the behavior resulting from
the charge is eliminated. A naive approach to the first tactic would
be to simply decrease the amount of decay required on the right hand
side of \eqref{basic_mor_maxwell}. For example, this can be done by
eliminating the use of the conformal field \eqref{mor_field}. As we
shall see in a moment, it does not seem possible to do this and
still retain the distinct weights (peeling) on the different
components on the left hand side of \eqref{basic_mor_maxwell}.
Furthermore, recall that the simple scaling argument used in the
introduction implies that it should not be possible to close a
global existence proof for the system \eqref{basic_MKG} under the
assumption of a total decay rate for the curvature less than $F\sim
\frac{1}{r^2}$ at space-like infinity. Therefore, we use
a different approach here.\\

A more sophisticated approach to the first tactic, which has been
used by other authors in similar contexts \cite{CK_Ein}, would be
to employ some Lie derivatives to the field strength
$F_{\alpha\beta}$ before putting it in weighted $L^2$ in order to
kill off the effect of the charge. For instance, the Lie derivatives
with respect to the rotation fields $\Omega_{ij}$ work well to
accomplish this because the leading order behavior of the charge is
spherically symmetric. However, this strategy runs into technical
complications when it comes to dealing with the boost vector fields
$\Omega_{i0}$, which we will make essential use of in our proof of
the sharp weighted $L^2(L^\infty)$. Therefore, we shall follow a tactic which
is both conceptually and technically much simpler. This is as follows:\
Solving the charge equation \eqref{charge_eq} explicitly,
we see that:
\begin{equation}
        \varphi(t,x) \ = \ -\ \frac{1}{4\pi}\ \int_{\RR} \frac{1}{|x-y|} J_0(t,y)\ dy \
        . \notag
\end{equation}
Therefore, at fixed time, we have the following asymptotic
behavior for $E^{cf}$ as $r\to\infty$:
\begin{equation}
        E^{cf}_i \ \sim  \  q\cdot \frac{\omega_i}{4\pi r^2} \ .
        \label{charge_asym}
\end{equation}
Since this behavior is sufficiently simple we shall subtract it off,
proving estimates for the remaining field strength. However, to do
this correctly it is necessary to take into account the behavior of
the charge in the wave zone $t\sim r$. Although it is not obvious
at first, there is a jump type behavior of the charged component of
$E$ across the cone $u=0$ as long as the current vector $J_\alpha$
satisfies certain weighted estimates which are consistent with the
right hand side of \eqref{basic_mor_maxwell} above. Keeping
\eqref{charge_asym} in mind, the approximation we make to the
charged component $E^{cf}_i$ is the following. We first define the
\emph{charge two-form}, denoted by $\overline{F}_{\alpha\beta}$, via
the exterior derivative (for $0\leqslant t$):
\begin{equation}
        \overline{F} \ = \ -\ \frac{q(F)}{4\pi}\ d\left(
        [ \int_0^r \frac{1}{s^2} \chi^+(s-t-2)\ ds  ]\,
        dt\right) \ . \label{bar_F_def}
\end{equation}
In the above formula, $\chi^+$ denotes a smoothed out version of the
Heaviside function; i.e. a non-decreasing $C^\infty$ function equal
to $0$ on $(-\infty , 0]$ and $1$ on $[1,\infty)$. Notice that by
design, the tensor $\overline{F}_{\alpha\beta}$ satisfies the
Bianchi identity \eqref{basic_Maxwell2} because it is exact. Also,
one has the explicit formulas for the electric-magnetic
decomposition \eqref{em_decomp} of $\overline{F}$:
\begin{align}
    \overline{E}_i \ &= \  q\cdot \frac{\omega_i}{4\pi r^2}
    \, \chi^+ (r-t-2) \ , &\overline{H}_i \ &= \ 0
    \ . \label{charge_tensor}
\end{align}
This last expression shows that $\overline{F}$ indeed represents a
charge with asymptotic \eqref{charge_asym} which propagates in the
exterior of the light-cone $u=0$. In terms of the null decomposition
\eqref{null_decomp}, the formula \eqref{charge_tensor} is even
simpler and becomes:
\begin{align}
    \overline{\rho}\ = \ \rho(\overline{F}) \ = \ \omega^i \overline{E}_i
    \ &= \ q\cdot \frac{1}{4\pi r^2}
    \, \chi^+ (r-t-2)\ , &\overline{\alpha} = \overline{
    \underline{\alpha}} = \overline{\sigma} = 0 \ , \label{null_charge_tensor}
\end{align}
and the associated current vector\footnote{That is, plugging
$\overline{F}$ into the field equation \eqref{basic_Maxwell}.} is
easily computed to be:
\begin{align}
    \overline{J}_{\bL} \ &= \ q\cdot \frac{1}{2\pi r^2}
    \, (\chi^+)' (r-t-2) \ , &\overline{J}_{L} =
    \overline{J}_A = 0 \ , \notag
\end{align}
In particular, we have that:
\begin{equation}
        \overline{J}_0(t,r) \ = \  q\cdot \frac{1}{4\pi r^2}
        \, (\chi^+)' (r-t-2) \ . \label{charge_J0}
\end{equation}
Notice that the current vector $\overline{J}_\alpha$ is only
supported in the region $u\sim 0$ and has only a $\bL$ component.
This will allow us to obtain good estimates when it appears as an error on the
right hand side of expressions such as \eqref{basic_mor_maxwell}.\\

The main estimates we need to prove now are for the remaining field
strength:
\begin{equation}
        \tF \ = \ F - \oF \ . \label{tF_def}
\end{equation}
This satisfies the field equations \eqref{basic_Maxwell} with
current vector:
\begin{equation}
        \widetilde{J}_\alpha \ = \ J_\alpha - \overline{J}_\alpha
        \ . \label{tJ_def}
\end{equation}
A simple calculation using the formula \eqref{charge_J0} shows that
the charge for the modified field strength $\tF$ vanishes:
\begin{equation}
        q(\tF) \ = \ \int_{\RR} \ \left(J_0 - \overline{J}_0 \right) \ dx
        \ = \ 0 \ . \notag
\end{equation}
This formula, used in conjunction with certain weighted elliptic
estimates  for gradients (see appendix \ref{appendix}), allows  control of
the right hand side of \eqref{basic_mor_maxwell} above with $F$
replaced by $\tF$. We shall return to this in a moment, after we
have first discussed the previously mentioned idea of decreasing the
powers occurring in the weights in estimate
\eqref{basic_mor_maxwell}.\\

In the sequel we shall consider Maxwell fields $F_{\alpha\beta}$
such that the remainder $\tF_{\alpha\beta}$ satisfies various rates
of decay; not necessarily the amount which comes from the usual
Morawetz estimate \eqref{basic_mor_maxwell}. This can be understood
as making different assumptions as to how ``asymptotically flat''
the initial data for $F_{\alpha\beta}$ is after one has subtracted
off the contribution of the charge \eqref{charge_tensor} at time
$t=0$. The least amount of decay which is consistent with a finite
value for the charge \eqref{charge_def} is essentially
$\widetilde{E} \ \sim \ r^{-(2+\gamma)}$  for some $0 < \gamma$.
Large values for $\gamma$ are available depending on the weighted
decay of $J_0$ as well as the vanishing of higher moments for this
quantity. In order to be able to take these various rates of decay
into account for the range $0 < \gamma \leqslant 1$, we would like a
version of estimate \eqref{basic_mor_maxwell} with weights of the
form $(1+r^2)^s$ on the right hand side for $\frac{1}{2} \leqslant s
\leqslant 1$. It will be essential for us to be able to do this in
such a way that distinct decay rates, also known as peeling, are
still available for the various components in the null decomposition
\eqref{null_decomp}. This will be accomplished below through a
suitable modification of the Morawetz field \eqref{mor_field}.\\

Finally, before stating and proving the main result of this section,
let us mention in words one further adjustment to the classical
Morawetz estimate \eqref{basic_mor_maxwell} which will be of great
use in the sequel. This involves adding some space-time energies
to the left hand
side of \eqref{basic_mor_maxwell}. These additional estimates
turn out to be more than just a mere technical convenience,
and are in fact essential for us
to avoid certain logarithmic divergences which would enter when trying
to use the usual
procedure of matching fixed time energy and pure $L^\infty$ bounds.\\

\noindent We now state the basic
general energy estimate we will use for the curvature
$F_{\alpha\beta}$ which is:\\

\begin{prop}[Generalized Morawetz estimate for the electro-magnetic field
with non-vanishing charge] \label{gen_mor_F_prop} Let
$F_{\alpha\beta}$ be a two-form which satisfies the system
\eqref{basic_Maxwell} with current vector $J_\alpha$, and let $q(F)$
and $\tF$ be the associated charge and remainder field strength
defined by formulas \eqref{charge_def} and \eqref{tF_def}
respectively. Furthermore, let $0< \gamma,\epsilon$ and $\frac{1}{2}
< s \leqslant 1$ be given parameters such that $s+\gamma <
\frac{3}{2}$ as well as $\epsilon \leqslant s-\frac{1}{2}$,
and define the weights:
\begin{align}
        \tz \ &= \ \tm/\tp \ , \label{tz_def} \\
        w_{\gamma}(t,r) \ &= \  \chi^+(r-t)\cdot \tm^{2\gamma} \ + \
        (1-\chi^+(r-t)) \ , \label{w_gamma_def1} \\
    w_{\gamma,\epsilon}(t,r) \ &= \  \chi^+(r-t)\cdot \tm^{2\gamma} \ + \
        (1-\chi^+(r-t))\cdot\tm^{2\epsilon} \ , \label{w_gamma_def1.5} \\
        w'_{\gamma,\epsilon}(t,r) \ &= \  \chi^+(r-t) \cdot \tm^{2\gamma-1}+
        (1-\chi^+(r-t))\cdot \tm^{-2\epsilon -1}\ .
        \label{w_gamma_def2}
\end{align}
Define the remainder null decomposition:
\begin{subequations}
\begin{align}
        \alpha_A \ &= \ F_{LA} \ = \ \tF_{LA} \ , &\balpha_A
        \ &= \ F_{\bL A} \ = \ \tF_{\bL A}\ , \label{prop_null_decomp1}\\
        \widetilde{\rho} \ &= \ \frac{1}{2}( F_{\bL L} - \overline{F}_{\bL L})
        \ =\ \frac{1}{2}\tF_{\bL L} \ , &\sigma \ &= \ \frac{1}{2}\in^{AB}
        F_{AB} \ = \ \frac{1}{2} \in^{AB} \tF_{AB} \ .
        \label{prop_null_decomp2}
\end{align}
\end{subequations}
Now define the charge modified generalized Morawetz type energy
content of $F$ in the time slab $\{0\leqslant t \leqslant t_0\}\cap\RR$ to be:
\begin{multline}
        E_0^{(s,\gamma,\epsilon)}(0,t_0)[F] \ = \ |q(F)|^2  \ \\
    + \ \sup_{0\leqslant t \leqslant t_0} \ \int_{ \RR \cap
      \{t\} }\
    \left( \tp^{2s} |\alpha|^2 + \tm^{2s}|\balpha|^2
    + \tp^{2s}( \td{\rho}^2  + \sigma^2) \right)\, w_{\gamma} \ dx\\
    + \ \sup_u \ \int_{ \{C(u) \}\cap \{0 \leqslant t \leqslant t_0\}}
    \left( \tp^{2s} |\alpha|^2 +
    \tm^{2s}( \td{\rho}^2 + \sigma^2) \right)\, w_{\gamma}\
    dV_{C(u)}\\
    + \ \int_0^{t_0} \int_{\RR}\ \left(
    \tp^{2s} |\alpha|^2 + \tz^{1+2\epsilon}\left(\tm^{2s}|\balpha|^2
    + \tp^{2s}\td{\rho}^2  + \tp^{2s}\sigma^2\right)
    \right)\, w'_{\gamma,\epsilon}\ dxdt \ . \label{cgF_s_energy_def}
\end{multline}
Then one has the following general weighted energy type estimate:
\begin{multline}
    E_0^{(s,\gamma,\epsilon)}(0,t_0)[F] \ \leqslant \\
     C_{\gamma,\epsilon} \ \Big[
    \int_0^{t_0}\int_{ \RR}\ \left(\tp^{2s+1 +2\epsilon}\tm^{-2\epsilon}
    |J_{L}|^2 + \tp^{1+2\epsilon-2s}\tm^{4s -2\epsilon}
    |J_{\bL}|^2 + \tp^{2s} \tm |\sJ|^2 \right)\, w_{\gamma,\epsilon} \ dxdt\\
    + \ \int_{ \RR \cap
    \{t=0\} }\ (1+r^2)^{s + \gamma}\left( |E^{df}|^2 + |H|^2
    \right)\ dx \ + \lp{(1+r)^{s+\gamma} J_0(0) }{L_x^\frac{6}{5}}^2 \ \Big] \ .
    \label{main_F_L2_est}
\end{multline}
Here, $C(u)$ denotes the forward facing light-cone $t-r = u$ for
fixed values of $u$. Furthermore, $|\sJ|^2 \ = \ \delta^{AB} J_A
J_B$ denotes the angular portion of $|J|^2$.
\end{prop}\ret

\begin{proof}[Proof of estimate \eqref{main_F_L2_est}]
The proof closely follows the general strategy for proving estimates
of the type \eqref{basic_mor_maxwell} with a few additions. We first
introduce a warped generalization of the conformal killing field
\eqref{mor_field}. This is:
\begin{equation}
        K_0^s \ = \ \frac{1}{2} \bu^{2s} L + \frac{1}{2} |u|^{2s} \bL
        \ . \label{frac_mor}
\end{equation}
Here we will allow $s$ to range as $\frac{1}{2} \leqslant s
\leqslant 1$. In order to gain some intuition as to the nature of
$K_0^s$, notice that this vector-field essentially interpolates
between the conformal Killing fields $K_0$ and $S$, which are its
endpoint values in the interior region $0 \leqslant u$ for $s=1$ and
$s=\frac{1}{2}$ respectively. In order for this analogy to have
value, it must be seen that the deformation tensor of
$K_0^s$ enjoys some positivity property
when contracted with a trace-free 2-tensor which satisfies the
positive energy condition. This is indeed seen to be the case
through the following calculations. We let:
\begin{align}
        \bv(\bu) \ &= \ \bu^{2s} \ , &v(u) \ &= \ |u|^{2s} \ , \notag
\end{align}
be short hand for the fractional weights. Then, with respect to the
Minkowski metric $g_{\alpha\beta}$, the deformation tensor of $K_0^s$ is
easily calculated using the identities \eqref{def_ten_id} and
\eqref{cov_ders}:
\begin{equation}
        {}^{(K_0^s)}\pi \ = \ \frac{\bv - v}{r} \, g
    + 2 \left(\frac{\bv - v}{r} - (\dot{\bv} + \dot{v})\right)
    (\theta^{\bL} \otimes \theta^L + \theta^L\otimes\theta^{\bL})\ .
    \label{K0s_def_tensor}
\end{equation}
In the above expression, the dot notation is used to indicate the
derivative of $v$ and $\bv$ as functions of a single variable. For
example, we have that $\dot{\bv} = 2s |u|^{2s-1}\cdot \hbox{sgn}(u)$.
The positivity claim now follows once we have shown that the factor
on the second term in \eqref{K0s_def_tensor} is non-negative:\\

\begin{lem}\label{positivity_lem}
Setting let $\frac{1}{2} \leqslant s \leqslant 1$ and $\bv =
\bu^{2s}$ and $v = u^{2s}$. Then one has that:
\begin{equation}
        0 \ \leqslant \ \left(\frac{\bv-v}{r} - (\dot{\bv} +
        \dot{v})\right) \ . \label{positivity_cond}
\end{equation}
\end{lem}\ret

\begin{proof}[Proof of \eqref{positivity_cond}]
This follows immediately from freshman calculus. For fixed $t$, we
define the function of $r$:
\begin{equation}
        f(r) \ = \ \bv(t+r) - v(t-r) \ . \notag
\end{equation}
Notice that we have $f(0)=0$. We also have that $f' = \dot{\bv} +
\dot{v}$. Therefore, from the mean value theorem there exists some
$r_0\in [0,r]$ such that:
\begin{equation}
       \frac{f(r)}{r} \ = \ f'(r_0) \ = \ \dot{\bv}(t+r_0)
       + \dot{v}(t-r_0) \ . \notag
\end{equation}
The claim \eqref{positivity_cond} then follows if we know that $f'$
is a non-increasing function of $r$. This function is computed to be:
\begin{equation}
        f'(r) \ = \ 2s(r+t)^{2s-1} +
        2s|t-r|^{2s-1}\cdot\hbox{sgn}(t-r) \ . \label{f_der}
\end{equation}
For $s=\frac{1}{2}$, this function is easily seen to be decreasing by
direct inspection. For the range $\frac{1}{2} < s \leqslant 1$,
\eqref{f_der} is continuous so we only need to show that $f'' \leqslant 0$
when $0 < r$. Calculating the second derivative we see that:
\begin{equation}
        f''(r) \ = \ (2s-1)2s\left(\bu^{2s-2} - |u|^{2s-2} \right)
    \ . \notag
\end{equation}
The desired result now follows from the fact that $0 < (2s-1)2s$
while $(2s-2) \leqslant 0$ and $|u| \leqslant \bu$.
\end{proof}\ret

We now form a momentum density associated with the charge free
portion $\tF$ of $F$, as according to equations \eqref{tF_def} and
\eqref{bar_F_def}, and the extra $\gamma,\epsilon$ weights defined in
the statement of Proposition \ref{gen_mor_F_prop} above:
\begin{equation}
        P_\alpha^{(s,\gamma,\epsilon)}[F] \ = \
        Q_{\alpha\beta}[\tF](\TK)^\beta \cdot \td{w}_{\gamma,\epsilon} ,
        \label{F_mom_density_def}
\end{equation}
where $\overline{K}^0_s = T + K_0^s$ and
the weight function $\td{w}_{\gamma,\epsilon}$ is defined as:
\begin{multline}\label{w_gam_eps_weight}
        \td{w}_{\gamma,\epsilon} \ = \
        \left( \left(1+ ( 2 -u )^{2\gamma}\right)\chi^+( -u) \ + \
        \left(1+ (2+u)^{-2\epsilon}\right)(1-\chi^+(-u)) \right)
        \ + \\
        ( 1 + \bu)^{-2\epsilon}\cdot\left(  ( 2 -u )^{2\gamma+2\epsilon}
        \chi^+( -u) \ + \
        1-\chi^+(-u) \right)\ .
\end{multline}
Here $\chi^+$ can be taken to be the same smoothed out Heaviside
function used in line \eqref{bar_F_def} above. In particular,
$(\chi^+)'$ is positive and supported on the interval $[0,1]$. This
last assumption assures that the term
$(2+u)^{-2\epsilon}(1-\chi^+(-u))$ on the right hand side above is
never singular. The weight function $\td{w}_{\gamma,\epsilon}$ and its
derivatives satisfy simple bounds with respect to the weight
functions $w_\gamma$ and $w'_{\gamma,\epsilon}$ defined on lines
\eqref{w_gamma_def1}--\eqref{w_gamma_def2} above. First, notice that
one has the bounds:
\begin{equation}
        C^{-1} w_\gamma \ \leqslant\ \td{w}_{\gamma,\epsilon} \ \leqslant
        \  C w_\gamma \ , \notag
\end{equation}
for a suitable positive constant $C$. We will restate this estimate
as follows:
\begin{equation}
        \td{w}_{\gamma,\epsilon} \ \sim \
        w_\gamma \ . \label{w_gam_eps_bound1}
\end{equation}
Furthermore, a brief calculation shows that:
\begin{align}
        - \frac{1}{2} L( \td{w}_{\gamma,\epsilon} ) \ &= \
        -\partial_{\bu}(\td{w}_{\gamma,\epsilon}) \ , \notag \\
        &= \ 2\epsilon\, ( 1 + \bu)^{-2\epsilon-1}\cdot\Big(  ( 2 -u )^{2\gamma+2\epsilon}
        \chi^+( -u) \notag \ + \
        1-\chi^+(-u) \Big) \ , \notag \\
        &\sim \ C^{-1}_{\gamma,\epsilon}\ \tz^{1 + 2\epsilon} w'_{\gamma,\epsilon} \ .
        \label{w_gam_eps_bound2}
\end{align}
A similar calculation shows that:
\begin{align}
        -\frac{1}{2} \bL(\td{w}_{\gamma,\epsilon}) \ &= \
    -\partial_u (\td{w}_{\gamma,\epsilon}) \ , \notag\\
        &= \ \left( (2 -u )^{2\gamma} - (2 + u)^{-2\epsilon}
        \right)(\chi^+)'(-u)  \ \ + \notag \\
        &\ \ \ \ \ (1+\bu)^{-2\epsilon}\left( (2 -u )^{2\gamma+2\epsilon}
        - 1 \right)(\chi^+)'(-u) \ \ + \notag \\
        &\ \ \ \  \ 2(\gamma + \epsilon)(1 + \bu)^{-2\epsilon}
        ( 2 -u )^{2\gamma+2\epsilon-1}\chi^+( -u) \ \ + \notag \\
    &\ \ \ \ \ 2\gamma\, ( 2 -u )^{2\gamma-1}\chi^+( -u) \ + \
        2\epsilon\, (2+u)^{-2\epsilon-1}(1-\chi^+(-u))  \ , \notag \\
        &\sim \ C^{-1}_{\gamma,\epsilon} \ w'_{\gamma,\epsilon} \ .
        \label{w_gam_eps_bound3}
\end{align}
The bound on the last line follows because one has:
\begin{equation}
        (2 + u)^{\delta_1} \ \leqslant \ (2 -u)^{\delta_2} \ ,
        \notag
\end{equation}
whenever $\delta_1\leqslant\delta_2$ and $u\in[-1,0]$. Note again that
$(\chi^+)'$ is a positive function. Also, it is
clear that the constant implicit in the $\sim$ notation on lines
\eqref{w_gam_eps_bound1}--\eqref{w_gam_eps_bound3} above can be
taken to be the same value uniform in $\epsilon$ and $\gamma$
because we kept the effect of these constants in our bounds with the
$C_{\gamma,\epsilon}$ notation. We
will continue to use the $\sim$ notation in the remainder of this
proof with the understanding that the same implicit constant is being
used throughout.\\

Now, calculating the space-time divergence of the quantity
\eqref{F_mom_density_def} we have that:
\begin{multline}
        \nabla^\alpha P_\alpha^{(s,\gamma,\epsilon)}[F] \ = \
        \tF_{\alpha\beta} \td{J}^\alpha (\TK)^\beta\cdot
        \td{w}_{\gamma,\epsilon} \ + \ \frac{1}{2} Q_{\alpha\beta}[\tF]
        \, {}^{(\TK)}\pi^{\alpha\beta}\cdot \td{w}_{\gamma,\epsilon} \\
        - \ \frac{1}{2} Q_{\bL \beta}[\tF] (\TK)^\beta \cdot
        L(\td{w}_{\gamma,\epsilon}) \ - \ \frac{1}{2} Q_{L \beta}[\tF]
        (\TK)^\beta \cdot \bL(\td{w}_{\gamma,\epsilon}) \ .
        \label{main_F_div_id}
\end{multline}
Integrating both sides of this last line over time slabs of the
form:
\begin{equation}
       \mathcal{R}(t_0,u_0) \ = \ \{0\leqslant t \leqslant
        t_0\}\cap\{u\leqslant u_0\} \ , \notag
\end{equation}
and applying the geometric version of the Stokes theorem (Gauss
theorem) we arrive at the following general integral identity:
\begin{multline}
        \int_{\{t=0\}\cap \{ -u_0 \leqslant r \}} \
        P_0^{(s,\gamma,\epsilon)}[F]\ dx \
        = \ \int_{\{t=t_0\}\cap\{ t_0 - u_0\leqslant r \}} \
        P_0^{(s,\gamma,\epsilon)}[F]\ dx \\
        + \ \int_{ C(u_0) \cap \{0\leqslant t \leqslant t_0 \} }
        \ P_L^{(s,\gamma,\epsilon)}[F]\ dV_{C(u_0)}
        \ + \ \int\int_{\mathcal{R}(t_0,u_0)} \
        \hbox{(R.H.S.)}\eqref{main_F_div_id} \ .
        \label{main_F_int_id}
\end{multline}
Notice that the above identity reduces to the usual energy type
estimate on the time slab $\{0\leqslant t \leqslant t_0\}\cap\RR$
when $t_0 \leqslant u_0$. In order to proceed, we now calculate each
of the terms in \eqref{main_F_int_id} individually.\\

When $t=0$, we compute that:
\begin{equation}
        P_0^{(s,\gamma,\epsilon)}[F] \ \sim \ (1 + r^2)^{s+\gamma}(
        |\td{E}|^2 + |H|^2) \ , \label{PF_0_est}
\end{equation}
where again the $\sim$ notation here means that the ratio of the two
terms above is bounded by on the left and right by $C^{-1}$ and $C$
respectively for a suitable \emph{positive} constant $C$. In
particular, if the right hand side of such an expression is
positive, then so is the left. Also, $\td{E}$ denotes the electric
part of the tensor $\tF$ while $H$ is the magnetic part of $F$. Note
that:
\begin{equation}
        \td{E} \ =\  E - \overline{E} \ , \notag
\end{equation}
where $\overline{E}$ is defined as on line \eqref{charge_tensor}
above. Likewise, with the help of the identities
\eqref{F_em_decomp}, we compute that at time $t=t_0$:
\begin{equation}
        P_0^{(s,\gamma,\epsilon)}[F] \ \sim \ \left( \tp^{2s}
        |\alpha|^2 + \tm^{2s}|\balpha|^2 + \tp^{2s}( \td{\rho}^2
        + \sigma^2) \right) \cdot w_{\gamma} \ . \label{PF_t0_est}
\end{equation}
We also compute the characteristic energy term:
\begin{equation}
        P_L^{(s,\gamma,\epsilon)}[F] \ \sim \ \left( \tp^{2s}
        |\alpha|^2 + \tm^{2s} ( \td{\rho}^2 + \sigma^2)\right)\cdot
        w_\gamma \ . \label{PF_L_est}
\end{equation}\ret

It remains to calculate the terms on the right hand side of
\eqref{main_F_div_id} above. Since the fist such term does not have
a sign, we simply expand it using the null decomposition
\eqref{prop_null_decomp1}--\eqref{prop_null_decomp2} and the
definition of $\td{J}$ given by \eqref{tJ_def} above as well as the
estimate \eqref{w_gam_eps_bound1} to bound:
\begin{multline}
        | \tF_{\alpha\beta} \td{J}^\alpha (\TK)^\beta\cdot
        \td{w}_{\gamma,\epsilon}| \ \lesssim \\
        \Big( \tp^{2s} (|\alpha|\cdot
        |\sJ| + |\td{\rho}|\cdot |J_L| ) \
        + \ \tm^{2s}(|\balpha|\cdot|\sJ|
        + |\td{\rho}|\cdot |J_{\bL}|) \Big)\cdot w_\gamma \ + \
        |\td{\rho}|\cdot|q|\cdot \tp^{-2} \tm^{-10} w_\gamma \ .
        \label{Q_F_div_exp1}
\end{multline}
Notice that the extra $\tm$ weight on the last term of the right
hand side above comes because $\overline{J}$ is supported in the
region where $\tm \sim 1$.\\

We now calculate the second term on the right hand side of
\eqref{main_F_div_id} above. Because $Q_{\alpha\beta}[F]$ is
trace-free and $T$ is Killing with respect to $g_{\alpha\beta}$, the
contraction need only be taken with the second term on the right
hand side of line \eqref{K0s_def_tensor}. This yields:
\begin{equation}
        Q_{\alpha\beta}[\tF]
        \, {}^{(\TK)}\pi^{\alpha\beta}\cdot \td{w}_{\gamma,\epsilon} \
        \sim \ \left(\frac{\bv - v}{r} - (\dot{\bv} +
        \dot{v})\right)\cdot( \td{\rho}^2 + \sigma^2 )\cdot w_\gamma
        \ . \label{Q_F_div_exp2}
\end{equation}
In particular, using Lemma \ref{positivity_lem} and the range
restriction $\frac{1}{2} \leqslant s \leqslant 1$ we have that this
last expression is non-negative.\\

Moving on, we calculate the third term on the right hand side of
\eqref{main_F_div_id}. Using the bound \eqref{w_gam_eps_bound2} and
the decomposition \eqref{F_em_decomp} we conclude that:
\begin{equation}
        - \ \frac{1}{2} Q_{\bL \beta}[\tF] (\TK)^\beta \cdot
        L(w_{\gamma,\epsilon}) \ \sim \ C^{-1}_{\gamma,\epsilon}\
    \tz^{1+2\epsilon}\left( \tm^{2s}
        |\balpha|^2 + \tp^{2s}( \td{\rho}^2 + \sigma^2)\right)
        \cdot w'_{\gamma,\epsilon} \ . \label{Q_F_div_exp3}
\end{equation}
Similarly, using the bound \eqref{w_gam_eps_bound3} we estimate:
\begin{equation}
        - \ \frac{1}{2} Q_{L \beta}[\tF] (\TK)^\beta \cdot
        \bL(w_{\gamma,\epsilon}) \ \sim \ C^{-1}_{\gamma,\epsilon}
    \ \left(\tp^{2s}|\alpha|^2 +
        \tm^{2s}( \td{\rho}^2 + \sigma^2)\right)\cdot w'_{\gamma,\epsilon} \ .
        \label{Q_F_div_exp4}
\end{equation}\ret

Therefore, collecting the positive terms
\eqref{PF_t0_est}--\eqref{PF_L_est} and
\eqref{Q_F_div_exp2}--\eqref{Q_F_div_exp4} onto the right hand side
of \eqref{main_F_int_id} above, and collecting \eqref{Q_F_div_exp1}
with \eqref{PF_0_est} on the left hand side and using the shorthand
notation \eqref{F_s_energy_def}, which will be introduced in a
moment,  we can estimate:
\begin{multline}
        \int_{\{t=0\}\cap \RR} \ (1+r^2)^{s+\gamma}( |\td{E}|^2 +
        |H|^2)\ dx \ + \ \int_0^{t_0}\int_{\RR}\ |\td{\rho}|\cdot|q|
        \cdot \tp^{-2} \tm^{-10} w_\gamma\ dxdt \  \\
        + \ \int_0^{t_0}\int_{\RR}\ \Big( \tp^{2s} (|\alpha|\cdot
        |\sJ| + |\td{\rho}|\cdot |J_L| ) \
        + \ \tm^{2s}(|\balpha|\cdot|\sJ|
        + |\td{\rho}|\cdot |J_{\bL}|) \Big)\cdot w_\gamma \ dxdt \\
        \geqslant \ C^{-1}_{\gamma,\epsilon}\ \Big[\
        E^{(s,\gamma,\epsilon)}(0,t_0)[\tF] \ + \
        \int_{0}^{t_0}\int_{\RR}\ \left(\frac{\bv - v}{r} - (\dot{\bv} +
        \dot{v})\right)\cdot( \td{\rho}^2 + \sigma^2 )\cdot w_\gamma
        \ dxdt \ \Big] \ . \label{main_F_L2_est_intermediate}
\end{multline}
We now make a preliminary reduction on estimate
\eqref{main_F_L2_est_intermediate} above by first discarding the
second positive term on the right hand side, and then using a
Cauchy--Schwartz to peel off a factor of
$(E^\frac{1}{2})^{(s,\gamma,\epsilon)}(0,t_0)[\tF]$ from the latter
two terms on the left hand side of the above expression. Setting:
\begin{multline}
        ||| J |||^2_{L^2[0,t_0](L^2)(s,\gamma,\epsilon) } \ = \\
        \int_0^{t_0}\int_{ \RR}\ \left(\tp^{2s} \tz^{-1-2\epsilon}\tm
        |J_{L}|^2 + \tz^{2s-1-2\epsilon}\tm^{2s + 1}
        | J_{\bL}|^2 + \tp^{2s} \tm |\sJ|^2 \right) \
        w_{\gamma,\epsilon}\ dxdt \ , \label{J_norm}
\end{multline}
by a simple use of H\"olders inequality, the condition $\epsilon\leqslant
s-\frac{1}{2}$, and the weight bound $w_\gamma^2\leqslant \tm
w'_{\gamma,\epsilon} w_{\gamma,\epsilon}$\
we may replace \eqref{main_F_L2_est_intermediate} with the estimate:
\begin{multline}
        E^{(s,\gamma,\epsilon)}(0,t_0)[\tF] \ \leqslant \
        C_{\gamma,\epsilon} \Big[
        (E^\frac{1}{2})^{(s,\gamma,\epsilon)}(0,t_0)[\tF]\cdot
        \left( ||| J |||_{L^2[0,t_0](L^2)(s,\gamma,\epsilon) } + |q| \right) \\
        + \ \int_{\{t=0\}\cap \RR} \ (1+r^2)^{s+\gamma}( |\td{E}|^2 +
        |H|^2)\ dx \ \Big] \ , \notag
\end{multline}\
from which easily follows the bound:
\begin{equation}
        E^{(s,\gamma,\epsilon)}(0,t_0)[\tF] \ \leqslant \
        C^2_{\gamma,\epsilon} \Big[
        ||| J |||^2_{L^2[0,t_0](L^2)(s,\gamma,\epsilon) }
        \ + \ |q|^2 \ + \
         \int_{\{t=0\}\cap \RR} \ (1+r^2)^{s+\gamma}( |\td{E}|^2 +
        |H|^2)\ dx \ \Big] \ , \label{main_F_L2_almost}
\end{equation}\ret

In order to deduce \eqref{main_F_L2_est} from
\eqref{main_F_L2_almost}, we need to prove the following bounds for
the range $\frac{1}{2} < s \leqslant 1$ and $0 < \gamma$,
where $s+\gamma < \frac{3}{2}$:
\begin{equation}
        |q|^2 \ \leqslant \ C_\gamma \lp{(1+r)^{s+\gamma}
        J_0(0)}{L^\frac{6}{5}_x}^2 \ , \label{charge_L65_est}
\end{equation}
and:
\begin{multline}
        \int_{\{t=0\}\cap\RR} \ (1 + r^2)^{s+\gamma} |\td{E}|^2\
        dx\\
        \leqslant \ C_\gamma\ \left[ \ \int_{\{t=0\}\cap\RR} \
        (1 + r^2)^{s+\gamma} |E^{df}|^2\ dx \ + \
        \lp{(1+r)^{s+\gamma}J_0(0)}{L^\frac{6}{5}_x}^2 \ \right] \ .
        \label{E_L65_J0_est}
\end{multline}\ret

The first of the above estimates, \eqref{charge_L65_est}, follows
from a simple application of H\"olders inequality with the weight
$(1+r)^{-s-\gamma}$ in the integral \eqref{charge_def} at time
$t=0$. Note that $\frac{1}{2} < s + \gamma$ because of our range
restrictions so the resulting factor integral converges with a bound
depending only on $\gamma$.\\

The second of the above estimates results from a Hodge decomposition
and the weighted elliptic estimate \eqref{weighted_grad_sob} from Appendix
\ref{appendix}. Expanding out $\td{E}$ on the left hand side we can
bound:
\begin{equation}
        \int_{\{t=0\}\cap\RR} \ (1 + r^2)^{s+\gamma} |\td{E}|^2\
        dx \ \leqslant \ 2\int_{\{t=0\}\cap\RR} \ (1 + r^2)^{s+\gamma} \left(|E^{df}|^2
        + |E^{cf} - \overline{E}|^2\right)\ dx
         \ . \notag
\end{equation}
Introducing a potential function $\varphi$ for $E^{cf}$ as in line
\eqref{charge_eq} above and recalling the definition
\eqref{charge_tensor} of $\overline{E}$, and using line
\eqref{charge_L65_est} to estimate:
\begin{equation}
        \int_{\{t=0\}\cap\RR} \ (1 + r^2)^{s+\gamma} \big| \frac{q(F)
        }{r}\nabla_x \chi^+(r-2) \big|^2 \ dx \ \leqslant \
        C_\gamma \ |q(F)|^2 \ \leqslant \ C_\gamma \ \lp{(1+r)^{s+\gamma}
        J_0(0)}{L^\frac{6}{5}_x}^2 \ , \notag
\end{equation}
we are reduced to proving the estimate:
\begin{multline}
        \int_{\{t=0\}\cap\RR} \ (1 + r^2)^{s+\gamma} \big|\nabla\big(
        \frac{1}{\Delta} J_0(0) + \frac{q(F)}{4\pi r}
        \chi^+(r-2)\big)\big|^2\ dx \\
    \lesssim \ C_\gamma \
        \lp{(1+r)^{s+\gamma}J_0(0)}{L^\frac{6}{5}_x}^2
    \ . \label{reduced_wt_grad_sob}
\end{multline}
By using the usual $L^\frac{6}{5}\hookrightarrow L^2$ Sobolev
embedding, as well as the bounds:
\begin{align}
        \int_{\{t=0\}\cap\RR} \ \big|\nabla\big( \frac{q(F)}{4\pi r}
        \chi^+(r-2)\big)\big|^2 \ dx \ &\lesssim \ |q(F)|^2 \ \leqslant
        \ C_\gamma \ \lp{(1+r)^{s+\gamma}
        J_0(0)}{L^\frac{6}{5}_x}^2 \ , \notag\\
    \int_{\{t=0\}\cap\RR} \ r^{2(s+\gamma)}\ \big|\nabla\big(\frac{q(F)}{4\pi
        r}(1-\chi^+(r-2))\big) \big|^2 \ dx \ &\leqslant  \ C_\gamma \
        |q(F)|^2 \ \leqslant
        \ C_\gamma \ \lp{(1+r)^{s+\gamma}
        J_0(0)}{L^\frac{6}{5}_x}^2 \ , \notag
\end{align}
the second of which follows from the condition $\frac{1}{2} <
s+\gamma$, we are reduced to showing that:
\begin{equation}
        \int_{\{t=0\}\cap\RR} \ r^{2(s+\gamma)} \big|\nabla\big(
        \frac{1}{\Delta} J_0(0) + \frac{q(F)}{4\pi r} \big)\big|^2\ dx \
    \lesssim \ C_\gamma \
        \lp{r^{s+\gamma}J_0(0)}{L^\frac{6}{5}_x}^2
    \ . \label{E_L65_J0_est_red}
\end{equation}
Using the definition \eqref{charge_def} of $q(F)$, this is precisely
the statement of \eqref{weighted_grad_sob} in Appendix
\ref{appendix}. Notice that the
condition $\frac{1}{2} < s+\gamma < \frac{3}{2}$ is enforced by
assumption. This ends the
proof of estimate \eqref{main_F_L2_est}.
\end{proof}\ret

We conclude this section with an important generalization of
Proposition \ref{gen_mor_F_prop}. This involves estimates for
derivatives of the field strength $F_{\alpha\beta}$. As it turns
out, the simple decomposition \eqref{tF_def} is remarkably robust
with respect to the operation of Lie differentiation. One can show
that $\mathcal{L}_X^I \tF_{\alpha\beta}$ satisfies bounds similar to
\eqref{main_F_L2_est} with the appropriate right hand side. This
will be used in the sequel to show that the various components of
the null decomposition of $\td{F}_{\alpha\beta}$ satisfy the
expected $L^\infty$ estimates. Together with the fact that
$\overline{F}_{\alpha\beta}$ is given explicitly, these bounds will
fully  determine the point-wise behavior of the original field
strength
$F_{\alpha\beta}$. We will prove:\\

\begin{prop}[Generalized Morawetz estimate for derivatives of the
electro-magnetic field ] \label{gen_mor_F_prop2}
Let $F_{\alpha\beta}$ be a two-form which satisfies the system
\eqref{basic_Maxwell} with current vector $J_\alpha$, and let $q(F)$
and $\tF$ be the associated charge and remainder field strength
defined by formulas \eqref{charge_def} and \eqref{tF_def}
respectively. Furthermore, let $0< \gamma,\epsilon$ and $\frac{1}{2}
< s \leqslant 1$ be given parameters such that $s+\gamma <
\frac{3}{2}$ and $\epsilon \leqslant  s-\frac{1}{2}$. Define the
weights $\tz, w_\gamma , w_{\gamma,\epsilon}, w'_{\gamma,\epsilon}$ as
on  lines \eqref{tz_def}--\eqref{w_gamma_def2} above.
Next, define the generalized (not charge modified) Morawetz type energy
content of a two-from $F$ in the time slab
$\{0\leqslant t \leqslant t_0\}\cap\RR$ to be:
\begin{multline}
        E^{(s,\gamma,\epsilon)}(0,t_0)[F] \ = \
     \ \sup_{0\leqslant t \leqslant t_0} \ \int_{ \RR \cap
      \{t\} }\
    \left( \tp^{2s} |\alpha|^2 + \tm^{2s}|\balpha|^2
    + \tp^{2s}( \rho^2  + \sigma^2) \right)\, w_{\gamma} \ dx\\
    + \ \sup_u \ \int_{ \{C(u) \}\cap \{0 \leqslant t \leqslant t_0\}}
    \left( \tp^{2s} |\alpha|^2 +
    \tm^{2s}( \rho^2 + \sigma^2) \right)\, w_{\gamma}\
    dV_{C(u)}\\
    + \ \int_0^{t_0} \int_{\RR}\ \left(
    \tp^{2s} |\alpha|^2 + \tz^{1+2\epsilon}\left(\tm^{2s}|\balpha|^2
    + \tp^{2s}\rho^2  + \tp^{2s}\sigma^2\right)
    \right)\, w'_{\gamma,\epsilon}\ dxdt \ . \label{F_s_energy_def}
\end{multline}
Now, define the $k^{th}$ weighted charge modified Morawetz energy:
\begin{equation}
        E_k^{(s,\gamma,\epsilon)}(0,t_0)[F] \ = \ q^2(F) \ + \
    \sum_{\substack{ |I|\leqslant k \\ X\in \mathbb{L}}}\
    E^{(s,\gamma,\epsilon)}(0,t_0)[\mathcal{L}^I_X \tF] \ . \label{k_F_eng}
\end{equation}
Then, recalling the vector norm $\llp{\cdot}{
L^2[0,t_0](L^2)(s,\gamma,\epsilon)}$ from line \eqref{J_norm} above,
we have the following general energy estimate:
\begin{multline}
    E_k^{(s,\gamma,\epsilon)}(0,t_0)[F] \ \leqslant \
    C_{\gamma,\epsilon}^2 \ \Big[ \ \sum_{\substack{
    |I|\leqslant k\\ X\in\mathbb{L}} } \ \Big( \ \
    \llp{\mathcal{L}^I_X J }{ L^2[0,t_0](L^2)(s,\gamma,\epsilon)}^2\\
    + \ \int_{ \RR \cap
    \{t=0\} }\ (1+r^2)^{s + \gamma + |I|}\left( |\nabla_x^I E^{df}|^2
    + |\nabla_x^I H|^2
    \right)\ dx \ \
    + \ \ \lp{(1+r)^{s+\gamma +|I|} \nabla_x^I
    J_0(0) }{L_x^\frac{6}{5}}^2 \  \Big) \\
    + \  \sum_{|I|\leqslant k-1}\ \lp{(1+r)^{s+\gamma +|I|+1} \nabla_{t,x}^I
    J(0) }{L_x^2}^2\ \Big] \ . \label{main_F_L2_est2}
\end{multline}
Here $\nabla_x$ and $\nabla_{t,x}$ denote the Lie derivatives of the
covector $J$ and the two-form $F$ with respect to the spatial
translation invariant fields $\{\partial_i\}$ and the full set of translation
invariant fields $\{\partial_\alpha\}$ respectively.
\end{prop}\ret

\begin{proof}[Proof of estimate \eqref{main_F_L2_est2}]
The proof is essentially that of estimate \eqref{main_F_L2_est}
applied to the field $\mathcal{L}_X^I \tF_{\alpha\beta}$ with some additional
calculations at the end to wrap things up.
First notice that by \eqref{F_div_iden} and the fact that for any Lie
derivative one has $[d,\mathcal{L}_X] = 0$, we have the following
formula for $\mathcal{L}^I_X \tF$ whenever $X\in \mathbb{L}\setminus \{S\}$:
\begin{align}
        \nabla^\beta (\mathcal{L}_X^I \tF)_{\alpha\beta} \ &= \
    (\mathcal{L}^I_X J)_\alpha - (\mathcal{L}^I_X
    \overline{J})_\alpha \ , \notag\\
    \nabla^\beta \, {}^*(\mathcal{L}_X^I \tF) \ &= \ 0 \ . \notag
\end{align}
The effect of the Lie derivative $\mathcal{L}_S$ is equally easy to account for. In
this case formula \eqref{F_div_iden} gives:
\begin{align}
        \nabla^\beta (\mathcal{L}_S \tF)_{\alpha\beta} \ &= \
    (J -\overline{J})_\alpha \ + \
    (\mathcal{L}_S J)_\alpha - (\mathcal{L}_S
    \overline{J})_\alpha \ , \notag\\
    \nabla^\beta \, {}^*(\mathcal{L}_S \tF) \ &= \ 0 \ . \notag
\end{align}
Using the above formulas and following the proof of
\eqref{main_F_L2_est} until line \eqref{main_F_L2_est_intermediate}
yields the inequality:
\begin{multline}
        E_k^{(s,\gamma,\epsilon)}(0,t_0)[F] \ \leqslant \
        C_{\gamma,\epsilon}\ \sum_{\substack{ |I|\leqslant k \\ X\in \mathbb{L}}}\
    \Big[ \ |q(F)|^2 \ + \ \llp{\mathcal{L}^I_X(J - \overline{J})}
    {L^2[0,t_0](L^2)(s,\gamma,\epsilon)}^2 \\
    + \ \int_{ \RR \cap \{t=0\} }\
    (1+r^2)^{s + \gamma }\left( |E(\mathcal{L}_X^I \tF)|^2
    + |H(\mathcal{L}_X^I \tF)|^2
    \right)\ dx \ \Big] \ . \notag
\end{multline}
Using the same steps as at the end of the proof of
\eqref{main_F_L2_est} to bound the quantity $|q(F)|$, we are done once
we have shown that:
\begin{equation}
        \llp{\mathcal{L}^I_X\overline{J}}
    {L^2[0,t_0](L^2)(s,\gamma,\epsilon)}^2 \ \lesssim \ |q(F)|^2 \ ,
    \label{oJ_bound}
\end{equation}
as well as the bound:
\begin{multline}
        \sum_{\substack{ |I|\leqslant k \\ X\in \mathbb{L}}}\
    \int_{ \RR \cap \{t=0\} }\
    (1+r^2)^{s + \gamma }\left( |E(\mathcal{L}_X^I \tF)|^2
    + |H(\mathcal{L}_X^I \tF)|^2
    \right)\ dx \\
    \lesssim \ \ \ \ \ \ \
    \sum_{ |I|\leqslant k }
    \ \Big[ \ \int_{ \RR \cap
    \{t=0\} }\ (1+r^2)^{s + \gamma + |I|}\left( |\nabla_x^I E^{df}|^2
    + |\nabla_x^I H|^2 \right)\ dx \\
    + \ \lp{(1+r)^{s+\gamma +|I|} \nabla_x^I
    J_0(0) }{L_x^\frac{6}{5}}^2 \ \Big]\\
    + \   \sum_{ |I|\leqslant k-1 }
    \lp{(1+r)^{s+\gamma +|I|+1} \nabla_{t,x}^I
    J(0) }{L_x^2}^2  \ . \label{k_F_data_est}
\end{multline}\ret

We begin with \eqref{oJ_bound}. This will follow from direct
computation using point-wise bounds on the quantity
$\mathcal{L}^I_X\overline{J}$ for $X\in\mathbb{L}$. These point-wise
bounds will be provided through induction on the value $|I|$. Notice
that when $|I|=0$ we may write:
\begin{align}
        J_{\bL} \ &= \ q\cdot \frac{\Omega_{\bL} (\omega) }{r^2}
    \cdot \chi(u) \ ,
    &J_{L} \ &= \ q\cdot \frac{\Omega_{L} (\omega) }{r^4}
    \cdot \chi(u)\ ,
    &J_{A} \ &= \ q\cdot \frac{\Omega_{A} (\omega) }{r^3}
    \cdot \chi(u)\ , \label{I0_line}
\end{align}
where $\chi$ is a $C^\infty$ and $O(1)$ bump function adapted to the
origin $u=0$ and the $\Omega_\alpha(\omega)$ are smooth function of
the angular variable only. In particular, these satisfy inductive
identities:
\begin{align}
        L(\Omega_\alpha) \ &= \ \bL(\Omega_\alpha) \ = \ 0 \ ,
    &e_A(\Omega_\alpha) \ = \ \frac{1}{r} \,
    \td{\Omega}_\alpha \ . \notag
\end{align}
where the $\td{\Omega}_\alpha$ have the same properties of the
$\Omega_\alpha$. Also, notice that by direct computation one can
substitute the decay rates \eqref{I0_line} into the norm $\llp{\cdot
} {L^2[0,t_0](L^2)(s,\gamma,\epsilon)}^2$ to achieve a bound in
terms of $|q(F)|$. We now show that the schematics on line
\eqref{I0_line} are preserved after each round of Lie
differentiating by inductively establishing the identities:
\begin{subequations}\label{I_line}
\begin{align}
        (\mathcal{L}^I_{X}J)_{\bL} \ &= \ q\cdot \sum_{k=0}^{|I|}
        \frac{\Omega^I_{\bL,k}
        (\omega) }{r^{2+k}} \cdot (\chi_k)_{\bL}^I(u) \ ,  \\
    (\mathcal{L}^I_X J)_{L} \ &= \ q\cdot \sum_{k=0}^{|I|}
    \frac{\Omega^I_{L,k}
        (\omega) }{r^{4+k}}\cdot(\chi_k)_L^I(u) \ , \\
    (\mathcal{L}^I_X J)_{A} \ &= \ q\cdot \sum_{k=0}^{|I|}
    \frac{\Omega^I_{A,k}
        (\omega) }{r^{3+k}}\cdot(\chi_k)_A^I(u) \ ,
\end{align}
\end{subequations}
where the $\Omega^I_{\alpha,k}$ and $(\chi_k)_\alpha^I$ depend on the
specific combination and number of vector-fields,
but have the same properties
as the $\Omega_\alpha$ and $\chi$ in line \eqref{I0_line} above.
Also, the product notation is symbolic and is used to denote a
finite sum of products of functions with these properties. Assuming
that \eqref{I_line} is true for $|I|=l-1$ it can be shown that
\eqref{I_line} holds for $|I|=l$ through direct use of the
identities \eqref{lor_null_decomp}--\eqref{lor_null_brack} as well
as the formula:
\begin{equation}
        (\mathcal{L}_X J)_\alpha \ = \ X(J_\alpha) - J([X,e_\alpha])
        \ . \notag
\end{equation}
We leave this as a straightforward,
although rather tedious, exercise for the reader. Note that everything
can be put in terms of the variables $u,r,\omega$ so the only thing to
verify is the homogeneity claim with respect to the power of $r$.\\

It remains to prove the second bound \eqref{k_F_data_est} above.
This is done in two separate steps. First of all, a direct
computation involving the formula:
\begin{equation}
        (\mathcal{L}_X \tF)_{\alpha\beta} \ = \ X(\tF_{\alpha\beta})
        - \tF([X,e_\alpha],e_\beta) - \tF(e_\alpha, [X,e_\beta])
    \ , \notag
\end{equation}
as well as the bracket identities \eqref{lor_bracket_rel1} and
\eqref{lor_bracket_rel3} shows that at time $t=0$ one has the
point-wise bound:
\begin{equation}
        \sum_{\substack{ |I|\leqslant k \\ X\in \mathbb{L}}}\
        \left( |E(\mathcal{L}_X^I \tF)|^2
    + |H(\mathcal{L}_X^I \tF)|^2
    \right) \ \lesssim \
    \sum_{ |I| \leqslant k}\
    (1+r^2)^{|I|} \left(
    |\nabla_{t,x}^I \td{E}|^2 + |\nabla_{t,x}^I H|^2
    \right) \ , \notag
\end{equation}
Using the field equations
\eqref{EH_Maxwell}, this last line can be further reduced to the
estimate:
\begin{multline}
         \sum_{\substack{ |I|\leqslant k \\ X\in \mathbb{L}}}\
        \left( |E(\mathcal{L}_X^I \tF)|^2
    + |H(\mathcal{L}_X^I \tF)|^2
    \right) \ \lesssim \ \sum_{ |I| \leqslant k}\
    (1+r^2)^{|I|} \left(
    |\nabla_{x}^I \td{E}|^2
    +|\nabla_x^I E^{df}|^2 + |\nabla_{x}^I H|^2 \right)\\
    + \ \sum_{ |I| \leqslant k-1} \ (1+r^2)^{|I|+1}\ \left(
    |\nabla_{t,x}^I J|^2 + |q(F)|^2\cdot |\nabla_{t,x}^I(\chi^+)'|^2
    \right) \ . \notag
\end{multline}
Multiplying through by $(1+r^2)^{s+\gamma}$ and integrating this last
line over $\RR$, we have achieved the bound \eqref{k_F_data_est} modulo the
estimate:
\begin{multline}
        \sum_{ |I|\leqslant k } \ \int_{ \RR \cap
    \{t=0\} }\ (1+r^2)^{s + \gamma + |I|}|\nabla_x^I \td{E}|^2
    \ dx  \ \lesssim \
     \sum_{ |I|\leqslant k } \ \Big( \int_{ \RR \cap
    \{t=0\} }\ (1+r^2)^{s + \gamma + |I|}|\nabla_x^I E^{df}|^2
    \ dx \\
    + \ \ \ \ \ \ \lp{(1+r)^{s+\gamma +|I|} \nabla_x^I
    J_0(0)}{L^\frac{6}{5}}^2 \ \Big) \ . \notag
\end{multline}
Using essentially the same steps which were employed to reduce
estimate \eqref{E_L65_J0_est} above to \eqref{E_L65_J0_est_red},
this last line is a consequence
of the following generalization of estimate \eqref{weighted_grad_sob}
in the Appendix:
\begin{multline}
        \sum_{|I|\leqslant k} \
    \int_{\RR} \ r^{2(s+\gamma+|I|)}\ \big|
    \nabla \ \nabla_x^I \big( \frac{1}{\Delta}J_0(0) + \frac{q}{4\pi
    r} \big) \big|^2 \ dx \ \\
    \leqslant \ C_\gamma \
    \sum_{|I|\leqslant k} \ \lp{r^{s+\gamma +|I|} \nabla_x^I
    J_0(0)}{L^\frac{6}{5}}^2 \ . \label{gen_weighted_grad_sob}
\end{multline}
Estimate \eqref{gen_weighted_grad_sob} can be reduced to estimate
\eqref{E_L65_J0_est_red} through a process of induction. To see
this, we assume that $0 < |I|$ and integrate by parts a couple of
times on the left hand
side of \eqref{gen_weighted_grad_sob}. This yields:
\begin{multline}
         \hbox{(L.H.S.)}\eqref{gen_weighted_grad_sob} \ = \
     \sum_{|I|\leqslant k} \ \Big( \ (s+\gamma+|I|+1)(s+\gamma+|I|)\
     \int_{\RR} \ r^{2(s+\gamma+|I|-1)}\ \big|
     \nabla_x^I \big( \frac{1}{\Delta}J_0(0) + \frac{q}{4\pi
     r} \big) \big|^2 \ dx \\
     - \  \int_{\RR} \ r^{2(s+\gamma+|I|)}\ (\nabla^I_x J_0(0))\cdot
     \nabla_x^I \big( \frac{1}{\Delta}J_0(0) + \frac{q}{4\pi
     r} \big)  \ dx \ \Big) \ . \notag
\end{multline}
Applying Cauchy-Schwartz and H\"olders inequality to the right hand
side of the above expression in conjunction with the
$L^2\hookrightarrow L^6$ Sobolev embedding and setting:
\begin{align}
        A(k) \ &= \ \sum_{|I|\leqslant k} \ \lp{r^{s+\gamma + |I|}
    \nabla \
    \nabla_x^J\big( \frac{1}{\Delta}J_0(0)
    + \frac{q}{4\pi r} \big) }{L^2} \ , \notag\\
    B(k) \ &= \ \sum_{|I|\leqslant k}\
     \lp{r^{s+\gamma +|I|} \nabla_x^I J_0(0)}{L^\frac{6}{5}}  \ , \notag
\end{align}
we have for $1\leqslant k$ the estimate:
\begin{equation}
        |\hbox{(L.H.S.)}\eqref{gen_weighted_grad_sob}| \ \sim \
    A^2(k) \ \lesssim \
    A^2(k-1) + B(k)\cdot \big( A(k) + A(k-1)\big) \ . \label{gwsob_induct}
\end{equation}
Notice that with this notation estimate \eqref{weighted_grad_sob}
reads $A(0) \lesssim \ B(0)$. This assumption, together with the
inductive estimate \eqref{gwsob_induct} shows that $A(k) \lesssim
B(k)$ for all $0\leqslant k$. This ends the proof of
\eqref{gen_weighted_grad_sob} and therefore the proof of
\eqref{main_F_L2_est2}.
\end{proof}\ret

\ret
%-------------------------------------------------------------------------
%%%%%%%%%%%%%%%%%%%%%%%%%%%%%%%%%%%%%%%%%%%%%%%%%%%%%%%%%%%%%%%%%%%%%%%%%%
%-------------------------------------------------------------------------

\section{Fixed Time and Space-time Energy Estimates for Complex
Scalar Fields: Conformal Modifications and the Morawetz Theory}\label{phi_L2_sect}

In this section, we prove estimates of the type
\eqref{main_F_L2_est} for the complex scalar field
\eqref{complex_field}. This is necessarily more involved than in
Section \ref{F_L2_section} because the tensor \eqref{phi_em} is not
trace-free  which is a reflection of the fact
that, unlike Maxwell's equations, the scalar wave equation is not
conformally invariant. However, it is well known that there is an
analog of the estimate \eqref{basic_mor_maxwell} for the wave
equation which is also known as the Morawetz estimate, and was first
used in \cite{Mor} to prove local energy decay for solutions to
\eqref{phi_em} in the case $F\equiv 0$. We record this estimate here
as:
\begin{multline}
        \int_{\{t=t_0\}\cap \RR}\ \tp^2|L\phi|^2 +
        \tm^2 |\bL\phi|^2 + (\tp^2+\tm^2)(
        |\snabla \phi|^2 + |\frac{\phi}{r}|^2) \ dx \\
        \lesssim \
        \int_{\{t=0\}\cap\RR}\ (1 + r^2) |\nabla\phi|^2 \ dx \ .
        \label{basic_mor_scalar}
\end{multline}
Notice that we have used the usual derivatives above because we are
assuming that the connection is flat. The usual procedure for
proving \eqref{basic_mor_scalar} involves first coming up with a
certain weighted $L^2$ type identity, again known as the Morawetz
identity, and then performing several integration by parts in order
to ultimately arrive at the estimate \eqref{basic_mor_scalar}. In
the context we are working in, where we wish to prove
energy estimates of the
type \eqref{main_F_L2_est} involving fractional weights as well as
characteristic and space-time energy estimates, such a procedure
would become unduly tedious. Therefore, we provide here a new
approach to Morawetz type estimates which provides
\eqref{basic_mor_scalar} directly in divergence form. Having done this, it
will be straight forward to modify our divergence identity to include
various weights of the type
\eqref{w_gamma_def1}--\eqref{w_gamma_def2} as well as the fractional
Morawetz field \eqref{frac_mor}. Furthermore, our procedure leads to a more
direct geometric insight as to the nature of
\eqref{basic_mor_scalar} which is obscured by the usual integration
by parts proof. \\

Our starting point is to conformally modify the Minkowski metric
$g_{\alpha\beta}$  in such a way that the scalar field equation
\eqref{complex_field} is preserved. As is well known, if one
performs the conformal change of metric:
\begin{equation}
        \widetilde{g} \ = \ \frac{1}{\Omega^2}\, g
    \ , \notag
\end{equation}
for some weight function $\Omega$ on $\mathcal{M}$, then any solution
$\phi$ to the inhomogeneous equation:
\begin{equation}
        \cBox \phi \ = \ G \ , \label{inhomog_cbox}
\end{equation}
will transform to $\psi = \Omega\phi$,
where $\psi$ is a solution to the inhomogeneous conformal scalar
field equation:
\begin{equation}
        \widetilde{\cBox} \psi - \frac{1}{6} \widetilde{R}\psi
    \ = \  \Omega^3 G \ . \label{inhomog_conformal_wave}
\end{equation}
Here $\widetilde{\cBox}$ is the covariant wave equation on
the space $\widetilde{M}\times \mathbb{C}$ with connections
$(\widetilde{\nabla},\widetilde{D})$ and $\widetilde{R}$ the
corresponding scalar curvature. One readily calculates that:
\begin{equation}
        \widetilde{R} \ = \ 6\,  \Omega^3\,
    \nabla^\alpha \nabla_\alpha \big(\frac{1}{\Omega}\big) \ .
    \label{conf_ricci_calc}
\end{equation}
In order that \eqref{inhomog_conformal_wave} match up with
\eqref{inhomog_cbox}, we require $\Omega$ to be such that the scalar
curvature vanishes, $\widetilde{R} \equiv 0$. In light of the
calculation \eqref{conf_ricci_calc}, this will be guaranteed if one
has $\Box (\frac{1}{\Omega})=0$ where $\Box$ is the usual
D'Lambertian on Minkowski space. There are two interesting and
useful choices of $\Omega$ which give a singular solution to this
problem and which are ultimately responsible for the estimate
\eqref{basic_mor_scalar}. These come from the fundamental solutions
to the Laplace and wave equation respectively:
\begin{align}
        {}^I\Omega \ &= \ r , &{}^{I\! I}\Omega \ = \ u\bu \ . \notag
\end{align}
We label the resulting conformal metrics by $\gi$ and $\gii$
respectively. Notice that both of these metrics are singular along
the varieties $r=0$ and $|t| = r$ respectively. However, these
singularities will not effect what we do here because there will
always be extra factors involving positive powers of $\Omega$ in all
the identities we use which will cancel the singularities off. A
striking property of the metrics $\gi$ and $\gii$ is that the
vector-field $K_0$ becomes Killing with respect to \emph{both} of
them (away from the singular set of course). This is seen simply
from the identities:
\begin{align}
        K_0 \Big(\frac{1}{r^2}\Big) \ &= \ -\frac{4t}{r^2} \ ,
        &K_0 \Big(\frac{1}{(u\bu)^2}\Big) \
        &= \ -\frac{4t}{(u\bu)^2} \ . \notag
\end{align}
Thus, for example one has that:
\begin{equation}
    \mathcal{L}_{K_0} \gi \ = \ K_0(\frac{1}{r^2})\, g
    + \frac{1}{r^2} \mathcal{L}_{K_0} g \
    \ = \ 0 , \label{gi_lie}
\end{equation}
with a similar calculation showing that:
\begin{equation}
        \mathcal{L}_{K_0} \gii \ = \ 0 \ . \label{gii_lie}
\end{equation}
Note that the vector-field $T$ is Killing with respect to $\gi$, but
is only conformal Killing with respect to $\gii$. The identities
\eqref{gi_lie}--\eqref{gii_lie} show that if one is to use the
vector-field $K_0$ to produce energy estimates, it is best done with
respect to the metrics $\gi$ and $\gii$ instead of the usual
Minkowski. Accordingly, we define the
\emph{conformal energy-momentum tensors of the first and second
kind} associated to \eqref{inhomog_cbox} to simply be the usual
energy-momentum tensors of the equation
\eqref{inhomog_conformal_wave} with respect to the metrics $\gi$ and
$\gii$ and the appropriately transformed solutions:
\begin{align}
        \cqi_{\alpha\beta}[\phi] \
    &= \  \Re ( D_\alpha (r\phi) \overline {D_\beta (r\phi)}) -
    \frac{1}{2}\gi_{\alpha\beta} \, \widetilde{D}^\gamma(r\phi) \overline{
    D_\gamma(r\phi)} \ , \label{1st_mor_ten}\\
    \cqii_{\alpha\beta}[\phi] \ &= \
    \Re ( D_\alpha (u\bu\phi) \overline {D_\beta (u\bu\phi)}) -
    \frac{1}{2}\gii_{\alpha\beta} \,
    \widetilde{D}^\gamma(u\bu\phi) \overline{
    D_\gamma(u\bu\phi)} \ . \label{2nd_mor_ten}
\end{align}
Here we have used the notation $\widetilde{D}^\gamma =
\widetilde{g}^{\alpha\gamma} D_\alpha$ for $\widetilde{g} =
\gi,\gii$ on lines \eqref{1st_mor_ten} and \eqref{2nd_mor_ten}
respectively. As a direct consequence of the equation
\eqref{inhomog_conformal_wave} and the divergence-fee property of
the energy-momentum tensor for scalar fields one has the divergence
laws:
\begin{align}
        \nablai^\alpha \cqi_{\alpha\beta}[\phi] \ &= \
        r^4 \left( \Re( G \cdot \overline{\frac{1}{r}
        D_\beta(r\phi)}) + F_{\beta\gamma}
        \Im (\phi\overline{\frac{1}{r}D^\gamma (r\phi)}) \right) \ , 
	\label{conf_mor_div_iden1} \\
        \nablaii^\alpha \cqii_{\alpha\beta}[\phi] \ &= \
        (u\bu)^4 \left( \Re( G \cdot \overline{\frac{1}{u\bu}
        D_\beta(u\bu\phi)}) +  F_{\beta\gamma}
        \Im (\phi\overline{\frac{1}{u\bu}D^\gamma (u\bu\phi)})\right) 
	\ , \label{conf_mor_div_iden2}
\end{align}
where $\nablai$ and $\nablaii$ are the Levi-Civita connections of
$\gi$ and $\gii$ respectively. Since the vector-field $K_0$ is
Killing with respect to both of these metrics we may contract it
with the tensors \eqref{1st_mor_ten}--\eqref{2nd_mor_ten} to obtain
momentum densities which satisfy the divergence laws:
\begin{align}
        \nablai^\alpha \left(\cqi_{\alpha\beta}[\phi](K_0)^\beta
        \right)\ &= \
        r^4 \left( \Re( G \cdot \overline{\frac{1}{r}
        D_{K_0}(r\phi)}) + (K_0)^\beta F_{\beta\gamma}
        \Im (\phi\overline{D^\gamma \phi}) \right) \ ,
        \label{mor_density1}\\
        \nablaii^\alpha
        \left(\cqii_{\alpha\beta}[\phi](K_0)^\beta\right) \ &= \
        (  u\bu)^4 \left( \Re( G \cdot \overline{\frac{1}{u\bu}
        D_{K_0}(u\bu\phi)}) +  (K_0)^\beta F_{\beta\gamma}
        \Im (\phi\overline{
        \frac{1}{u\bu}D^\gamma(u\bu \phi)})\right) \ . \label{mor_density2}
\end{align}
These last two identities can now be integrated over various
space-time regions to obtain positive  quantities for the scalar
field $\phi$ at the cost of estimating the source terms on the right
hand side of \eqref{mor_density1}--\eqref{mor_density2}. To
calculate these, notice that the volume forms of the metrics $\gi$
and $\gii$ are:
\begin{align}
        dV_I \ &= \ \frac{1}{r^4} dV_\mathcal{M} ,
    &dV_{I\! I} \ &= \ \frac{1}{(u\bu)^4} dV_\mathcal{M}
    \ , \label{volume_forms}
\end{align}
while $r T$ and $u\bu T$ are the respective Lorentzian unit normals
to the time slices $t=const.$ Furthermore, notice that the
vector-fields $r L $ and $u\bu L$ are the respective Lorentzian unit
normal to the cones $u=const.$ Therefore, applying the geometric
Stokes theorem to with respect to these quantities, we arrive at the
basic Morawetz type energy estimates for the complex scalar field
\eqref{inhomog_cbox}:\\

\begin{prop}[First and Second Morawetz Estimates for Complex Scalar
Fields]\label{basic_scalar_mor_prop} Let  $\Omega= r,u\bu$. Then one
has the following estimates for solutions to the inhomogeneous
equation \eqref{inhomog_cbox}:
\begin{multline}
       \int_{ \{t=t_0\}\cap\RR } \frac{1}{4}\left(
       \bu^2|\frac{1}{\Omega} D_L(\Omega\phi)|^2 +
       u^2|\frac{1}{\Omega}D_{\bL}(\Omega\phi)|^2 +
       (\bu^2 + u^2) |\sD \phi|^2 \right) dx \\
       + \ \sup_u\
       \int_{ C(u)\cap \{0 < t < t_0\}}
       \frac{1}{4}\left(
       \bu^2|\frac{1}{\Omega} D_L(\Omega\phi)|^2 +
       u^2 |\sD \phi|^2 \right) dV_{C(u)} \ \\
       \leqslant \ \ \ \
       \int_0^{t_0}\int_{\RR} | G\cdot \frac{1}{\Omega}
       D_{K_0}(\Omega\phi)|
       + |(K_0)^\beta F_{\beta\gamma}
       \Im (\phi\overline{\frac{1}{\Omega}D^\gamma (\Omega\phi)})| \ dxdt \\
       + \int_{\{t=0\}\cap\RR } \frac{r^2}{2}
       |\frac{1}{\Omega}D(\Omega\phi)|^2 \ dx \ . \label{basic_scalar_mor_est}
\end{multline}
\end{prop}\ret

To derive \eqref{basic_mor_scalar} from estimate
\eqref{basic_scalar_mor_est} above is a simple matter. Assume now
that both $G=F\equiv 0$. Notice first that it suffices to prove
\eqref{basic_mor_scalar} with the $\tp,\tm$ weights replaced by
$\bu,u$ and $(1+r^2)$ on the right hand side replaced by $r^2$. This
follows at once from expanding the usual energy identity in a
null frame, and using the fixed time Poincare inequality:
\begin{equation}
        \int_{\{t=const.\}\cap\RR } \
        |\frac{\phi}{r}|^2 \ dx \ \lesssim \
        \int_{\{t=const.\}\cap\RR } \
        |\partial_r\phi|^2 \ dx \ , \label{basic_poincare}
\end{equation}
to deal with the non-differentiated term. To deal with the homogeneous
version of estimate \eqref{basic_mor_scalar}, we simply use the identities:
\begin{align}
        \bu^2|\frac{1}{r} L(r\phi)|^2 \ &= \ |\bu L\phi + \frac{\bu \phi}{r}|^2 \ ,
        &u^2 |\frac{1}{r} \bL(r\phi)|^2 \ &= \ |u\bL\phi -
        \frac{u\phi}{r}|^2 \ , \notag \\
        \bu^2|\frac{1}{u\bu} L(u\bu\phi)|^2 \ &= \ |\bu L\phi +
        2\phi|^2 \ , &u^2|\frac{1}{u\bu} \bL(u\bu\phi)|^2 \ &= \ |u \bL\phi +
        2\phi|^2 \ . \notag
\end{align}
together with:
\begin{align}
        \frac{-u\phi}{r} \ &= \ \big(\bu L\phi +
        2\phi\big) - \big(\bu L\phi + \frac{\bu
        \phi}{r}\big) \ , \label{aux_mor_id1}\\
        \frac{\bu\phi}{r} \ &= \ \big(u \bL\phi +
        2\phi\big) - \big(u\bL\phi -
        \frac{u\phi}{r}\big) \ , \label{aux_mor_id2}
\end{align}
which collectively imply the two pointwise estimates:
\begin{align}
        \bu^2|L\phi|^2 + u^2|\bL\phi|^2  \ &\leqslant \
    2\, \big ( \bu^2 |\frac{1}{r} L(r \phi)|^2 +
        u^2 |\frac{1}{r} \bL(r \phi)|^2 \big) +
    2 \, (\bu^2 + u^2)|\frac{\phi}{r}|^2 \ , \notag\\
    (\bu^2 + u^2)|\frac{\phi}{r}|^2 \ &\leqslant  \
    2 \ \sum_{\Omega = r,u
        \bu}\ \bu^2 |\frac{1}{\Omega} L(\Omega \phi)|^2 +
        u^2 |\frac{1}{\Omega} \bL(\Omega \phi)|^2 \ . \notag
\end{align}
Combining these last two, we have completed
our demonstration of \eqref{basic_mor_scalar}.\\

We now prove a fractionally weighted and space-time generalization
of \eqref{basic_scalar_mor_est} which will be our analog of estimate
\eqref{main_F_L2_est} for the complex scalar field
\eqref{inhomog_cbox}. This is:\\

\begin{prop}[Generalized Morawetz Estimate for Complex Scalar Fields]
\label{gen_mor_prop} Let $\phi$ be a complex scalar field which
satisfies the equation \eqref{inhomog_cbox}, and let $F$ be the
curvature of the corresponding connection. Let $0< \epsilon ,
\gamma$ be given parameters, and let $s$ be chosen so that
$\frac{1}{2} \leqslant s \leqslant 1$ as well as
$\epsilon \leqslant s-\frac{1}{2}$. Define the weights
$\tz,w_\gamma,w_{\gamma,\epsilon}$, and $w'_\gamma$ as on lines
\eqref{tz_def}--\eqref{w_gamma_def2}. Also, define the generalized
Morawetz type energy content of $\phi$ in the time slab
$\{0\leqslant t \leqslant t_0\}\cap\RR$:
\begin{align}
        &E_0^{(\epsilon,\gamma,s)}(0,t_0)[\phi] \ ,   \label{phi_slab_eng}\\
        = \ &\sup_{0\leqslant t \leqslant t_0}\ \int_{ \{t\}\cap\RR } \left(
        \tp^{2s}|\frac{1}{r} D_L(r\phi)|^2 +
        \tm^{2s} |D_{\bL}\phi|^2 +
        \tp^{2s} \left( |\sD \phi|^2 + |\frac{\phi}{r}|^2 \right)\right)
        w_\gamma  \ dx \ , \notag \\ + \
        &\sup_{u}\ \int_{ C(u)\cap \{0 \leqslant t \leqslant t_0\}}
        \left(
        \tp^{2s}|\frac{1}{r} D_L(r\phi)|^2
        + \tm^{2s} |\sD \phi|^2 + \tp^{2s}|\frac{u\phi}{\bu r}|^2
        \right)w_\gamma \  dV_{C(u)} \ , \notag \\
        + \ &\int_{0}^{t_0} \int_{ \RR } \ \left(
        \tp^{2s} |\frac{1}{r} D_L(r\phi)|^2
        + \tz^{1+2\epsilon} \left( \tm^{2s} |D_{\bL}\phi|^2 +
        \tp^{2s} |\sD \phi|^2 + \tp^{2s} |\frac{\phi}{r}|^2 \right)
        \right)w'_{\gamma,\epsilon}  \ dxdt \ . \notag
\end{align}
Then one has the following general weighted energy type inequality:
\begin{multline}
        E_0^{(\epsilon,\gamma,s)}(0,t_0)[\phi] \ \leqslant \
        C^2_{\epsilon,\gamma} \Big[\llp{F}
    {L^\infty[0,t_0](\epsilon)}^2\cdot
        E_0^{(\epsilon,\gamma,s)}(0,t_0)[\phi] \\
        \ + \ \int_{0}^{t_0}
        \int_{\RR } \tp^{2s} \tm \ | G|^2 \ w_{\gamma,\epsilon} \ dxdt  \
        + \ \int_{ \{t=0\}\cap\RR } (1 + r^2)^{s + \gamma}
        |D(\phi)|^2 \ dx \ \Big] \ , \label{main_phi_L2_est}
\end{multline}
where we have set $\llp{F}{L^\infty[0,t_0](\epsilon)}$ equal to the
time-slab $L^\infty$ type norm:
\begin{equation}
        \llp{F}{L^\infty[0,t_0](\epsilon)} \ = \
    \sup_{\substack{ 0 \leqslant \ t \
    \leqslant t_0\\ x\in \{t\}\times\RR}}
    \ \left( \tp^{2}\tz^{-2\epsilon}|\alpha| \, + \,
    \tp^\frac{3}{2} \tm^\frac{1}{2} \tz^{-\epsilon} |\rho|
    \, + \, \tp\tm |\balpha| \right)\cdot w_{\gamma,\epsilon}
    (w_\gamma)^{-1} \ . \label{F_Linfty_sharp_norm}
\end{equation}
In particular, if we further define the multiindexed
$L^\infty$ norm for electro-magnetic fields:
\begin{multline}
    \llp{F}{L^\infty[0,t_0](s,\gamma,\epsilon)}^2 \ = \
    \sup_{\substack{ 0 \leqslant \ t \
    \leqslant t_0\\ x\in \{t\}\times\RR}}
    \ \left( \tp^{2s+3}|\alpha|^2 + \tp^2\tm^{2s+1}|\balpha|^2 +
    \tp^{2s+2}\tm (\td{\rho}^2 + \sigma^2)\right)\cdot w_\gamma \\
    + \ |q(F)|^2 \  + \ \lp{\tp^{s+1}\tm^\frac{1}{2}
    (w')^\frac{1}{2}_{\gamma,\epsilon} \alpha}{L^2(L^\infty)[0,t_0]}^2
    \ , \label{F_Linfty_norm}
\end{multline}
where the quantity $q(F)$ is defined as on line \eqref{charge_def},
then under the more restrictive condition
that $4\epsilon \leqslant s-\frac{1}{2}$ we have the bound:
\begin{equation}
        \llp{F}{L^\infty[0,t_0](\epsilon)} \ \leqslant \
    \llp{F}{L^\infty[0,t_0](s,\gamma,\epsilon)} \ .
    \label{linfty_to_linfty_bound}
\end{equation}
In particular, the estimate \eqref{main_phi_L2_est} remains valid
with \eqref{F_Linfty_sharp_norm} replaced by \eqref{F_Linfty_norm}.
\end{prop}\ret

\begin{rem}\label{analogy_rem}
Comparing estimate \eqref{main_F_L2_est} with
\eqref{main_phi_L2_est} we see that there is a close analogy, at
least with respect to energy estimates, between the various null
components of $F$ and the gradient $D\phi$. Schematically, these
are:
\begin{align}
        \alpha \ &\sim \ \frac{1}{r}D_L(r\phi) \ , &\balpha \ &\sim \
        D _{\bL}\phi \ , \label{analogy1}\\
        \sigma \ &\sim \ \sD\phi \ , &\td{\rho} \ &\sim \
        \frac{\phi}{r} \ . \label{analogy2}
\end{align}
As we shall see in the following two sections, this analogy persists
in the discussion of $L^\infty$ type estimates. We will make solid
of this in the sequel, where we will use the analogy
\eqref{analogy1}--\eqref{analogy2} to reduce the field equations
\eqref{basic_MKG} and its commutators to a equation for a single
abstract vector quantity $\Psi$ which has weighted energy and
peeling properties equal to the null decompositions of $F$ and
$D\phi$.
\end{rem}\ret

\begin{rem}
Note that the $L^\infty$ norm \eqref{F_Linfty_norm} contains
several terms not present in the original $L^\infty$ norm
\eqref{F_Linfty_sharp_norm}. Since the more involved norm
will be used exclusively in the sequel, we have stated it
here for future reference. Also, notice that on line
\eqref{linfty_to_linfty_bound} we have used the formulas
\eqref{tF_def} and \eqref{null_charge_tensor}
for the subtracted charge $\td{\rho}$. Note that this bound
(for the charge) becomes sharp in the far exterior $r-t\sim r$.
\end{rem}\ret

\begin{proof}[Proof of estimate \eqref{main_phi_L2_est}]
The proof is nearly identical to that of \eqref{main_F_L2_est} in
the previous section. Our first step is to come up with a
fractionally weighted momentum density associated with the first
conformal metric $\gi$. We define this to be:
\begin{equation}
        {}^I\td{P}_\alpha^{(s,\gamma,\epsilon)}[\phi] \ = \
        \cqi_{\alpha\beta}[\phi](\overline{K}_0^s)^\beta\cdot
        \td{w}_{\gamma,\epsilon} \ , \label{phi_mom_density}
\end{equation}
where $\td{w}_{\gamma,\epsilon}$ is the weight function from line
\eqref{w_gam_eps_weight}. This satisfies the divergence law:
\begin{multline}
        {}^I\td{\nabla}^\alpha  {}^I\td{P}_\alpha^{(s,\gamma,\epsilon)}[\phi] \ = \
        \ r^4 \left( \Re( G \cdot \overline{\frac{1}{r}
        D_{\overline{K}^s_0}(r\phi)}) + (\overline{K}^s_0)^\beta F_{\beta\gamma}
        \Im (\phi\overline{D^\gamma \phi}) \right)\cdot \td{w}_{\gamma,\epsilon}\\
        + \ \frac{1}{2} \ \cqi_{\alpha\beta}[\phi]
        \, {}^{(\overline{K}_0^s)}\td{\pi}^{\alpha\beta}\cdot \td{w}_{\gamma,\epsilon} \\
        - \ \frac{r^2}{2}\ \cqi_{\bL \beta}[\phi] (\overline{K}_0^s)^\beta \cdot
        L(\td{w}_{\gamma,\epsilon}) \
    - \ \frac{r^2}{2}\ \cqi_{L\beta}[\phi]
        (\overline{K}_0^s)^\beta \cdot \bL(\td{w}_{\gamma,\epsilon}) \ .
        \label{main_phi_div_id}
\end{multline}
Integrating both sides of the above expression with respect to the
volume $dV_I \ = \ r^{-4} dV_\mathcal{M}$ over regions of the form:
\begin{equation}
        \mathcal{R}(t_0,u_0) \ = \ \{0\leqslant t \leqslant
        t_0\}\cap\{u\leqslant u_0\} \ , \notag
\end{equation}
and noting that $rT$ and $rL$ are the respective Lorentzian unit
normals to the time slices $t=const.$ and cones $u=const.$, we
arrive at the following integral identity:
\begin{multline}
        \int_{ \{t=0\}\cap\{-u_0 \leqslant r \}} \ r^{-2}\
        {}^I\td{P}_0^{(s,\gamma,\epsilon)}[\phi]\ dx \ = \
        \int_{\{t=t_0\}\cap\{t_0 - u_0 \leqslant r\}} \ r^{-2}\
        {}^I\td{P}_0^{(s,\gamma,\epsilon)}[\phi]\ dx \\
        + \ \int_{C(u_0)\cap\{0\leqslant t \leqslant t_0\}}\ r^{-2}\
        {}^I\td{P}_L^{s,\gamma,\epsilon}[\phi]\ dV_{C(u_0)} \ + \
        \int\int_{\mathcal{R}(t_0,u_0)} \
        r^{-4}\ \hbox{(R.H.S.)}\eqref{main_phi_div_id}\ dxdt \ .
        \label{main_phi_int_id}
\end{multline}
We now calculate each term in \eqref{main_phi_div_id} individually.
We will use the same $\sim$ notation as in the proof of
\eqref{main_F_L2_est}. At time $t=0$, we compute that:
\begin{equation}
        r^{-2}\ {}^I\td{P}_0^{(s,\gamma,\epsilon)}[\phi] \ \sim \
        (1+r^2)^{s+\gamma} |\frac{1}{r} D(r\phi)|^2 \ .
        \label{Pphi_0_est}
\end{equation}
Likewise, at time $t=t_0$ we have that:
\begin{equation}
        r^{-2}\ {}^I\td{P}_0^{(s,\gamma,\epsilon)}[\phi] \ \sim \
        \left(\tp^{2s}|\frac{1}{r}D_L(r\phi)|^2 +
        \tm^{2s}|\frac{1}{r}D_{\bL}(r\phi)|^2 +
        \tp^{2s} |\sD\phi|^2\right)\cdot w_\gamma \ , \label{Pphi_t0_est}
\end{equation}
and on the cone $u_0=const.$ we have:
\begin{equation}
        r^{-2}\ {}^I\td{P}_L^{(s,\gamma,\epsilon)}[\phi] \ \sim \
        \left(\tp^{2s}|\frac{1}{r}D_L(r\phi)|^2 + \tm^{2s}|\sD\phi|^2
        \right)\cdot w_\gamma \ .
        \label{Pphi_u0_est}
\end{equation}\\

It remains to calculate the expression $r^{-4}
\hbox{(R.H.S.)}\eqref{main_phi_div_id}$. We do this for each term
separately. Since the first such term does not have a sign, we put
absolute value signs around it and estimate:
\begin{align}
        \big| \Re( G \cdot \overline{\frac{1}{r}
        D_{\overline{K}^s_0}(r\phi)}) \big| \ &\lesssim \
        |G|\cdot\left(\tp^{2s} |\frac{1}{r} D_L(r\phi)| + \tm^{2s}|
        D_{\bL}\phi| + \tm^{2s}|\frac{\phi}{r}| \right)
        \ , \label{Q_phi_div_exp1}\\
        \big|(\overline{K}^s_0)^\beta F_{\beta\gamma}
        \Im (\phi\overline{D^\gamma \phi})\big| \ &\lesssim \
        \tp^{2s+1}\left(  |\alpha|\cdot|\frac{\phi}{r}|\cdot|\sD\phi| +
        |\rho|\cdot|\frac{\phi}{r}|\cdot|\frac{1}{r}D_L(r\phi)|\right)
        \label{Q_phi_div_exp2}\\
        &\ \ \ \ \ \ + \ \tp\tm^{2s}\left(|\balpha|\cdot|\frac{\phi}{r}|\cdot|\sD\phi|
        + |\rho|\cdot|\frac{\phi}{r}|\cdot |D_{\bL} \phi| \right) \
        . \notag
\end{align}
To calculate the second term on the right hand side of expression
\eqref{main_phi_div_id} we need to compute the deformation tensor
${}^{(\overline{K}_0^s)}\td{\pi}_{\alpha\beta}$. Using the
identities \eqref{lie_def_tensor} and \eqref{K0s_def_tensor} this
is:
\begin{align}
        {}^{(\overline{K}_0^s)}\td{\pi} \ &= \
        \overline{K}_0^s(\frac{1}{r^2})\, g + \frac{1}{r^2}\,
        {}^{(\overline{K}_0^s)}\pi \ , \notag\\
        &= \ \frac{2}{r^2}\left(\frac{\bv-v}{r} - (\dot{\bv} + \dot{v})\right)
        (\theta^{\bL} \otimes \theta^L +
        \theta^L\otimes\theta^{\bL}) \ . \notag
\end{align}
Contracting this last line with $\frac{1}{2}r^{-4}\cqi[\phi]\cdot
\td{w}_{\gamma,\epsilon}$ yields:
\begin{equation}
        \frac{1}{2}r^{-4}\, \cqi_{\alpha\beta}[\phi]\,
        {}^{(\overline{K}_0^s)}\td{\pi}^{\alpha\beta}
        \cdot \td{w}_{\gamma,\epsilon} \ \sim \
        \left( \frac{\bv-v}{r} -
        (\dot{\bv} + \dot{v})\right)
        \cdot|\sD\phi|^2\cdot w_\gamma
        \ . \label{Q_phi_div_exp3}
\end{equation}
In particular, we see that this term only contributes a positive
addition to the right hand side of \eqref{main_phi_div_id}. To
compute the last two terms of $r^{-4}\hbox{(R.H.S.)}
\eqref{main_phi_div_id}$, we simply use the calculations
\eqref{w_gam_eps_bound2}--\eqref{w_gam_eps_bound3} and the
expansions \eqref{phi_em_decomp} to conclude that:
\begin{align}
        -\ \frac{1}{2}r^{-2}\ \cqi_{\bL \beta}[\phi]
        (\overline{K}_0^s)^\beta \cdot
        L(\td{w}_{\gamma,\epsilon})\ &\sim \ C^{-1}_{\gamma,\epsilon}\
    \tz^{1+\epsilon}\left(\tm^{2s}
        |\frac{1}{r} D_{\bL}(r\phi)|^2 +
        \tp^{2s}|\sD\phi|^2 \right)\cdot w'_{\gamma,\epsilon} \ ,
        \label{Q_phi_div_exp4}\\
        -\ \frac{1}{2}r^{-2}\ \cqi_{L \beta}[\phi]
        (\overline{K}_0^s)^\beta \cdot
        \bL(\td{w}_{\gamma,\epsilon})\ &\sim \
        C^{-1}_{\gamma,\epsilon}\ \left(\tp^{2s}|\frac{1}{r}D_L(r\phi)|^2 +
        \tm^{2s}|\sD\phi|^2\right)\cdot w'_{\gamma,\epsilon} \ .
        \label{Q_phi_div_exp5}
\end{align}\ret

Now, inserting the terms
\eqref{Pphi_0_est}--\eqref{Q_phi_div_exp5}
in the identity \eqref{main_phi_int_id} and removing the positive
contribution \eqref{Q_phi_div_exp3} form the right hand side of the
resulting estimate, we arrive at our first preliminary estimate to
\eqref{main_phi_L2_est}:
\begin{multline}
        \int_{\{t=0\}\cap\RR}\ (1+r^2)^{s+\gamma} |\frac{1}{r}
        D(r\phi)|^2 \ dx \ + \
        \int_{0}^{t_0}\int_{\RR}\ |G|\cdot\left(\tp^{2s}
    |\frac{1}{r} D_L(r\phi)| + \tm^{2s}|
        D_{\bL}\phi| + \tm^{2s}|\frac{\phi}{r}| \right)\cdot w_\gamma \ dxdt \\
        \int_0^{t_0} \int_{\RR} \ \tp^{2s+1}\left(
        |\alpha|\cdot|\frac{\phi}{r} |\cdot|\sD\phi| +
        |\rho|\cdot|\frac{\phi}{r}|\cdot|\frac{1}{r}D_L(r\phi)|\right)\cdot
        w_\gamma \ dxdt\\
        + \ \int_0^{t_0}\int_{\RR}\
        \tp\tm^{2s}\left(|\balpha|\cdot|\frac{\phi}{r}|\cdot|\sD\phi|
        + |\rho|\cdot|\frac{\phi}{r}|\cdot |D_{\bL} \phi| \right)
    \cdot w_\gamma \ dxdt\\
    \geqslant \ C^{-1}_{\gamma,\epsilon}\ \Bigg[ \sup_{0\leqslant t
        \leqslant t_0}\ \int_{\{t\}\cap\RR}\ \left(\tp^{2s}|\frac{1}{r}D_L(r\phi)|^2 +
        \tm^{2s}|\frac{1}{r}D_{\bL}(r\phi)|^2 +
        \tp^{2s} |\sD\phi|^2\right)\cdot w_\gamma  \ dx \\
    + \ \sup_{u}\ \int_{C(u)\cap \{0\leqslant t \leqslant t_0\}}\
    \left(\tp^{2s}|\frac{1}{r}D_L(r\phi)|^2 + \tm^{2s}|\sD\phi|^2
        \right)\cdot w_\gamma \ dV_{C(u)}\\
    + \ \int_0^{t_0} \int_{\RR} \
        \left(\tp^{2s}|\frac{1}{r}D_L(r\phi)|^2 + \tz^{1+\epsilon}\left(\tm^{2s}
        |\frac{1}{r} D_{\bL}(r\phi)|^2 + \tp^{2s}|\sD\phi|^2
        \right)\right)\cdot w'_{\gamma,\epsilon}
    \ dxdt \ \Bigg] \ . \label{prilim_main_phi_est1}
\end{multline}\ret

By repeating the above steps with the density \eqref{phi_mom_density}
replaced by the following weighted version of the usual
energy density for $\phi$:
\begin{equation}
        P^{(\gamma,\epsilon)}_\alpha[\phi] \ = \
        Q_{\alpha\beta}[\phi](T)^\beta\cdot
    \td{w}_{\gamma,\epsilon} \ , \notag
\end{equation}
and using the bounds:
\begin{align}
        \big| \Re( G \cdot \overline{
        D_{T}\phi} ) \big| \ &\lesssim \
    \hbox{(L.H.S)}\eqref{Q_phi_div_exp1} \ , \notag\\
    |F_{0\gamma}\Im (\phi\overline{D^\gamma \phi})\big| \ &\lesssim \
    \hbox{(L.H.S)}\eqref{Q_phi_div_exp2}
\end{align}
we also have the estimate:
\begin{multline}
        \hbox{(R.H.S.)}\eqref{prilim_main_phi_est1}\
    \geqslant \ C^{-1}_{\gamma,\epsilon}\ \Bigg[ \sup_{0\leqslant t
        \leqslant t_0}\ \int_{\{t\}\cap\RR}\ |D_{\bL}\phi|^2 \cdot w_\gamma\ dx \\
    + \ \sup_{u}\ \int_{C(u)\cap \{0\leqslant t \leqslant t_0\}}\
    |D_L\phi|^2 \cdot w_\gamma \ dV_{C(u)}\
    + \ \int_0^{t_0} \int_{\RR} \
        \tz^{1+\epsilon}\, |D_{\bL}\phi|^2 \cdot w'_{\gamma,\epsilon}
    \ dxdt \ \Bigg] \ . \label{prilim_main_phi_est1.5}
\end{multline}\ret

By adding together estimates
\eqref{prilim_main_phi_est1}--\eqref{prilim_main_phi_est1.5}, we
have almost achieved the  statement of estimate
\eqref{main_phi_L2_est}. What is missing from the right hand side
are terms involving the weighted quantity $\bu^{2s}\, \big|
\frac{\phi}{r}\big|$. As with estimate \eqref{basic_scalar_mor_est},
bounds on this quantity will come from integrating a divergence
identity involving the tensor \eqref{2nd_mor_ten} and then combining
the resulting estimate with \eqref{prilim_main_phi_est1} above. To
do this, we form the second weighted momentum density:
\begin{equation}
        {}^{I\! I}\td{P}_\alpha^{(s,\gamma,\epsilon)}[\phi] \ = \
        \cqii_{\alpha\beta}[\phi](K_0)^\beta\cdot
        \td{w}_{s,\gamma,\epsilon}\ , \label{phi_mom_density2}
\end{equation}
where the weight function $\td{w}_{s,\gamma,\epsilon}$ is given by:
\begin{multline}
        \td{w}_{s,\gamma,\epsilon} \ = \ (1 + \bu)^{2s-2}\cdot\left(
    (2-u)^{2\gamma}\chi^+(-u) \ + \ 1-\chi^+(-u)\right)\\
    + \ (1 + \bu)^{2s-2-2\epsilon}\cdot\left(
    (2-u)^{2\gamma+2\epsilon}\chi^+(-u) \ + \ 1-\chi^+(-u)\right)
    \label{w_s_gam_eps_weight}
\end{multline}
We now follow through the same steps used to prove
\eqref{prilim_main_phi_est1} above. Computing the divergence of
\eqref{phi_mom_density2} we have that:
\begin{multline}
        {}^{I\! I}\td{\nabla}^\alpha  {}^{I\! I}\td{P}_\alpha^{(s,
    \gamma,\epsilon)}[\phi] \ = \
        \ (u\bu)^4 \left( \Re( G \cdot \overline{\frac{1}{u\bu}
        D_{K_0}(u\bu\phi)}) + (K_0)^\beta F_{\beta\gamma}
        \Im (\phi\overline{D^\gamma \phi}) \right)\cdot \td{w}_{s,\gamma,\epsilon}\\
        - \ \frac{(u\bu)^2}{2}\ \cqii_{\bL \beta}[\phi] (K_0)^\beta \cdot
        L(\td{w}_{s,\gamma,\epsilon}) \
    - \ \frac{(u\bu)^2}{2}\ \cqii_{L\beta}[\phi]
        (K_0)^\beta \cdot \bL(\td{w}_{s,\gamma,\epsilon}) \ .
        \label{main_phi_div_id2}
\end{multline}
Integrating this expression with respect to the volume $dV_{I\! I} =
(u\bu)^{-4} dV_\mathcal{M}$ and using the Stokes theorem we arrive
at the integral identity:
\begin{multline}
        \int_{ \{t=0\}\cap\{-u_0 \leqslant r \}} \ (u\bu)^{-2}\
        {}^{I\! I}\td{P}_0^{(s,\gamma,\epsilon)}[\phi]\ dx \ = \
        \int_{\{t=t_0\}\cap\{t_0 - u_0 \leqslant r\}} \ (u\bu)^{-2}\
        {}^{I\! I}\td{P}_0^{(s,\gamma,\epsilon)}[\phi]\ dx \\
        + \ \int_{C(u_0)\cap\{0\leqslant t \leqslant t_0\}}\ (u\bu)^{-2}\
        {}^{I\! I}\td{P}_L^{s,\gamma,\epsilon}[\phi]\ dV_{C(u_0)} \ + \
        \int\int_{\mathcal{R}(t_0,u_0)} \
        (u\bu)^{-4}\ \hbox{(R.H.S.)}\eqref{main_phi_div_id2}\ dxdt \ .
        \label{main_phi_int_id2}
\end{multline}
We now compute in turn each term in \eqref{main_phi_int_id2} above.
When $t=0$ we have that:
\begin{equation}
        (u\bu)^{-2}\ {}^{I\! I}\td{P}_0^{(s,\gamma,\epsilon)}[\phi]
        \ \sim \ r^{2(s + \gamma)} |\frac{1}{r^2} D(r^2 \phi)|^2
    \ . \label{Pphi_0_est2}
\end{equation}
Likewise, at time $t=t_0$ and along the cone $u=const.$ we compute
that:
\begin{align}
        (u\bu)^{-2}\ {}^{I\! I}\td{P}_0^{(s,\gamma,\epsilon)}[\phi]
        \ &\sim \ \left( \tp^{2s-2}\bu^2\, | \frac{1}{u\bu}D_L(u\bu\phi)|^2
    + \tp^{2s-2}u^2\, |\frac{1}{u\bu}D_{\bL}(u\bu\phi) |^2 +
    \tp^{2s-2}\bu^2\, |\sD \phi|^2 \right)\cdot w_\gamma \ ,
        \label{Pphi_t0_est2}\\
    (u\bu)^{-2}\ {}^{I\! I}\td{P}_L^{(s,\gamma,\epsilon)}[\phi]
        \ &\sim \ \left( \tp^{2s-2}\bu^2\, | \frac{1}{u\bu}D_L(u\bu\phi)|^2 +
    \tp^{2s-2}u^2\, |\sD\phi|^2 \right)\cdot w_\gamma
    \ . \label{Pphi_L_est2}
\end{align}
We now turn our attention to the terms in the expression
$(u\bu)^{-4}\ \hbox{(R.H.S.)}\eqref{main_phi_div_id2}$. First, using  the identities
\eqref{aux_mor_id1}--\eqref{aux_mor_id2} with $L,\bL$ replaced
by $D_L, D_{\bL}$ we conclude that:
\begin{align}
         |\bu D_L \phi + 2\phi| \ &\leqslant \ \bu |\frac{1}{r} D_L(r\phi)|
         + u |\frac{\phi}{r}| \ , \notag \\
     |uD_{\bL}\phi + 2\phi| \ &\leqslant \ u|D_{\bL}\phi| +
         2\bu|\frac{\phi}{r}| \ . \notag
\end{align}
This then implies that one may estimate:
\begin{align}
        \Big|\Re( G \cdot \overline{\frac{1}{u\bu}
        D_{K_0}(u\bu\phi)})\Big|\cdot \td{w}_{s,\gamma,\epsilon}
    \ &\lesssim \ |G|\cdot\left( \tp^{2s-1} |\bu D_L \phi + 2\phi|
    + \tp^{2s-2} u  |uD_{\bL} \phi + 2\phi|
    \right)\cdot w_\gamma  \ , \notag \\
    &\lesssim \  |G|\cdot\left(\tp^{2s}|\frac{1}{r}D_L(r\phi)|
    + \tm^{2s}|D_{\bL}\phi| + \tp^{2s}|\frac{u\phi}{\bu r}|
    \right)\cdot w_\gamma \ . \label{Q_phi_div_exp2.1}
\end{align}
To bound the second term on the right hand side of \eqref{main_phi_div_id2} above,
we simply use the estimate \eqref{Q_phi_div_exp2} noting that one may trade
$\tp^{2s-2}\tm^2 \lesssim \tm^{2s}$:
\begin{align}
        \Big|(K_0)^\beta F_{\beta\gamma}\Im (\phi\overline{D^\gamma \phi})
    \Big|\cdot \td{w}_{s,\gamma,\epsilon} \ \lesssim \
    &\tp^{2s+1}\left(  |\alpha|\cdot|\frac{\phi}{r}|\cdot|\sD\phi| +
        |\rho|\cdot|\frac{\phi}{r}|\cdot|\frac{1}{r}D_L(r\phi)|\right)
    \cdot w_\gamma \label{Q_phi_div_exp2.2} \ , \\
    &\ \ \ \
    + \ \ \
    \ \tp\tm^{2s}\left(|\balpha|\cdot|\frac{\phi}{r}|\cdot|\sD\phi|
        + |\rho|\cdot|\frac{\phi}{r}|\cdot |D_{\bL} \phi| \right)
    \cdot w_\gamma  \ . \notag
\end{align}
To compute the second two terms on the right hand side of
\eqref{main_phi_div_id2}, first notice that a simple calculation
similar to that used in lines
\eqref{w_gam_eps_bound2}--\eqref{w_gam_eps_bound3}
above, which we will omit, shows that:
\begin{align}
        -\frac{1}{2} L(\td{w}_{\gamma,\epsilon}) \ &\geqslant \
    C^{-1}_{\gamma,\epsilon} \ \tp^{2s-2} \tz^{1 + 2\epsilon}
    w'_{\gamma,\epsilon}\ , \notag \\
    -\frac{1}{2} \bL(\td{w}_{\gamma,\epsilon}) \ &\geqslant \ 0 \
    . \notag
\end{align}
Notice that we only care that the second term above is non-negative as
this will be its only use in our proof of \eqref{main_phi_L2_est}.
In particular, the above two lines allow us to conclude that:
\begin{align}
        - \ \frac{1}{2(u\bu)^2}\ \cqii_{\bL \beta}[\phi] (K_0)^\beta \cdot
        L(\td{w}_{\gamma,\epsilon}) \ &\geqslant \ C^{-1}_{\gamma,\epsilon}
    \tp^{2s-2}\tz^{1+2\epsilon}\left(u^2\,
        |\frac{1}{u\bu}D_{\bL}(u\bu\phi)|^2 + \bu^2\, |\sD \phi|^2
    \right)\cdot w'_{\gamma,\epsilon} \ , \label{Q_phi_div_exp2.3} \\
    -\ \frac{1}{2(u\bu)^2}\ \cqii_{L\beta}[\phi]
        (K_0)^\beta \cdot \bL(\td{w}_{\gamma,\epsilon}) \ &\geqslant \
        0 \ . \label{Q_phi_div_exp2.4}
\end{align}
Combining the estimates
\eqref{Pphi_0_est2}--\eqref{Q_phi_div_exp2.4} into the identity
\eqref{main_phi_int_id2} and excluding all the terms on the right
hand side of the resulting inequality except those which involve the
$D_{\bL}$ derivative and the characteristic term involving the $D_L$
derivative, we arrive at our compliment to estimate
\eqref{prilim_main_phi_est1}:
\begin{multline}
        \int_{\{t=0\}\cap\RR}\ r^{2(s+\gamma)} |\frac{1}{r^2}
        D(r^2 \phi)|^2 \ dx \\
        + \ \int_{0}^{t_0}\int_{\RR}\ |G|\cdot\left(\tp^{2s}
    |\frac{1}{r} D_L(r\phi)| + \tm^{2s}|
        D_{\bL}\phi| + \tp^{2s}|\frac{u\phi}{\bu r}| \right)\cdot w_\gamma \ dxdt \\
        + \ \int_0^{t_0} \int_{\RR} \ \tp^{2s+1}\left(
        |\alpha|\cdot|\frac{\phi}{r} |\cdot|\sD\phi| +
        |\rho|\cdot|\frac{\phi}{r}|\cdot|\frac{1}{r}D_L(r\phi)|\right)\cdot
        w_\gamma \ dxdt\\
        + \ \int_0^{t_0}\int_{\RR}\
        \tp\tm^{2s}\left(|\balpha|\cdot|\frac{\phi}{r}|\cdot|\sD\phi|
        + |\rho|\cdot|\frac{\phi}{r}|\cdot |D_{\bL} \phi| \right)
    \cdot w_\gamma \ dxdt\\
    \geqslant \ C^{-1}_{\gamma,\epsilon}\ \Bigg[
    \sup_{0\leqslant t \leqslant t_0}\ \int_{\{t\}\cap\RR}\
        \tp^{2s-2}u^2\, |\frac{1}{u\bu}D_{\bL}(u\bu\phi)|^2 \ w_\gamma
        \ dx \\
    + \ \sup_{u}\ \int_{C(u)\cap \{0\leqslant t \leqslant t_0\}}\
    \tp^{2s-2}\bu^2\, |\frac{1}{u\bu}D_L(u\bu \phi)|^2\ w_\gamma \
        dV_{C(u)}\\
    + \ \int_0^{t_0}\int_{\RR}\ \tp^{2s-2}\tz^{1+2\epsilon}u^2\,
    |\frac{1}{u\bu}D_{\bL}(u\bu\phi)|^2 \ w'_{\gamma,\epsilon} \
        dxdt\ \Bigg] \ . \label{prilim_main_phi_est2}
\end{multline}\ret
We now form estimate \eqref{main_phi_L2_est} by adding together
estimates \eqref{prilim_main_phi_est1},
\eqref{prilim_main_phi_est1.5}, and \eqref{prilim_main_phi_est2},
while also using some further identities to simplify the resulting
right and left hand sides. First, notice that by again using the
expressions \eqref{aux_mor_id1}--\eqref{aux_mor_id2} we can
estimate:
\begin{align}
        \tm^{2s}|D_{\bL}\phi|^2 + \tp^{2s}|\frac{\phi}{r}|^2 \ &\lesssim
        \ \tp^{2s-2}u^{2}\, |\frac{1}{u\bu}D_{\bL}(u\bu\phi)|^2 +
        |D_{\bL}\phi|^2 + \tm^{2s}|\frac{1}{r}D_{\bL}(r\phi)|^2 \ , \notag \\
    \tp^{2s}|\frac{1}{r} D_L(r\phi)|^2 + \tp^{2s}
    |\frac{u\phi}{\bu r}|^2 \ &\lesssim \
    \tp^{2s-2}\bu^2\, |\frac{1}{u\bu}D_L(u\bu\phi)|^2
    + |D_L\phi|^2 + \tp^{2s}|\frac{1}{r} D_L(r\phi)|^2
    \ . \notag
\end{align}
Thus, using the bound:
\begin{equation}
    \tp^{2s} |\frac{u\phi}{\bu r}|^2 \ \lesssim \
    \tm^{2s} |\frac{\phi}{r}|^2 \ , \notag
\end{equation}
and the energy notation \eqref{phi_slab_eng}, and expanding out the
$t=0$ energy expressions on the left hand side of the resulting
estimate we conclude that:
\begin{multline}
        E_0^{(s,\gamma,\epsilon)}[\phi] \ \leqslant \
        C_{\gamma,\epsilon} \Bigg[
        \int_{0}^{t_0}\int_{\RR}\ |G|\cdot\left(\tp^{2s}
    |\frac{1}{r} D_L(r\phi)| + \tm^{2s}|
        D_{\bL}\phi| + \tm^{2s}|\frac{\phi}{r}| \right)\cdot w_\gamma \ dxdt \\
        \int_0^{t_0} \int_{\RR} \ \tp^{2s+1}\left(
        |\alpha|\cdot|\frac{\phi}{r} |\cdot|\sD\phi| +
        |\rho|\cdot|\frac{\phi}{r}|\cdot|\frac{1}{r}D_L(r\phi)|\right)\cdot
        w_\gamma \ dxdt\\
        + \ \int_0^{t_0}\int_{\RR}\
        \tp\tm^{2s}\left(|\balpha|\cdot|\frac{\phi}{r}|\cdot|\sD\phi|
        + |\rho|\cdot|\frac{\phi}{r}|\cdot |D_{\bL} \phi| \right)
    \cdot w_\gamma \ dxdt\\
    + \ \int_{\{t=0\}\cap\RR}\ (1+r^2)^{s+\gamma} \left( |D\phi|^2
    + |\frac{\phi}{r}|^2 \ \right) \ dx \
    \Bigg] \ . \notag
\end{multline}
Setting:
\begin{equation}
        ||| G |||^2_{L^2[0,t_0](L^2)(s,\gamma)} \ = \
    \int_0^{t_0}\int_{\RR} \ \tp^{2s}\tm |G|^2 \
    w_{\gamma,\epsilon}\ dxdt \ , \notag
\end{equation}
and using the notation \eqref{F_Linfty_sharp_norm}, the
condition that $\epsilon \leqslant s-\frac{1}{2}$, and again
using the bound $w_\gamma^2 \leqslant w_{\gamma,\epsilon}
\tm w'_{\gamma,\epsilon}$ along with several
different instances of the Cauchy-Schwartz inequality we have:
\begin{multline}
         E_0^{(s,\gamma,\epsilon)}[\phi] \ \leqslant \
     C_{\gamma,\epsilon} \Bigg[  \ ||| G
     |||_{L^2[0,t_0](L^2)(s,\gamma)}\cdot
     (E_0^\frac{1}{2})^{(s,\gamma,\epsilon)}[\phi] \ + \
     \llp{F}{L^\infty[0,t_0](\epsilon)}\cdot
     E_0^{(s,\gamma,\epsilon)}[\phi] \notag \\
     + \ \int_{\{t=0\}\cap\RR}\ (1+r^2)^{s+\gamma} \left( |D\phi|^2
     + |\frac{\phi}{r}|^2 \ \right) \ dx \
     \Bigg] \ . \notag
\end{multline}
Dividing this last expression through by
$(E^\frac{1}{2})^{(s,\gamma,\epsilon)}[\phi]$ and squaring, we easily
achieve the bound:
\begin{multline}
         E^{(s,\gamma,\epsilon)}[\phi] \ \leqslant \
     C^2_{\gamma,\epsilon} \Bigg[  \ ||| G
     |||^2_{L^2[0,t_0](L^2)(s,\gamma)} \ + \
     \llp{F}{L^\infty[0,t_0](\epsilon)}^2\cdot
     E^{(s,\gamma,\epsilon)}[\phi] \notag \\
     + \ \int_{\{t=0\}\cap\RR}\ (1+r^2)^{s+\gamma} \left( |D\phi|^2
    + |\frac{\phi}{r}|^2 \ \right) \ dx \
     \Bigg] \ . \notag
\end{multline}
The estimate \eqref{main_phi_L2_est} now follows from this last line
and the following gauge covariant Poincare type estimate which follows
from the same reasoning used to produce \eqref{first_poincare}
which we shall proved in the sequel:
\begin{equation}
        \int_{\{t=0\}\cap\RR}\ (1+r^2)^{s+\gamma}|\frac{\phi}{r}|^2 \
        dx \ \lesssim \ \int_{\{t=0\}\cap\RR}\
    (1+r^2)^{s+\gamma} \ |D\phi|^2 \ dx \ . \notag
\end{equation}
\end{proof}

\ret
%-------------------------------------------------------------------------
%%%%%%%%%%%%%%%%%%%%%%%%%%%%%%%%%%%%%%%%%%%%%%%%%%%%%%%%%%%%%%%%%%%%%%%%%%
%-------------------------------------------------------------------------

\section{$L^\infty$ Estimates for Electro-Magnetic Fields}\label{F_Linfty_sect}

This section is the sequel to Section \ref{F_L2_section} above in
that it contains the second main set of estimates for the field
strength $F_{\alpha\beta}$. Here we will prove $L^\infty$ type
estimates for this quantity without assuming that it necessarily
satisfies the field equations \eqref{basic_MKG}. That is, in this
section we will only assume that $F_{\alpha\beta}$ is an arbitrary
two-form and will prove $L^\infty$ type estimates for it at the cost
of bounds on the energies \eqref{F_s_energy_def}. This is done through
the well known procedure of applying global Sobolev inequalities in
conjunction with geometric identity \eqref{otf_lie_der2} for the
component decomposition of the Lie derivative of a two-form. This
procedure was first explored in \cite{CK_Fields}, and
our treatment adds little except a bound on the $L^2(L^\infty)$ term
in the norm \eqref{F_Linfty_norm}. This turns out to be a direct
consequence of the space-time $L^2$ estimate \eqref{main_F_L2_est2}
and the weighted Sobolev estimates
\eqref{ext_sob1}--\eqref{ext_sob2} and \eqref{int_sob}. However, we
choose to prove the entire $L^\infty$ estimate from scratch here
primarily for the sake of completeness and because we use the
fractional weights
$\tp^{2s},\tm^{2s}$.\\

Before we proceed with the main estimate of this section, we
introduce two notational devices which will allow us to reduce things
to a more simple form where the estimates
\eqref{ext_sob1}--\eqref{ext_sob2}
and \eqref{int_sob} can be used directly. The first of these is a
simple algebraic tool which will make dealing with various commutators
more straight forward. This is the introduction of the ``radial boost''
field:
\begin{equation}
        \Omega_{0r} \ = \ \omega^i \Omega_{0i} \
    = \ -\frac{1}{2} \left(\bu L - u \bL\right) \
    . \label{radial_boost}
\end{equation}
It has the advantage of commuting with the rotation and
scaling fields $\{\Omega_{ij},S\}$.
Furthermore, it turns out that $\Omega_{0r}$ preserves
the null decomposition \eqref{null_decomp} even after one
reexpands $\Omega_{0r}$ and passes back to the usual boosts
$\Omega_{i0}$. This remarkable
fact\footnote{ It is not a-priori clear that this should happen given
the fact that the boosts $\Omega_{i0}$ cause quite a bit of
permutation in the null decomposition \eqref{null_decomp}.}
as well as the rule that we only need one instance of the boosts
in the exterior according to \eqref{ext_sob1} will help to streamline
the steps we take in the sequel.\\

The second tool we use here is a dyadic decomposition of the
distance to the standard light-cone $u=0$, as well as the distance
along its translates $u=const.$ We will use this tool to localize
all of our $L^\infty$ estimates to shells adapted to these
distances. This allows us to add extra powers of $\tp,\tm$
at will because they just appear as constants in our
localized estimates. The language for this is as follows:\\

We first isolate on each fixed time slice $t=const.$ the extended
exterior region $t < 2r$. We then chop this into dyadic pieces depending
on the distance from the cone $u=0$.
We decompose the extended exterior region into the set of double spherical shells:
\begin{align}
        \mathcal{J}_i \ &= \ \{x\in\{t\}\times\mathbb{R}^3 \
    , t<2r \ \big| \
    2^{i} \leqslant  \big|t-|x|\big|  + 1 \leqslant 2^{i+1}\}
    \ , &i\in\mathbb{N}  \ . \notag
\end{align}
For fixed $0 < i$, each of these sets can be further decomposed
into exterior and inter portions $\mathcal{J}_i = \mathcal{J}^+_i
\cup \mathcal{J}^-_i$, where $ \mathcal{J}^\pm_i$ contains only the
points where $0 < \mp u$. Note that these sets are connected, and
that  $\mathcal{J}^-_i$ empty for $i$ sufficiently large. In general,
we will write $\mathcal{J}$ for one of these connected sets. We also
form a partition of unity adapted to the decomposition
$\mathcal{J}_i$ which has the homogeneity property:
\begin{equation}
        |\partial_r \chi_{\mathcal{J}_i}| \ \lesssim \
         \tm^{-1}(\mathcal{J}_i) \ , \notag
\end{equation}
where:
\begin{align}
        \tp(\mathcal{J}_i) \ &= \ \min_{x\in\mathcal{J}_i} \ \tp \ ,
    &\tm(\mathcal{J}_i) \ &= \ \min_{x\in\mathcal{J}_i} \ \tm \ .
    \label{tptm_J}
\end{align}
Furthermore, we note that it is clear
that we can also choose the largest $J^-_i$ so its cutoff
is supported where $r \leqslant \frac{1}{4} t $. When the context is understood,
we will simply substitute $\chi$ for $\chi_{\mathcal{J}_i}$
and $\tp,\tm$ for $\tp(\mathcal{J}),\tm(\mathcal{J})$.\\

Next, we deal with the light-cones $u=const.$ The notation we use
here is similar to what was outlined in the preceding paragraph, and
we will only consider the portion of $C(u)$ in the region where
$t<2r$. For each such cone (omitting the dependence on u), we define
the dyadic conical spherical shells:
\begin{align}
        \mathcal{I}_i \ &= \ \{x\in C(u) \ , \ 0\leqslant t
    < 2r \ \ \big| \
    2^{i} \leqslant  t + |x| + 1 \leqslant 2^{i+1}\}
    \ , &i\in\mathbb{N}  \ . \notag
\end{align}
We use the notation $\tp(\mathcal{I}_i)$,$\tm(\mathcal{I}_i)$ in
analogy with line \eqref{tptm_J} above. We also introduce a smooth
partition of unity $\chi_{\mathcal{I}_i}$ adapted to the $\mathcal{I}_i$
which satisfies the bounds:
\begin{equation}
        |L(\chi_{\mathcal{I}_i})| \ \lesssim \
    \tp^{-1}(\mathcal{I}_i) \ . \notag
\end{equation}
In particular, notice that we have not cut off strictly along the
truncated cone $C(u)\cap \{0 \leqslant t \leqslant t_0\}$, but we may
assume that each of these $\chi$ are supported where $t < 4r$.\\

Using these notations we now glue together the estimates
\eqref{ext_sob1}--\eqref{ext_sob2}, \eqref{int_sob}, and
\eqref{char_sob1}--\eqref{char_sob2}  of the Appendix
into a form which will be used in the sequel:\\

\begin{lem}[Weighted global Sobolev estimates]
Let $f$ be a test function on $\mathbb{R}\times\RR$, and let
$w_\gamma$ be the weight function defined on line
\eqref{w_gamma_def1}.
Then on each time slice $t=const.$ with $1 \leqslant t$
one has the following (non-uniform) estimates for arbitrary weights
$\delta_\pm \in \mathbb{R}$:
\begin{align}
        \sup_{t \leqslant 2r}  \ |\tp^{1+ \delta_+ }
    \tm^{\frac{1}{2}+\delta_-} f|^2 \, w_\gamma \ &\lesssim \
    \sum_{\substack{|I| \leqslant 1 , |J|\leqslant 2\\
    |I|+|J| \leqslant 2\\
    X\in \{\partial_r,S,\Omega_{0r}\} , Y\in \mathbb{O} }}
    \lp{\tp^{\delta_+}\tm^{\delta_-} w^\frac{1}{2}_\gamma\,
    X^IY^J (f)}{L_x^2(t < 4r)}^2 \ ,
    \label{ext_slice_sob}\\
    \sup_{2r<t}  \ |    \tp^{\frac{3}{2} + \delta_+ }
    \,     f|^2 \ &\lesssim \
    \sum_{\substack{|I| \leqslant 2 \\
    X\in \{S,\Omega_{0i}\} }}
    \lp{\tp^{\delta_+}\,
    X^I (f)}{L_x^2(r< \frac{3}{4}t )}^2 \ ,
    \label{int_slice_sob}
\end{align}
Furthermore, defining the truncated cones:
\begin{align}
        \overline{C}(u) \ &= C(u)\cap\{1\leqslant t\leqslant t_0\}
    \cap\{ t < 2r \} \ ,\notag\\
    \td{C}(u) \ &= C(u)\cap\{1\leqslant t\leqslant t_0\}
    \cap\{ t < 4r \} \ ,\notag
\end{align}
we have the characteristic $L^\infty$ estimate:
\begin{equation}
    \sup_{(t,x)\in \overline{C}(u)}
    \ | \tp^{\frac{3}{2} + \delta_+ }\,     f|^2 \ \lesssim \
    \sum_{\substack{|I| \leqslant 1 , |J|\leqslant 2 \\
    |I|+|J| \leqslant 2\\
    X\in \{ \Omega_{0r},S\} ,  Y\in\mathbb{O}       }}
    \lp{\tp^{\delta_+}\,
    X^I Y^J (f)}{L_x^2\big(\td{C}(u)\big)}^2 \ .
    \label{char_slice_sob}
\end{equation}
In the above estimates, $\mathbb{O}=\{\Omega_{ij}\}$
denotes the Lie algebra of the rotation group.
\end{lem}\ret

\begin{proof}[Proof of \eqref{ext_slice_sob}--\eqref{char_slice_sob}]
These follow almost directly from their local versions
\eqref{ext_sob1}--\eqref{ext_sob2}, \eqref{int_sob}, and
\eqref{char_sob1}--\eqref{char_sob2}. We'll just deal with
\eqref{ext_slice_sob} here as the other two estimates follow from
similar reasoning and are left to the reader. For a dyadic shell
$J$, by combining estimates \eqref{ext_sob1}--\eqref{ext_sob2} in
tandem (for example with the $L^q$ set to $q=3$) and using the
estimate:
\begin{equation}
        |\tm \partial_r \chi f| \ \lesssim \ |\td{\chi} \, f| +
    |\td{\chi}\, S(f)| + |
        \td{\chi}\, \Omega_{0r}(f)| + |\td{\chi}\, \partial_r f| \ , \notag
\end{equation}
where $\td{\chi}$ is a cutoff on $\hbox{supp}\{\chi\}$,
we arrive at:
\begin{equation}
        \sup_{x}  \ \tp^{2}
    \tm|\chi f|^2 \ \lesssim \
    \sum_{\substack{|I| \leqslant 1 , |J|\leqslant 2\\
    |I|+|J| \leqslant 2\\
    X\in \{\partial_r,S,\Omega_{0r}\} , Y\in \mathbb{O} }}
    \lp{\td{\chi} X^I Y^J (f)}{L^2_x}^2 \ . \notag
\end{equation}
Since both sides of this last estimate are restricted to a doubling
of the shell $\mathcal{J}$, we can multiply both sides by the
constant $\tp^{2\delta_+}(\mathcal{J})\tm^{2\delta_-
+2\gamma}(\mathcal{J})$ if $t<r$, and by the constant
$\tp^{2\delta_+}(\mathcal{J})\tm^{2\delta_-}(\mathcal{J})$ if $r<t$.
Doing this and then summing the result over all shells $\mathcal{J}$
in the extended exterior $t\leqslant 2r$ yields the desired result.
\end{proof}\ret

\noindent We now state and prove the main result of this section:\\

\begin{prop}[General $L^\infty$ estimate for electro-magnetic fields]
\label{main_ext_est}
Let $F_{\alpha\beta}$ be an arbitrary two-form on
$\mathbb{R}\times\RR$, and let $w_\gamma$ be the weight function
defined on line \eqref{w_gamma_def1}.
Then in terms of null decomposition \eqref{null_decomp}
and the general energy norm
\eqref{F_s_energy_def}, one has the following $L^\infty$ estimate:
\begin{multline}
     \sup_{\substack{0\leqslant t \leqslant t_0\\ x\in\{t\}\times\RR}}
     \left(\tp^{2s+3}|\alpha|^2 + \tp^2\tm^{2s+1} |\balpha|^2
     + \tp^{2s+2}\tm( \rho^2 + \sigma^2)\right)\cdot w_\gamma\\
     + \  \lp{\tp^{s+1}\tm^\frac{1}{2}
     (w')^\frac{1}{2}_{\gamma,\epsilon}\, \alpha }
     {L_t^2(L^\infty_x)[0,t_0]}^2
     \  \lesssim \ \sum_{\substack{|I|\leqslant 2\\
     X\in\mathbb{L}}}E^{(s,\gamma,\epsilon)}(0,t_0)[\mathcal{L}_X^I F] \
     . \label{main_F_linfty_est}
\end{multline}
In particular, using the charge modified energy norm \eqref{k_F_eng}, and the
general weighted $L^\infty$ norm \eqref{F_Linfty_norm} we have that:
\begin{equation}
        \llp{F}{L^\infty[0,t_0](s,\gamma,\epsilon)} \ \lesssim \
    E_2^{(s,\gamma,\epsilon)}(0,t_0)[F] \ . \label{charge_decomp_F_linfty}
\end{equation}
\end{prop}\ret

\begin{proof}[Proof of estimate \eqref{main_F_linfty_est}]
In what follows, we shall only consider the estimate
\eqref{main_F_linfty_est} in the region $1\leqslant t$. For
$0\leqslant t \leqslant 1$ the result follows from a weighted estimate
on the electric-magnetic decomposition \eqref{em_decomp} in
conjunction with a simple weighted version of the usual Sobolev
embedding. We leave the details of this to the reader.\\

The proof can now be broken down into three sections. All of the
components in the first term on the left hand side of
\eqref{main_F_linfty_est} except $\alpha$ can be treated with fixed
time energy estimates. The pure $L^\infty$ estimate for $\alpha$ is
a consequence of characteristic energies, while the mixed norm
estimate follows from
the space-time energy estimate contained in \eqref{main_F_L2_est2}.\\

\subsection*{Bounds on $\{\balpha,\rho,\sigma\}$.}
We first prove this assertion for the  $\balpha$ component in both
the extended exterior region $t<2r$ and its compliment $r <
\frac{1}{2} t$. In fact, it is just as easy to demonstrate the full
bound:
\begin{multline}
        \sup_{x} \ \tp^2 \tm^{2s+1}\left(
        |E|^2 + |H|^2
        \right)\cdot w_{\gamma} \ \lesssim \\
        \sum_{\substack{|I| \leqslant 1 , |J|\leqslant 2\\
        |I|+|J| \leqslant 2\\
        X\in \{\partial_r,S,\Omega_{0r}\} , Y\in \mathbb{O} }}
        \lp{\tm^{s} w^\frac{1}{2}_\gamma\,
        \big(E(\mathcal{L}_X^I \mathcal{L}_Y^J F ) , H(\mathcal{L}_X^I \mathcal{L}_Y^J
        F)\big)
        }{L_x^2}^2 \ . \label{full_bound_for_alpha}
\end{multline}
To prove \eqref{full_bound_for_alpha}, we first combine estimates
\eqref{ext_slice_sob} and \eqref{int_slice_sob} with $\delta^+
=0$,$\delta^-=s$ and $\delta^+=s$ respectively. This yields the
following bound for scalar components:
\begin{multline}
        \sup_{x} \ \tp^2 \tm^{2s+1}\left(
        |E|^2 + |H|^2
        \right)\cdot w_{\gamma} \ \lesssim \\
        \sum_{\substack{|I| \leqslant 1 , |J|\leqslant 2\\
        |I|+|J| \leqslant 2\\
        X\in \{\partial_r,S,\Omega_{0r}\} , Y\in \mathbb{O} }}
        \lp{\tm^{s} w^\frac{1}{2}_\gamma\,
        \big(X^I Y^J (E)  , X^I Y^J (H) \big) }{L_x^2}^2 \ . \notag
\end{multline}
To finish things up, we need only show the following point-wise
bounds:
\begin{multline}
        \sum_{\substack{|I| \leqslant 1 , |J|\leqslant 2\\
        |I|+|J| \leqslant 2\\
        X\in \{\partial_r,S,\Omega_{0r}\} , Y\in \mathbb{O} }}
        \big|  \big(X^I Y^J (E)  , X^I Y^J (H) \big) \big|^2
        \ \lesssim \\
        \sum_{\substack{|I| \leqslant 1 , |J|\leqslant 2\\
        |I|+|J| \leqslant 2\\
        X\in \{\partial_i,S,\Omega_{0i}\} , Y\in \mathbb{O} }}
        \ \big|\big(E(\mathcal{L}_X^I \mathcal{L}_Y^J F ) ,
        H(\mathcal{L}_X^I \mathcal{L}_Y^J F)\big) \big|^2 \ .
        \notag
\end{multline}
This, in turn, follows immediately from expanding out the
vector-fields $\partial_r$ and $\Omega_{0r}$ in terms of elements of
$\mathbb{L}$, and then using the Lie derivative formula
\eqref{otf_lie_der2} in conjunction with the commutator identity
\eqref{lor_bracket_rel1}. The important thing to note here is that
in this latter formula one always ends up with constant linear
combinations of elements of the translation invariant frame
$\{\partial_\alpha\}$. This allows one to write the bracket terms on
the right hand side of \eqref{otf_lie_der2} directly in terms of
constant coefficient linear combinations of $(E,H)$.\\

Our next step is to bound the terms $\rho$ and $\sigma$ in the
extended exterior $t<2r$. This, combined with estimate
\eqref{full_bound_for_alpha} above in the deep interior
$r<\frac{1}{2}t$, will demonstrate \eqref{main_F_linfty_est} for
these components. We now need to be a bit more careful regarding how
the operation of Lie differentiation permutes null components. To
take this into account, we will again make use of
\eqref{otf_lie_der2} as well as the expansion identities
\eqref{lor_null_brack} for  Lie derivatives, with some additional
help from the frame identities on line \eqref{omega_relations}.
First, notice that to expand vector-fields of the form $X_r =
\omega^i X_i$, the formula \eqref{otf_lie_der2} implies:
\begin{equation}
        (\mathcal{L}_{X_r}F_{\alpha\beta}) \
        = \ \omega^i (\mathcal{L}_{X_i} F)_{\alpha\beta} +
        e_\alpha(\omega^i)X_i^\gamma F_{\gamma \beta}
        + e_\beta(\omega^i) X_i^\gamma F_{\alpha\gamma}
        \ . \label{omega_otf_lie_der2}
\end{equation}
Therefore, after some simple computations using the first of the
special identities (the second will be used in a moment):
\begin{align}
        e_A(\omega^i)\omega_i \ &= \ 0 \ ,
        &e_A(\omega^i)\omega_i^B \ &= \ \frac{1}{r}\delta_A^{\ B} \ ,
        \label{special_omega_iden}
\end{align}
we arrive at the following equalities:
\begin{subequations}\label{rho_lie_component}
\begin{align}
        \partial_r (\rho) \ &= \
        \omega^i \rho(\mathcal{L}_{\partial_i} F) \ ,
        &\Omega_{ij} (\rho) \ &= \ \rho(\mathcal{L}_{\Omega_{ij}}
        F) \ , \\
        S(\rho) \ &= \ \rho(\mathcal{L}_{S}
        F) - 2\rho \ ,
        &\Omega_{0r}(\rho ) \ &= \ \omega^i
        \rho(\mathcal{L}_{\Omega_{0i}} F) \ .
\end{align}
\end{subequations}
and:
\begin{subequations}\label{sigma_lie_component}
\begin{align}
        \partial_r(\sigma) \ &= \ \omega^i
        \sigma(\mathcal{L}_{\partial_r} F) \ ,
        &\Omega_{ij}(\sigma) \ &= \ \sigma(\mathcal{L}_{\Omega_{ij}}F) \ ,
        \\
        S(\sigma) \ &= \ \rho(\mathcal{L}_{S}
        F) - 2\sigma \ ,
        &\Omega_{0r}(\sigma ) \ &= \ \omega^i
        \sigma(\mathcal{L}_{\Omega_{0i}} F) \ .
\end{align}
\end{subequations}
Notice that the set \eqref{rho_lie_component} is easily proved using
the identities \eqref{otf_lie_der2}, \eqref{omega_otf_lie_der2},
\eqref{omega_relations}, and \eqref{lor_null_brack}, while the
second set \eqref{sigma_lie_component} follows from
\eqref{rho_lie_component} and the duality formulas
\eqref{star_lie_comm} and \eqref{dual_null_decomp}.\\

Next, we see that by differentiating both sides of the expressions
\eqref{rho_lie_component}--\eqref{sigma_lie_component} with respect
to the rotations $\Omega_{ij}$ and making use of the  homogeneity
bound:
\begin{equation}
        |\Omega_{ij}(\omega_k)| \ \lesssim \
    1 \ , \label{first_omega_homog_bound}
\end{equation}
we arrive at the point-wise bound:
\begin{equation}
        \sum_{\substack{|I| \leqslant 1 , |J|\leqslant 2\\
        |I|+|J|\leqslant 2\\
        X\in \{\partial_r,S,\Omega_{0r}\} , Y\in \mathbb{O} }}
        |X^IY^J(\rho,\sigma)|^2 \ \lesssim \
        \sum_{\substack{|I| \leqslant 1 , |J|\leqslant 2\\
        |I|+|J|\leqslant 2\\
        X\in \{\partial_i,S,\Omega_{0i}\} , Y\in \mathbb{O} }}
        \ \rho^2(\mathcal{L}_X^I\mathcal{L}_Y^J F) +
        \sigma^2(\mathcal{L}_X^I\mathcal{L}_Y^J F ) \ .
        \label{rs_lie_exp_bound}
\end{equation}
This last estimate, combined with the scalar $L^\infty$ estimate:
\begin{equation}
        \sup_{t<2r} \ \tp^{2s +2} \tm w_{\gamma} (\rho^2 + \sigma^2) \ \lesssim \
        \sum_{\substack{|I| \leqslant 1 , |J|\leqslant 2\\
        |I|+|J|\leqslant 2\\
        X\in \{\partial_r,S,\Omega_{0r}\} , Y\in \mathbb{O} }}
        \lp{\tp^{s} w^\frac{1}{2}_\gamma\, X^IY^J(\rho,\sigma)}{L_x^2(t\leqslant
        4r)}^2 \ , \label{ext_rho_sigma_linfty_bound}
\end{equation}
which follows from  the bound \eqref{ext_slice_sob}
with $\delta_+ = s$ and $\delta_- =0$, implies the
estimate \eqref{main_F_linfty_est} for the components $\rho$ and
$\sigma$ in the exterior region $t < 2r$. We have now completed the
proof of \eqref{main_F_linfty_est} for the components $\balpha$,
$\rho$, and $\sigma$.\\

\subsection*{The $L^2(L^\infty)$ estimate for $\alpha$}
We now prove the mixed norm estimate:
\begin{equation}
        \lp{\tp^{s+1} \tm^\frac{1}{2}
    (w')^\frac{1}{2}_{\gamma,\epsilon}
    \alpha}{L^2(L^\infty)}^2 \ \lesssim \
    \sum_{\substack{|I|\leqslant 2\\
     X\in\mathbb{L}}}E^{(s,\gamma,\epsilon)}(0,t_0)[\mathcal{L}_X^I F]
    \ . \label{F_L2Linfty_est}
\end{equation}
Due to the space-time energy estimate for $\alpha$ contained in
the right hand side of \eqref{F_L2Linfty_est}, it suffices to prove
the two fixed time estimates:
\begin{align}
        \sup_{t<2r} \tp^{2s+2} \tm |\alpha|^2
        \, w'_{\gamma,\epsilon} \ &\lesssim \
        \sum_{\substack{|I|\leqslant 1 , |J|\leqslant 2\\
        |I|+|J|\leqslant 2\\
        X\in\{\partial_i , S, \Omega_{0i}\} , Y \in \mathbb{O}}}
        \ \lp{\tp^s
        (w')^\frac{1}{2}_{\gamma,\epsilon}
        \alpha(\mathcal{L}^I_X\mathcal{L}^J_Y F)}{L^2(t < 4r)}^2
         \ , \label{alpha_fixed_time1}\\
        \sup_{r < \frac{1}{2}t } \tp^{2s+2 - 2\epsilon} |\alpha|^2
        \ &\lesssim \ \sum_{\substack{|I|\leqslant 2 \\
        X\in\{S, \Omega_{0i}\} }}\
        \lp{\tp^{s-1-\epsilon} \Big(E(\mathcal{L}^I_X F)
        ,H(\mathcal{L}^I_X F)\Big)}
        {L^2(r < \frac{3}{4}t)}^2
        \ . \label{alpha_fixed_time2}
\end{align}
These in turn follow from the scalar estimates\footnote{It is clear
that one can replace $w_\gamma$ in estimate \eqref{ext_slice_sob}
with $w'_{\gamma,\epsilon}$.}
\eqref{ext_slice_sob}--\eqref{int_slice_sob} and the same procedure
used above to prove the bounds \eqref{full_bound_for_alpha} and
\eqref{ext_rho_sigma_linfty_bound} above. Specifically, notice that
by the same types of calculations used to produce
\eqref{rho_lie_component}--\eqref{rho_lie_component}, in particular
\eqref{otf_lie_der2}, \eqref{lor_null_brack},
\eqref{omega_otf_lie_der2}, and the identities
\eqref{special_omega_iden}, we have that:
\begin{subequations}\label{alpha_lie_component}
\begin{align}
        \partial_r (\alpha_A) \ &= \
        \omega^i \alpha_A(\mathcal{L}_{\partial_i} F) \ , \\
        \Omega_{ij} (\alpha_A) \ &= \ \alpha_A(\mathcal{L}_{\Omega_{ij}}
        F) + [\Omega_{ij},e_A]^B \alpha_{B} \ , \\
        S(\alpha_A) \ &= \ \alpha_A(\mathcal{L}_{S}
        F) - 2 \alpha_A \ , \\
        \Omega_{0r}(\alpha_A) \ &= \ \omega^i
        \alpha_A(\mathcal{L}_{\Omega_{0i}} F) + \alpha_A \ .
\end{align}
\end{subequations}
These allow us to pass from a scalar estimate of the form
\eqref{ext_slice_sob} (with $\delta^+ = s$ and $\delta^- =0$) in
\eqref{alpha_fixed_time1} to the right hand side of this estimate
which contains Lie derivatives. To do this we simply differentiate
both sides of the expressions \eqref{alpha_lie_component} with
respect to rotations, and then express the resulting right hand side
in terms of the collection \ $\{\alpha, \alpha(\mathcal{L}_X F),
\alpha(\mathcal{L}_X \mathcal{L}_Y F)\}$ \ for $X\in\mathbb{L}$ and
$Y\in \mathbb{O}$. In doing this, if the extra angular derivative
lands on a $\alpha$ component on the right hand side of
\eqref{alpha_lie_component}, we can use these formulas to re-expand
the result using the above remarks to bound the resulting
coefficients. If, on the other hand, the extra angular derivative
lands on a coefficient on the right hand side of
\eqref{alpha_lie_component} we can bound things by using
\eqref{first_omega_homog_bound} and its companion estimate:
\begin{equation}
        |\Omega_{jk}([\Omega_{lm},e_A]^B)|
        \ \lesssim \ \ 1 \
        , \label{second_omega_homog_bound}
\end{equation}
which easily follows from the homogeneity property \eqref{ang_frame}
of the fields $\{e_A\}$, and the fact that $\Omega^A_{ij} = r
\td{\Omega}^A_{ij}$ where $\td{\Omega}^A_{ij}$ is a function of the
angular variable only.\\

\subsection*{The $L^\infty$ bound for $\alpha$}
We need only consider the case $t < 2r$, as the desired estimate in
the compliment was proved on line \eqref{full_bound_for_alpha}
above. That is, we only need to prove this bound along the cones
$\overline{C}(u)$ which were introduced just before estimate
\eqref{char_slice_sob} above. Using that estimate with $\delta_+ =
s$ , we have the scalar bound:
\begin{equation}
        \sup_{(t,x)\in\overline{C}(u)}
    \ \tp^{2s +3} w_{\gamma} |\alpha|^2 \ \lesssim \
    \sum_{\substack{|I| \leqslant 1 , |J|\leqslant 2\\
        |I|+|J|\leqslant 2\\
    X\in \{S,\Omega_{0r}\} , Y\in \mathbb{O} }}
    \lp{\tp^{s} w^\frac{1}{2}_\gamma\,
    X^IY^J(\alpha)}{L_x^2\big(\td{C}(u)\big)}^2 \ . \notag
\end{equation}
Notice that the $w_\gamma$ factor can be added in at will because it
is a constant on the cones $C(u)$. Using this last line, we see that
is suffices to prove the point-wise bound:
\begin{equation}
        \sum_{\substack{|I| \leqslant 1 , |J|\leqslant 2\\
        |I|+|J|\leqslant 2\\
    X\in \{S,\Omega_{0r}\} , Y\in \mathbb{O} }}
    |X^IY^J(\alpha)|^2 \ \lesssim \
    \sum_{\substack{|I| \leqslant 1 , |J|\leqslant 2\\
        |I|+|J|\leqslant 2\\
    X\in \{S,\Omega_{0i}\} , Y\in \mathbb{O} }}
    \ |\alpha(\mathcal{L}_X^I\mathcal{L}_Y^J F)|^2 \ .
    \notag
\end{equation}
This in turn follows from the decomposition formulas
\eqref{alpha_lie_component} and the coefficient bounds
\eqref{first_omega_homog_bound} and
\eqref{second_omega_homog_bound}. The details are left to the reader.
This completes our proof of estimate \eqref{main_F_linfty_est}.
\end{proof}\ret

We end this section with a discussion of the peeling properties of
Lie derivatives of the charge portion $\overline{F}$ of the
curvature $F$. This will be used in the following sections wherever
the quantities $\mathcal{L}^I_X \overline{F}_{\alpha\beta}$ occur
because, while there are  nice explicit formulas for these objects,
they obey no useful space-time estimates.\\

\begin{prop}[Peeling properties of the pure charge fields
$\overline{F}_{\alpha\beta}$] Let $F_{\alpha\beta}$ be an
electro-magnetic field, and let $\overline{F}_{\alpha\beta}$ be its
pure charge component as defined on the line
\eqref{null_charge_tensor}. Then for each multiindex $I$ we have the
following $L^\infty$ estimates:
\begin{align}
        |\alpha(\mathcal{L}^I_X \overline{F})| \ &\leqslant \
        C_I\
        |q(F)|\cdot \tz \tp^{-2} \chi_{t < r+1 } \ ,
        \label{alpha_bF_peel}\\
        |\balpha(\mathcal{L}^I_X \overline{F})|
        \ , \ |\rho(\mathcal{L}^I_X \overline{F}))| \ ,
        \ |\sigma(\mathcal{L}^I_X \overline{F}))| \ , \
        &\leqslant \ C_I\ |q(F)|\cdot \tp^{-2}
        \chi_{t < r+1 } \ ,
        \label{other_bF_peel}
\end{align}
where all $X\in\mathbb{L}$.
\end{prop}\ret

\begin{proof}[Proof of the estimates
\eqref{alpha_bF_peel}--\eqref{other_bF_peel}] The proof is a simple
inductive procedure similar to what was done for
$\overline{J}_\alpha$ in section \ref{F_L2_section} starting on line
\eqref{I0_line}. We begin by making the inductive hypothesis that:
\begin{align}
        \alpha(\mathcal{L}^I_X \overline{F}) \ &= \
        q\cdot \sum_{k=1}^{|I|}\
        \ \frac{
        \Omega^I_\alpha(\omega) }{r^{2+k}} \cdot (\chi_k^+)_\alpha^I(r-t-2) \ , \label{bF_indct1}\\
        \balpha(\mathcal{L}^I_X \overline{F}) \ &= \ q\cdot \sum_{k=0}^{|I|}\
        \ \frac{\Omega^I_{\balpha}(\omega)
        }{r^{2+k}} \cdot (\chi_k^+)_{\balpha}^I(r-t-2) \ , \label{bF_indct2}\\
        \rho(\mathcal{L}^I_X \overline{F}) \ &= \ q\cdot \sum_{k=0}^{|I|}\
        \ \frac{\Omega^I_\rho(\omega)
        }{r^{2+k}} \cdot (\chi_k^+)_\rho^I(r-t-2) \ , \label{bF_indct3}\\
        \sigma(\mathcal{L}^I_X \overline{F}) \ &= \ q\cdot \sum_{k=0}^{|I|}\
        \ \frac{\Omega^I_\sigma(\omega)
        }{r^{2+k}} \cdot (\chi_k^+)_\sigma^I(r-t-2) \ , \label{bF_indct4}
\end{align}
where the $\Omega^I_{\bullet}$ are smooth functions of the angular
variable whose $C^\infty$ bounds depend only on $I$ and the specific
component being considered, and the $(\chi_k^+)^I_\bullet(s)$ are
smooth functions of the single variable $s$ which vanish for
$s\leqslant -1$ and satisfy the homogeneity bound:
\begin{equation}
        |\partial^j_s (\chi_k^+)^I_\bullet| \ \leqslant \ C_I \
        s^{k-j} \ , \notag
\end{equation}
for $1\leqslant s$. Here, as in previous discussions, the product
notation $\Omega\cdot \chi^+$ is symbolic for a sum of products of
functions with these properties.\\

Now, from the formulas \eqref{null_charge_tensor}, it is clear that
the assumption \eqref{bF_indct1}--\eqref{bF_indct4} holds for
$|I|=0$. Furthermore, simple explicit calculations using the
formulas \eqref{lor_null_decomp}, and which we leave to the reader,
show that differentiating as scalars the quantities on the right
hand side of the formulas \eqref{bF_indct1}--\eqref{bF_indct4}
yields quantities with these same properties. This follows simply
from expressing everything in the derivative in terms of the
variables $u,r,\omega$. Therefore, it remains to show that one can
inductively reproduce identities of this form  after the bracket
terms on the right hand side of formula \eqref{otf_lie_der2} are
taken into account.\\

First notice that if the Lie derivative is with respect to the
fields $\{\Omega_{ij} , S\}$ then the null decomposition
\eqref{null_decomp} is preserved, so the claim follows from the fact
that the coefficients in the expansions \eqref{lor_null_brack} can
all be put in the form $\Omega(\omega)\cdot \varphi( \frac{u}{r})$
for some affine function $\varphi$.\\

Next, notice that while the translation invariant fields
$\{\partial_\alpha\}$ do not preserve the decomposition
\eqref{null_decomp}, their coefficients in the formulas
\eqref{lor_null_brack} introduce a factor which is bounded by
$\frac{1}{r}$. Therefore, it remains to deal with the bracket
portion when dealing boosts $\{\Omega_{0i}\}$. Notice that it
suffices to limit discussion to the component $\alpha$ because the
remaining components have the same general form in
\eqref{bF_indct2}--\eqref{bF_indct4} above and the homogeneity of
the coefficients in the formulas \eqref{lor_null_brack} guarantee
that this general form is preserved after multiplication by them.
The claim now follows from the explicit formulas
\eqref{lor_null_brack6} and \eqref{lor_null_brack8}, because the only
thing that can disturb the structure of \eqref{bF_indct1} is when one of
the other components $\balpha(\mathcal{L}^I_X
\overline{F}),\rho(\mathcal{L}^I_X\overline{F})$ or
$\sigma(\mathcal{L}^I_X \overline{F})$ gets moved to the $\alpha$
position. Notice that this in turn can only happen when an extra factor of
$\frac{u}{r}$ is introduced, having the effect of causing the sum in
formulas \eqref{bF_indct2}--\eqref{bF_indct4} to start with $k=1$.
This completes the demonstration of
\eqref{alpha_bF_peel}--\eqref{other_bF_peel}.
\end{proof}

\ret
%-------------------------------------------------------------------------
%%%%%%%%%%%%%%%%%%%%%%%%%%%%%%%%%%%%%%%%%%%%%%%%%%%%%%%%%%%%%%%%%%%%%%%%%%
%-------------------------------------------------------------------------

\section{$L^\infty$ estimate for complex scalar fields}\label{phi_Linfty_sect}

This section is the companion to Section \ref{phi_L2_sect}. It
contains $L^\infty$ type estimates comparable to
\eqref{main_F_linfty_est} for complex scalar fields in terms of bounds
on energy norms of the kind \eqref{phi_slab_eng}. These will be proved
as in the last section through the application of global Sobolev
estimates and several geometric identities which reduce bounds on
scalar derivatives to bounds on covariant derivatives. The device
which ultimately enables us to do this the following simple, although far
reaching, estimate of Kato:\\

\begin{lem}[Kato's ``diamagnetic'' inequality]\label{Kato_lem}
Let $\phi$ be a section of $V$ over $\mathcal{M}$, and denote by
$|\phi|^2 = \langle \phi , \phi \rangle_V$. Then one has the
point-wise inequality:
\begin{equation}
    \big|X(|\phi|)\big| \ \leqslant \ |D_X \phi | \ ,
    \label{Kato_ineq}
\end{equation}
for any vector-field $X$ on $\mathcal{M}$.
\end{lem}\ret

\begin{proof}[Proof of \eqref{Kato_ineq}]
This is a straight forward consequence of the compatibility
condition \eqref{bundle_metric}. We compute:
\begin{align}
    X(|\phi|) \ &= \ X \left(   \langle \phi ,
    \phi \rangle_V^\frac{1}{2}\right) \ , \notag \\
    &= \ \frac{ X \langle \phi , \phi \rangle_V }{2\langle \phi ,
    \phi \rangle_V^\frac{1}{2}} \ , \notag\\
    &= \frac{\langle D_X \phi,\phi \rangle_V
    + \langle\phi , D_X\phi \rangle_V}{2|\phi|} \ , \notag\\
    &= \ \Re \left( \langle\frac{\phi}{|\phi|}
    , D_X\phi\rangle_V \right) . \notag
\end{align}
Applying absolute values to both sides of the above equation and using the
Cauchy--Schwartz inequality, the desired result follows.
\end{proof}\ret

As a first application of estimate \eqref{Kato_ineq}, we prove the
following generalization of estimates \eqref{ext_slice_sob}--\eqref{char_slice_sob}
to complex scalar fields:\\

\begin{lem}[Weighted global Sobolev estimates for complex scalar
fields]\label{general_wighted_phi_sob}
Let  and let $\phi$ be a
section of $V$ over $\mathcal{M}$, and let $w_\gamma$ be the weight
function defined on line \eqref{w_gamma_def1}. Then on each time
slice $t=const.$, with $1 \leqslant t$, one has the estimates:
\begin{align}
        \sup_{t \leqslant 2r}  \ |\tp^{1+ \delta_+ }
    \tm^{\frac{1}{2}+\delta_-} \phi|^2 \, w_\gamma \ &\lesssim \
        \sum_{\substack{|I| \leqslant 1 , |J|\leqslant 2\\
        |I|+|J|\leqslant 2\\
        X\in \{\partial_r,S,\Omega_{0r}\} , Y\in \mathbb{O} }}
        \lp{\tp^{\delta_+}\tm^{\delta_-} w^\frac{1}{2}_\gamma\,
        D_X^I D_Y^J\phi}{L_x^2(t < 4r)}^2 \ ,
        \label{cov_ext_slice_sob}\\
        \sup_{2r<t}  \ |    \tp^{\frac{3}{2} + \delta_+ }
        \,     \phi |^2 \ &\lesssim \
        \sum_{\substack{|I| \leqslant 2 \\
        X\in \{S,\Omega_{0i}\} }}
        \lp{\tp^{\delta_+}\,
        D_X^I \phi }{L_x^2(r< \frac{3}{4}t )}^2 \ ,
        \label{cov_int_slice_sob}
\end{align}
Furthermore, defining the truncated cones:
\begin{align}
        \overline{C}(u) \ &= C(u)\cap\{1\leqslant t\leqslant t_0\}
        \cap\{ t < 2r \} \ ,\notag\\
        \td{C}(u) \ &= C(u)\cap\{1\leqslant t\leqslant t_0\}
        \cap\{ t < 4r \} \ ,\notag
\end{align}
we have the characteristic $L^\infty$ estimate:
\begin{equation}
        \sup_{(t,x)\in \overline{C}(u)}
        \ | \tp^{\frac{3}{2} + \delta_+ }\, \phi |^2 \ \lesssim \
        \sum_{\substack{|I| \leqslant 1 , |J|\leqslant 2 \\
        |I|+|J|\leqslant 2\\
        X\in \{ S - \Omega_{0r}\} ,  Y\in\mathbb{O}       }}
        \lp{\tp^{\delta_+}\,
        D_X^I D_Y^J \phi }{L_x^2\big(\td{C}(u)\big)}^2 \ .
        \label{cov_char_slice_sob}
\end{equation}
In the above estimates, $\mathbb{O}=\{\Omega_{ij}\}$
denotes the Lie algebra of the rotation group.
\end{lem}\ret

\begin{rem}
Notice that we have used the vector-field $S - \Omega_{0r}$ instead of
the individuals $\{S,\Omega_{0r}\}$ in the statement of estimate
\eqref{cov_char_slice_sob} above. This additional structure will be
important in the sequel.
\end{rem}\ret

\begin{proof}[Proof of estimates
\eqref{cov_ext_slice_sob}--\eqref{cov_char_slice_sob}]
The proof of these are virtually identical to the proof of
\eqref{ext_slice_sob}--\eqref{char_slice_sob}, with the added twist that
one uses \eqref{Kato_ineq} after each application of the scalar
estimates \eqref{ext_sob1}--\eqref{ext_sob2}, \eqref{int_sob},
and \eqref{char_sob1}--\eqref{char_sob2} to expand the scalar
derivatives to covariant derivatives. The details are left to the
interested reader.
\end{proof}\ret

It will also be useful for us to have a version of the estimate
\eqref{cov_ext_slice_sob} which is tailored to deal with the initial
data \eqref{initial_data2} in the context of the norms which
appear on the left hand side of \eqref{initial_smallness}.\\

\begin{lem}[Weighted Sobolev estimates for the initial data]
Let  and let $\phi$ be a section of the hyperplane bundle
$\{0\}\times \RR \times \mathbb{C}$ and let  $\frac{1}{2} < s_0$ be
a given parameter. Then the following weighted $L^\infty$ estimates
hold assuming that their right hand side is finite:
\begin{align}
        \sup_x \ (1+r)^{s_0 + \frac{1}{2}}\, |\phi| \ &\lesssim \
    \sum_{\substack{ 1\leqslant |I| \leqslant 2\\
        X\in\{\partial_i\}}}\
    \lp{(1+r)^{s_0 - 1 + |I|}    D_X^I \phi}{L^2(\RR)} \ ,
    \label{phi_0_Linfty}\\
    \sup_x \ (1+r)^{s_0 + \frac{3}{2}}\, |\phi| \ &\lesssim \
    \sum_{\substack{ |I| \leqslant 2\\
        X\in\{\partial_i\}}}\
    \lp{(1+r)^{s_0  + |I|}    D_X^I \phi}{L^2(\RR)} \ .
    \label{dotphi_0_Linfty}
\end{align}
\end{lem}\ret

\begin{proof}[Proof of the estimates
\eqref{phi_0_Linfty}--\eqref{dotphi_0_Linfty}]
The proof of these is almost identical to the proof of estimate
\eqref{cov_char_slice_sob} above. We simply apply the estimate
\eqref{ext_sob1}--\eqref{ext_sob2} of the appendix in order to the functions $\chi
\phi$, $\chi$ now being the cutoff on a spherical shell of dyadic
distance from the origin, using the Kato inequality
\eqref{Kato_ineq} after each step. Then, using the differential bounds:
\begin{equation}
        |\partial_i^I \chi| \ \lesssim \ (1+r)^{-|I|}\, \td{\chi} \ ,
         \notag
\end{equation}
where $\td{\chi}$ is some cutoff on the support of $\chi$, and adding
together over all the (finitely overlapping) $\chi$ we arrive at the
set of estimates:
\begin{align}
        \sup_x \ (1+r)^{s_0 + \frac{1}{2}}\, |\phi| \ &\lesssim \
    \sum_{\substack{ |I| \leqslant 2\\
        X\in\{\partial_i\}}}\
    \lp{(1+r)^{s_0 - 1 + |I|}    D_X^I \phi}{L^2(\RR)} \ ,
    \notag\\
    \sup_x \ (1+r)^{s_0 + \frac{3}{2}}\, |\phi| \ &\lesssim \
    \sum_{\substack{ |I| \leqslant 2\\
        X\in\{\partial_i\}}}\
    \lp{(1+r)^{s_0  + |I|}    D_X^I \phi}{L^2(\RR)} \ .
    \notag
\end{align}
The second of the above estimates is already of the form
\eqref{dotphi_0_Linfty}. To conclude \eqref{phi_0_Linfty} from the
first estimate above we simply need to eliminate the zero order
derivative term. This can be done through the use of the following
covariant Poincare type estimate which holds for $\frac{1}{2} < s_0$
assuming that $\phi$ is both smooth and compactly supported (which
we may do via a standard density argument):
\begin{equation}
        \int_{\RR}\ (1+r)^{2s_0 -2}\ |\phi|^2\ dx \ \lesssim \
    \int_{\RR}\ (1+r)^{2s_0}\ |D_r \phi|^2\ dx \ . \label{first_poincare}
\end{equation}
By an application of the Kato estimate \eqref{Kato_ineq}, this last
estimate follows from the corresponding statement with $\phi$
replaced by a real valued (smooth and compactly supported) test
function $\varphi$, and the covariant derivative $D_r$ replaced by
$\partial_r$. This in turn is achieved by integrating both sides of
the following identity with respect to the measure $dr\, d\omega$
and then applying a Cauchy-Schwartz:
\begin{multline}
        \partial_r\left( (1+r)^{2s_0 -1} r^2\, \varphi^2 \right) \ = \
    (2s_0 -1) (1+r)^{2s_0 -2} r^2 \, \varphi^2 \ + \\
     2 (1+r)^{2s_0 -1} r\, \varphi^2 \ + \
    2(1+r)^{2s_0 -1} r^2\, \varphi\,  \partial_r\varphi \ . \notag
\end{multline}
This completes the proof of \eqref{phi_0_Linfty}--\eqref{dotphi_0_Linfty}.
\end{proof}\ret

Before proceeding to the main estimate of this section, we first
prove another set of preliminary estimates which will be extremely
useful in the sequel. These are conjugated versions of the estimate
\eqref{first_poincare} above which also involve the optical weight
$\tm$, and which we also call Poincare estimates:\\

\begin{lem}[Covariant Poincare estimates for time slices]\label{ponc_lem}
Let $\phi$ be a section to $V$ over $\mathcal{M}$ which  has finite
$H^{1,\frac{p+q}{2}}$ norm (see line \eqref{scalar_initial_sob}),
and let $D_r$ denote the corresponding radial covariant derivative.
Then for constants $p,q$ such that both $|q| < p + 1 $ and $-1 < p$
one has the following estimate:
\begin{equation}
        \int_{\RR}\ \tm^p \tp^q  |\phi|^2 \ dx
        \ \lesssim \ \int_\RR \ \tm^{p + 2} \tp^q
        |\frac{1}{r} D_r (r\phi) |^2 \ dx \ , \label{radial_poincare}
\end{equation}
where the implicit constant depends on both $p$ and $q$.
More specifically, in the exterior region one has the estimate:
\begin{equation}
         \int_{t<r}\ \tm^p \tp^q  |\phi|^2 \ dx
        \ \lesssim \ \int_{t<r} \ \tm^{p + 2} \tp^q
        |\frac{1}{r} D_r (r\phi) |^2 \ dx \ , \label{radial_poincare1}
\end{equation}
whenever $0 < p +1 + q$ and $-1 <  p$, and in the interior one has:
\begin{equation}
         \int_{r<t}\ \tm^p \tp^q  |\phi|^2 \ dx
        \ \lesssim \ \int_{r<t} \ \tm^{p + 2} \tp^q
        |\frac{1}{r} D_r (r\phi) |^2 \ dx \ , \label{radial_poincare2}
\end{equation}
for values $q < p+1$ and $-1 <  p$.
\end{lem}\ret

\begin{proof}[Proof of estimate \eqref{radial_poincare}]
Keeping in mind the bound \eqref{Kato_ineq}, it is clear that we
only need to consider the incarnation of \eqref{radial_poincare}
with $\phi$ replaced by the real valued test function
$\varphi=|\phi|$ and $D_r$ replaced by the usual radial derivative
$\partial_r$. Notice that $\varphi$ is both bounded and has finite
$H^{1,\frac{p+q}{2}}$ (non-covariant) norm. We now consider the
regions $t<r$ and $r<t$ separately.\\

In the first case, it suffices to prove \eqref{radial_poincare} with
$\tm$ replaced by $ 1 - u$ and $\tp$ replaced by $1 + \bu $.
Furthermore, since we are assuming that the left hand side of
\eqref{radial_poincare1} is finite we will proceed via a density
argument that  assumes $\varphi$ is both smooth \emph{and} compactly
supported. This is because in order get things to work below we need
to use the truncated  mollification:
\begin{equation}
       \varphi_\epsilon \ = \  \chi_{{\epsilon}}\cdot f_\epsilon *
       \varphi , \label{varphi_mollification}
\end{equation}
where $\chi_{{\epsilon}}= \chi(\epsilon^{-1}\, \cdot)$ for some
$O(1)$ supported smooth unit function $\chi$, and $f_\epsilon =
\epsilon^{-3} f (\epsilon^{-1}\, \cdot)$ for some positive smooth
function $f$ with unit mass. Now, a standard argument shows this
mollification is such that $\varphi_\epsilon \to \varphi$ in the
space $H^{1,\frac{p+q}{2}}$. Therefore, we are reduced to showing
that for each fixed value of $\epsilon$ we have the bound
\eqref{radial_poincare1} for $\varphi_\epsilon$. In particular, we
may assume from the start that $\varphi$ is smooth and compact. With
this in mind, we integrate the identity:
\begin{multline}
        \partial_r\left( (1 - u)^{p+1} (1 + \bu)^q (r\varphi)^2\right)
    \ = \ (p+1)(1 - u)^{p} (1 + \bu)^q (r\varphi)^2\\
    + \ q (1 - u)^{p+1} (1 + \bu)^{q-1} (r\varphi)^2 \ + \
    2 (1 - u)^{p+1} (1 + \bu)^q \partial_r( r\phi) (r\phi) \ .
    \notag
\end{multline}
with respect to the measure $dr d\omega$ over the region $t<r$ to achieve:
\begin{multline}
        \int_{t<r} \ h^{p,q}(t,r) |\varphi|^2\ dx \ + \
    \int_{\mathbb{S}_{t=r}^2}\ (1 + 2r)^q |\varphi|^2
    dV_{\mathbb{S}_{t=r}^2 } \\ = \
    2 \int\int \ (1 - u)^{p+1} (1 + \bu)^q \partial_r( r\phi)
    (r\phi)\ dr d\omega \ , \label{phi_poincare_first_step}
\end{multline}
where $h^{p,q}(t,r)$ is the weight function given by:
\begin{align}
        h^{p,q}(t,r) \ &= \ (p+1)(1 - u)^{p} (1 + \bu)^q +
    q (1 - u)^{p+1} (1 + \bu)^{q-1} \ , \notag\\
    &= \ (1 - u)^{p} (1 + \bu)^q \left( (p+1) +
    q \frac{1-u}{1+\bu} \right) \ . \notag
\end{align}
The conditions $-1 < p$ and $-q < p+1$ now easily guarantee the
existence of a constant $C_{p,q}$ such that we have the point-wise
bounds:
\begin{equation}
        C_{p,q}^{-1} (1 - u)^{p} (1 + \bu)^q \ \leqslant \
    h^{p,q}(t,r) \ \leqslant \
    C_{p,q} (1 - u)^{p} (1 + \bu)^q  \ . \notag
\end{equation}
Applying this last bound and a Cauchy--Schwartz to the square of
\eqref{phi_poincare_first_step} with the angular integral on the left
hand side discarded, we arrive at the estimate:
\begin{multline}
        \left( \int_{t<r} \ (1 - u)^{p} (1 + \bu)^q |\varphi|^2\ dx\right)^2
    \\ \lesssim  \ \int_{t<r} (1 - u)^{p} (1 + \bu)^q|\varphi|^2\
    dx \ \cdot \ \int_{t<r} (1 - u)^{p+2} (1 + \bu)^q|\frac{1}{r}
    \partial_r(r\varphi)|^2\ dx \ . \notag
\end{multline}
This proves the assertion \eqref{radial_poincare1} by moving back to
the weights $\tp,\tm$. \\

We now prove the exterior estimate \eqref{radial_poincare2}. We will
again proceed via a density argument in which we assume $\varphi$ is
smooth\footnote{The only thing one needs to check here is that
$\lp{\chi_{\lesssim 1}\cdot  r^{-1} ( \varphi_\epsilon -
\varphi)}{L^2} \to 0$ as $\epsilon \to 0$, but this follows from the
vanilla Poincare estimate \eqref{basic_poincare} which holds in the
spaces $H^{1,k}$.} and hence bounded at the origin $r=0$. The
estimate follows essentially the same steps as above. First, we
replace the weights $\tp,\tm$ by $(1 + u)$ and $(1+\bu)$
respectively. The analog of \eqref{phi_poincare_first_step} now
reads:
\begin{multline}
        \int_{r<t} \ h^{p,q}(t,r) |\varphi|^2\ dx \ - \
    \int_{\mathbb{S}_{t=r}^2}\ (1 + 2r)^q |\varphi|^2
    dV_{\mathbb{S}_{t=r}^2 } \\ = \
    2 \int\int \ (1 + u)^{p+1} (1 + \bu)^q \partial_r( r\phi)
    (r\phi)\ dr d\omega \ , \label{int_phi_poincare_first_step}
\end{multline}
where this time  $h^{p,q}(t,r)$ is given by the expression:
\begin{equation}
        h^{p,q}(t,r) \ = \ (1 + u)^{p} (1 + \bu)^q \left( -(p+1) +
    q \frac{1+u}{1+\bu} \right) \ . \notag
\end{equation}
The bounds $-1 < p$ and $q < p+1$ now imply the existence of a
constant $C_{p,q}$ such that:
\begin{equation}
         C_{p,q}^{-1} (1 + u)^{p} (1 + \bu)^q \ \leqslant \
    -\, h^{p,q}(t,r) \ \leqslant \
     C_{p,q} (1 + u)^{p} (1 + \bu)^q  \ . \notag
\end{equation}
The remainder of the proof now proceeds by squaring and using
Cauchy-Schwartz as above.
\end{proof}\ret

\noindent We are now ready to state and prove the main estimate of this
section:\\

\begin{prop}[$L^\infty$ estimates for complex scalar fields]
Let $\phi$ be a section to $V$ over $\mathcal{M}$ with norm
$|\cdot|^2$ and compatible connection $D$. Let $F_{\alpha\beta}$
be the curvature tensor of $D$. Define the $k^{th}$ weighted
generalized energy content of $\phi$ in the time slab $0\leqslant t
\leqslant t_0$ to be:
\begin{equation}
        E^{(s,\gamma,\epsilon)}_k(0,t_0)[\phi] \ = \
    \sum_{\substack{|I|\leqslant k\\ X\in\mathbb{L}}}\
    E_0^{(s,\gamma,\epsilon)}(0,t_0)[D^I_X\phi] \ , \label{k_phi_eng}
\end{equation}
where $E_0^{(s,\gamma,\epsilon)}(0,t_0)[\phi]$ is defined on line
\eqref{phi_slab_eng}. Now define the time-slab $L^\infty$ type norm:
\begin{multline}
        \llp{\phi}{L^\infty[0,t_0](s,\gamma,\epsilon)}^2 \ = \
    \lp{\tp^{s+1} \tm^\frac{1}{2}
     (w')^\frac{1}{2}_{\gamma,\epsilon} \,
     \big(\frac{1}{r}D_L(r \phi) \chi_{t<2r} +
    D_L\phi \, \chi_{r<\frac{1}{2}t} \big)}{L^2(L^\infty)[0,t_0]}^2 \\
    + \ \tp^{2s+3}|\frac{1}{r}D_L(r\phi)\chi_{t<2r} +
    D_L\phi \, \chi_{r<\frac{1}{2}t}|^2\cdot w_\gamma \\
    \Big(   \tp^2 \tm^{2s-1} |\phi|^2 +
    \tp^2\tm^{2s+1} |D_{\bL}\phi|^2 +
    \tp^{2s+2}\tm\big(|\sD\phi|^2 + \sum_{\substack{|I|\leqslant 1\\
        X\in \mathbb{L}  }} |\frac{D_X^I\phi}{\tp}|^2\big)\Big)\cdot
    w_\gamma \ . \label{phi_Linfty_norm}
\end{multline}
Then recalling the definition of the energy norm \eqref{k_F_eng}
for the curvature $F_{\alpha\beta}$, we have the
following nonlinear estimate for the field $\phi$ under the additional
assumption  that $\epsilon < s-\frac{1}{2} $:
\begin{multline}
        \llp{\phi}{L^\infty[0,t_0](s,\gamma,\epsilon)}^2 \ \lesssim \
    E^{(s,\gamma,\epsilon)}_2(0,t_0)[\phi] \cdot\left( 1 +
      E_2^{(s,\gamma,\epsilon)}(0,t_0)[F] \right) \\
    + E_2^{(s,\gamma,\epsilon)}(0,t_0)[F]\cdot
    \llp{\phi}{L^\infty[0,t_0](s,\gamma,\epsilon)}^2
    \ . \label{phi_Linfty_est}
\end{multline}
\end{prop}\ret

\begin{proof}[Proof of estimate \eqref{phi_Linfty_est}]
To begin with we may assume that $1\leqslant t$, as the
complimentary case can easily be dealt with by a straightforward
application of radially weighted Sobolev estimates and energy bounds
similar to what we will use in the deep interior region $r <
\frac{1}{2}t$. We now split cases according to whether $t < 2r$ or
$r < \frac{1}{2}t$.\\

\subsection*{Estimate \eqref{phi_Linfty_norm} for $D_{\bL} \phi$ in the region $t < 2r$,
and for all components in the region $r < \frac{1}{2} t$} As in the
proof of Proposition \ref{main_ext_est} above, we will deal with the
worst component estimate together with the deep interior estimate.
These are both a consequence of the single fixed time bound (where
we leave the time dependence  on the left hand side implicit):
\begin{multline}
        \sup_{x} \
        \tp^{2} \tm^{2s+1} \big( |D \phi|^2 +
        \sum_{\substack{|I|\leqslant 1 \\ X\in \mathbb{L}}}\big|
         \frac{D_X^I \phi}{\tp} \big|^2 \big)\cdot w_\gamma
        \ \lesssim \\
        E_2^{(s,\gamma,\epsilon)}(0,t_0)[\phi]\cdot\left(
        1 + \llp{F}{L^\infty[0,t_0](s,\gamma,\epsilon)}^2\right) \
        + \ E_1^{(s,\gamma,\epsilon)}(0,t_0)[F]\cdot \llp{\phi}{
        L^\infty[0,t_0](s,\gamma,\epsilon)}^2 \
        . \label{general_weak_phi_decay}
\end{multline}
To prove this estimate, we first make a preliminary reduction.
Notice that by the homogeneity property of any $X\in \mathbb{L}$ we
easily have the bound:
\begin{equation}
        \sum_{\substack{|I|\leqslant 1 \\ X\in \mathbb{L}}}\big|
         \frac{D_X^I \phi}{\tp} \big|^2 \ \lesssim \
         |D \phi |^2 + \big| \frac{\phi}{\tp} \big|^2 \ .
         \label{simple_homog_bound}
\end{equation}
We will now prove the bound \eqref{general_weak_phi_decay} for the
right hand side of \eqref{simple_homog_bound} above. We begin by
using estimates \eqref{cov_ext_slice_sob} and
\eqref{cov_int_slice_sob} together with weights $\delta^+ = 0 ,
\delta^- = s$ and $\delta^+ =s$ (respectively) for the first term on
the right hand side of \eqref{simple_homog_bound}, and weights
$\delta^+ = -1 , \delta^- = s$ and $\delta^+ =s-1$ (respectively)
for the second term, which allows us to give the full bound:
\begin{multline}
        \sup_{x} \
        \tp^{2} \tm^{2s+1} \big(|D \phi|^2 + \big|\frac{\phi}{\tp}
        \big|^2
        \big)\cdot  w_\gamma \ \lesssim \\
        \sum_{\substack{|I|\leqslant 2 \\ X\in \mathbb{L}}}\
        \lp{\tm^s w_\gamma^\frac{1}{2}\, D_X^I D\phi }{L^2_x}^2
        \ + \ \sum_{\substack{|I|\leqslant 2 \\ X\in \mathbb{L}}}\
        \lp{\tp^{-1} \tm^s w_\gamma^\frac{1}{2}\, D_X^I \phi }{L^2_x}^2
        \ .
        \label{weak_precom_phi_Linfty}
\end{multline}
The second term on the right hand side of
\eqref{weak_precom_phi_Linfty} is immediately bounded by  R.H.S.
\eqref{phi_Linfty_est}, so we will concentrate on bounding the
second. To do this we will need to deal with the commutator $[D_X^I
, D]$. By using the formula \eqref{curvature_def} and the basic Lie
algebra formulas \eqref{lor_bracket_rel}, we arrive at the
point-wise bound:
\begin{equation}
        \sum_{\substack{|I| \leqslant 2\\
        X\in \mathbb{L}}}\ \big|[D_X^I , D] \phi\big| \ \ \lesssim
        \ \ \tp \cdot \sum_{\substack{|I| + |J| \leqslant 1\\
        X,Y\in \mathbb{L}}}\ |\mathcal{L}_X^I F|\cdot |D_Y^J \phi| \
        . \label{general_comm_phi_pointwise}
\end{equation}
Here we have used the notation $|\mathcal{L}_X^I F| =
|E(\mathcal{L}_X^I F)| + |H(\mathcal{L}_X^I F)|$. Using the
decomposition \eqref{tF_def} and the $L^\infty$ estimates
\eqref{alpha_bF_peel}--\eqref{other_bF_peel}, as well as the norm
definitions \eqref{F_Linfty_norm} and \eqref{phi_Linfty_norm} we can
fold estimate \eqref{general_comm_phi_pointwise} into
\eqref{weak_precom_phi_Linfty} above and produce the bound:
\begin{multline}
        \sum_{\substack{|I|\leqslant 2 \\ X\in \mathbb{L}}}\
        \lp{\tm^s w_\gamma^\frac{1}{2}\, D_X^I D\phi }{L^2_x}^2
        \ \lesssim \ E_2^{(s,\gamma,\epsilon)}[\phi]\cdot\left(
        1 + \llp{F}{L^\infty[0,t_0](s,\gamma,\epsilon)}^2\right)  \\
        + \ E_1^{(s,\gamma,\epsilon)}[F]\cdot \llp{\phi}{
        L^\infty[0,t_0](s,\gamma,\epsilon)}^2 \ + \
        \llp{F}{L^\infty[0,t_0](s,\gamma,\epsilon)}^2\cdot
        \sum_{\substack{|I|\leqslant 1\\ X\in \mathbb{L}}}
        \lp{ \tm^{-\frac{1}{2}} D_X^I\phi}{L^2(\RR)}^2 \
        . \notag
\end{multline}
The first two terms on the right hand side of this last expression
are exactly what we see on the right hand side of
\eqref{general_weak_phi_decay}. Therefore, we only need to bound the
last term. Unfortunately, this expression cannot be directly bounded
by the energy \eqref{k_phi_eng}. This is where we need to make use
of the fact that there is the additional room $0 < s-\frac{1}{2}$,
which allows us to replace the $\tm^{-\frac{1}{2}}$ weight with the
factor $\tm^{s -1}$. We now prove the  following estimate:
\begin{multline}
        \sum_{\substack{|I|\leqslant 1\\ X\in \mathbb{L}}}
    \lp{ \tm^{s-1} D_X^I\phi}{L^2(\RR)}^2
    \ \lesssim \ \sum_{\substack{|I|\leqslant 1\\ X\in \mathbb{L}}}
    \Big( \lp{ \tp^s (w)^\frac{1}{2}_\gamma \frac{1}{r}D_L(r
    D_X^I\phi)}{L^2(\RR)}^2 \\
    + \ \ \lp{ \tm^s (w)^\frac{1}{2}_\gamma D_{\bL} D_X^I\phi}{L^2(\RR)}^2 +
    \lp{  \tp^s (w)^\frac{1}{2}_\gamma \frac{1}{r} D_X^I\phi}{L^2(\RR)}^2 \Big)
    \label{Dr_poincare_applied} .
\end{multline}
This last line follows from Poincare estimate\footnote{We may assume
that the $H^{1,s-1}$ norm of $\phi$ and $D_X\phi$ is finite or there
is nothing to prove.} \eqref{radial_poincare} with $p = 2s -2 $ and
$q=0$ along with the simple point-wise bound (note that the extra
$w_\gamma$ factor is just added in to match with R.H.S.
\eqref{general_weak_phi_decay} and is not needed for things to be
sharp):
 \begin{equation}
         \tm^s | \frac{1}{r} D_r (r\phi)| \ \lesssim \ \left(
        \tp^s | \frac{1}{r} D_L(r\phi)| + \tm^s | D_{\bL}\phi|
        + \tp^s |\frac{\phi}{r}| \right)\cdot (w)^\frac{1}{2}_\gamma\ .
        \label{Dr_phi_bound}
\end{equation}
We have now completed our proof of estimate
\eqref{weak_precom_phi_Linfty} above.\\

\subsection*{Estimate \eqref{phi_Linfty_norm} for $\frac{\phi}{r}$ and
$\sD\phi$ in the region $t < 2r$} We begin with the undifferentiated
estimate for $\frac{\phi}{r}$. Applying estimate
\eqref{cov_ext_slice_sob} with $\delta_+ = s-1$ and $\delta_-=0$,
and using the point-wise bound:
\begin{equation}
        \sum_{X\in \{\partial_r , S, \Omega_{ij}, \Omega_{0r}\}}
    \ |D_X \phi|  \ \lesssim \
    \tp |\frac{1}{r} D_{L}(r\phi)| + \tm |D_{\bL} \phi| +
    \tp \left( |\sD \phi| + |\frac{\phi}{r}|\right) \ ,
    \label{phi_grad_bound}
\end{equation}
which follows easily from the expansions \eqref{lor_null_decomp},
we arrive at the $L^\infty$ estimate:
\begin{multline}
        \sup_{ 1\leqslant t < 2r}\
    \tp^{2s+2}\tm |\frac{\phi}{r}|^2w_\gamma \ \lesssim  \
    \sum_{\substack{ |I|\leqslant 1 \\
    Y\in \mathbb{O}}}\ \int_{t<4r} \ \Big(
    \tp^{2s}|\frac{1}{r} D_L(r D^I_Y \phi)|^2 \ + \
    \tm^{2s} |D_{\bL}(D^I_Y\phi)|^2 \\ + \
    \tp^{2s} \big( |\sD (D^I_Y\phi)|^2 + |\frac{D^I_Y \phi}{r}|^2 \big)\Big)
    w_\gamma \ dx \ . \label{phi_to_energy}
\end{multline}
Notice that one can add another derivative to $\phi$ on the
right hand side of \eqref{phi_to_energy} and still remain bounded by
the energy $E_2^{(s,\gamma,\epsilon)}[\phi]$. Therefore we also have
the desired bound for $\frac{D_X\phi}{r}$.\\

To bound the term $\sD\phi$ we use the following point-wise estimate
which is valid in the region $t\leqslant 2r$:
\begin{equation}
        \tp |\sD\phi| \ \lesssim \ \sum_{i<j}\ |D_{\Omega_{ij}} \phi|
        \ , \notag
\end{equation}
and then employ  estimate \eqref{phi_to_energy} with
$\frac{\phi}{r}$ replaced by $\frac{ D_{\Omega_{ij}}\phi}{r}$.\\

\subsection*{Improved $L^\infty$ estimates for the term
$\frac{\phi}{r}$ in the region $t<2r$} This estimate follows from a
simple combination of the global Sobolev and Poincare estimates as
we now explain. First, we see that by applying
\eqref{cov_ext_slice_sob} with $\delta_+=0$ and $\delta_-=s-1$, and
expanding out the fields $\partial_r$ and $\Omega_{0r}$, we have
that:
\begin{equation}
         \sup_{1\leqslant t < 2r} \tp^2 \tm^{2s -1} |\phi|^2 w_\gamma \
         \lesssim \ \sum_{\substack{|I| \leqslant 1 , |J|\leqslant
        2\\ |I|+|J| \leqslant 2\\
        X\in \{\partial_i,S,\Omega_{0i}\} , Y\in \mathbb{O} }} \
        \lp{\tm^{s-1}(w)^\frac{1}{2}_\gamma D^I_X D^J_Y \phi
        }{L^2(t<4r)} \ . \notag
\end{equation}
We are now finished by applying estimate \eqref{radial_poincare} in
the same way as used to prove \eqref{Dr_poincare_applied} above,
with the help of the  point-wise bound \eqref{Dr_phi_bound} to
convert the result into energy.\\

\subsection*{The $L^2(L^\infty)$ estimate for $\frac{1}{r}D_L(r\phi)$}
We'll begin here by estimating things in the region $1 < t < 2r$.
Using estimate \eqref{cov_ext_slice_sob} with $\delta_+ = s$ ,
$\delta_-=0$ and the weight $w_\gamma$ replaced by
$w'_{\gamma,\epsilon}$ on each fixed time slice, we directly have
the bound:
\begin{multline}
        \lp{\tp^{s+1} \tm^\frac{1}{2} (w')^\frac{1}{2}
    _{\gamma,\epsilon} \frac{1}{r} D_L (r\phi)}
    {L^2(L^\infty)[1,t_0]\cap\{t<2r\}} \ \lesssim \\
    \sum_{\substack{|I|\leqslant 1 , |J|\leqslant 2\\ |I|+|J|
    \leqslant 2\\ X\in\{\partial_r , S , \Omega_{0r}\}
    , Y\in \mathbb{O} }}\
    \lp{\tp^s (w')^\frac{1}{2}_{\gamma,\epsilon} D^I_X D^J_Y
    \frac{1}{r} D_L(r\phi)
    }{L^2(L^2)[1,t_0]\cap\{t<4r\}} \ . \label{DL_phi_L2Linfty_prelim}
\end{multline}
Keeping in mind now the space-time energy norm contained in
$E^{(s,\gamma,\epsilon)}_2[\phi]$, we are finished once we have
taken into account the commutator $[D_X^ID_Y^J , \frac{1}{r} D_L (r
\cdot)]$. Because we are dealing with the component with the best
decay properties, these commutators will need to be dealt with
carefully. First, we use the formulas \eqref{curvature_def} and
\eqref{lor_bracket_rel2} to compute the effect of the rotations (we
also use the decomposition \eqref{tF_def}):
\begin{subequations}\label{rot_alpha_comm_comp}
\begin{align}
        [D_{\Omega_{ij}}, \frac{1}{r} D_{L}(r\cdot)] \
    &= - \ \sqrt{-1} \Omega_{ij}^A \alpha_A  \ ,  \\
    \Omega_{kl}( [D_{\Omega_{ij}}, \frac{1}{r} D_{L} (r\cdot)] ) \ &= \ - \
    \sqrt{-1} \left(\Omega_{ij}^A \alpha_A(\mathcal{L}_{\Omega_{kl}}\tF)
    + \delta_{(lk} \Omega^A_{ij)} \, \alpha_A \right) \ .
\end{align}
\end{subequations}
Similarly, using the formulas \eqref{alpha_lie_component}  and the
fact that $\partial_r(\Omega_{ij}^A) = r^{-1} \Omega_{ij}^A$, we
compute the effect of the fields $\{\partial_r , S, \Omega_{0r}\}$
to be:
\begin{subequations}\label{alpha_S0r_comm_comp}
\begin{align}
        [D_r,\frac{1}{r} D_L(r\cdot)] \ &= \ - \ \sqrt{-1} \rho -
    \frac{1}{r^2}\ , \\
    [D_S,\frac{1}{r} D_L(r\cdot)] \ &= \ \sqrt{-1} u \, \rho
    \ - \ \frac{1}{r}D_L(r\cdot) \ , \\
    [D_{\Omega_{0r}}, \frac{1}{r} D_L(r\cdot)] \ &= \ \sqrt{-1}
    u \, \rho  \ + \ \frac{1}{r} D_L(r\cdot) \ + \
    \frac{u}{r^2} \ , \\
    \partial_r ( [D_{\Omega_{ij}},\frac{1}{r} D_L(r\cdot)] ) \ &= \
    -\ \sqrt{-1} \Omega_{ij}^A\left(
    \omega^k \alpha_A(\mathcal{L}_{\partial_k} \tF) + \frac{1}{r}
    \alpha_A \right)  \ , \\
    S( [D_{\Omega_{ij}}, \frac{1}{r}D_L(r\cdot)] ) \
    &= \ - \ \sqrt{-1}\Omega_{ij}^A
    \left( \alpha_A (\mathcal{L}_S \tF) - \alpha_A \right)\\
    \Omega_{0r}( [D_{\Omega_{ij}}, \frac{1}{r}D_L(r\cdot)] ) \ &= \
    - \ \sqrt{-1} \Omega_{ij}^A \left( \omega^k
    \alpha_A (\mathcal{L}_{\Omega_{0i}}\tF) - \frac{u}{r} \alpha_A\right)\ .
\end{align}
\end{subequations}
Combining the identities
\eqref{rot_alpha_comm_comp}--\eqref{alpha_S0r_comm_comp}, and
expanding out the fields $\partial_r$ and $\Omega_{0r}$ using the
identity $\bL(\omega_i)=0$ to move the contractions past the
$\frac{1}{r}D_L(r\cdot)$ derivative, we arrive at the following
point-wise bounds for the commutator in the region $1\leqslant t <
2r$:
\begin{multline}
        \sum_{\substack{|I| \leqslant 1 , |J|\leqslant
    2\\ |I|+|J| \leqslant 2\\
    X\in \{\partial_r,S,\Omega_{0r}\} , Y\in \mathbb{O} }}
    \ |D_X^I D_Y^J \frac{1}{r} D_{L}(r\phi)| \ \lesssim \
    \sum_{\substack{|I| \leqslant 1 , |J|\leqslant
    2\\ |I|+|J| \leqslant 2\\
    X\in \{\partial_i,S,\Omega_{0i}\} , Y\in \mathbb{O} }}\
    |\frac{1}{r} D_{L} (r D_X^I D_Y^J \phi)  | \ \\
    + \ \ \sum_{\substack{|I| + |J| \leqslant 1\\
    X,Y \in \mathbb{L}}}\
    \tp |\alpha(\mathcal{L}^I_X \tF)|\cdot |D^J_Y \phi|
    \ + \ \sum_{\substack{|I|\leqslant 1\\ X\in \mathbb{O}}}\
    \left( \tm |\rho | \cdot |D_X^I \phi | \
    + \  \tz |\frac{D_X^I\phi}{r}|\right) \ . \label{DL_comm_bound}
\end{multline}
We now insert each term on the right hand side of
\eqref{DL_comm_bound} into the right hand side of estimate
\eqref{DL_phi_L2Linfty_prelim} and bound the result using the norms
\eqref{k_F_eng}, \eqref{k_phi_eng}, \eqref{F_Linfty_norm}, and
\eqref{phi_Linfty_norm}. Since there are some choices involved and
the condition $\epsilon\leqslant s - \frac{1}{2}$ makes its
appearance here we do this explicitly. Note that the first term on
the right hand side of \eqref{DL_comm_bound} is already of the right
form so we move on to the second. For this term, we split cases
depending on wether the derivatives are falling on $\alpha$ or
$\phi$. In the first case, we estimate:
\begin{align}
        &\sum_{\substack{|I|\leqslant 1\\ X\in \mathbb{L}}}\
        \lp{\tp^{s+1} (w')^\frac{1}{2}_{\gamma,\epsilon} \cdot |\alpha
        (\mathcal{L}_X^I \tF)|\cdot
        \phi }{L^2(L^2)[0,t_0]}^2 \ , \notag\\
        \lesssim \ &\sum_{\substack{|I|\leqslant 1\\ X\in \mathbb{L}}}\
        \lp{\tp^{s} (w')^\frac{1}{2}_{\gamma,\epsilon} \ \alpha
        (\mathcal{L}_X^I \tF)}{L^2(L^2)}^2\cdot
        \lp{\tp \phi }{L^\infty(L^\infty)[0,t_0]}^2 \ , \notag \\
        \lesssim \ &E_1^{(s,\gamma,\epsilon)}(0,t_0)[F]\cdot
        \llp{\phi}{L^\infty [0,t_0](s,\gamma\epsilon)}^2
        \ . \notag
\end{align}
In the second case, we use H\"olders inequality
to conclude that:
\begin{align}
        &\sum_{\substack{|I|\leqslant 1\\ X\in \mathbb{L}}} \
        \lp{\tp^{s+1} (w')^\frac{1}{2}_{\gamma,\epsilon} \cdot |\alpha
        |\cdot D^I_X\phi }{L^2(L^2)[0,t_0]}^2 \ , \notag\\
        \lesssim \ &\lp{\tp^{s+1}\tm^\frac{1}{2}
        (w')^\frac{1}{2}_{\gamma,\epsilon} \ \alpha}
        {L^2(L^\infty)[0,t_0]}^2\cdot
        \sum_{\substack{|I|\leqslant 1\\ X\in \mathbb{L}}}\
        \lp{\tm^{-\frac{1}{2}} D_X^I\phi}{L^\infty(L^2)[0,t_0]}^2
        \ . \notag
\end{align}
Notice that the first factor on this last line is contained in the
norm \eqref{F_Linfty_norm}.
To bound the second factor on the last line above, we again use the
extra room $\frac{1}{2} < s$ and apply the estimate
\eqref{Dr_poincare_applied}. \\

It remains to deal with the last term on the right had side of
\eqref{DL_comm_bound}. To handle this we expand $\rho = \td{\rho} +
\overline{\rho}$. Using the identity \eqref{null_charge_tensor}, we
see that we may combine the term $\tm \overline{\rho}\, D^I_X\phi$
with the second term in this expression. Notice that this is
automatically included in the energy  \eqref{k_phi_eng}. Thus, it
suffices to use the bound:
\begin{multline}
        \sum_{\substack{|I|\leqslant 1\\ X\in \mathbb{L}}} \
        \lp{\tp^s \tm (w')^\frac{1}{2}_{\gamma,\epsilon}\
    \td{\rho} \, D_X^I \phi}{L^2(L^2)[0,t_0]}^2 \
    \lesssim \\
    \llp{F}{L^\infty[0,t_0](s,\gamma,\epsilon)}^2\cdot
    \sum_{\substack{|I|\leqslant 1\\ X\in \mathbb{L}}} \
    \lp{\tm^\frac{1}{2} (w')^\frac{1}{2}_{\gamma,\epsilon}\
     \, \frac{1}{r} D_X^I \phi}{L^2(L^2)[0,t_0]}^2 \ . \notag
\end{multline}
To control the second factor on the right hand side of this last
expression, we simply use the condition $\epsilon \leqslant s-
\frac{1}{2}$ which guarantees that $\tm \leqslant \tp^{2s} \tz^{1 +
2\epsilon}$. This concludes the $t<2r$ estimate for
$\frac{1}{r}D_L(r\phi)$ in $L^2(L^\infty)$.\\

To wrap things up for this subsection, we need to deal with the
$L^2(L^\infty)$ estimate in the deep interior region
$r<\frac{1}{2}t$. In this case, it is just as easy to deal with the
entire gradient $D\phi$ as we did in the first subsection of this
proof. Therefore, the bound we seek to prove is the following:
\begin{equation}
        \lp{\tp^{s+1 - \epsilon} D\phi
        }{L^2(L^\infty)[1,t_0](r<\frac{1}{2}t)}^2
        \ \lesssim \ E_2^{(s,\gamma,\epsilon)}(0,t_0)[\phi]
        \ . \label{far_int_L2Linfty}
\end{equation}
Using the  $r < \frac{1}{2}t$ range restriction in conjunction with
the point-wise bound:
\begin{equation}
        C^{-1}\ \sum_{X\in\mathbb{L}}\
        |D_X\phi| \ \leqslant \
        \tm |D\phi| \ \leqslant \ C\ \sum_{X\in\mathbb{L}}\
        |D_X\phi| \ , \label{double_exp_bound}
\end{equation}
as well as the following interior Sobolev estimate on each fixed
time slice (from \eqref{cov_int_slice_sob} with
$\delta^+=s-\frac{3}{2} - \epsilon$):
\begin{equation}
        \sup_{r < \frac{1}{2}t}\
        \sum_{X\in\mathbb{L}}\ \tp^{s-\epsilon} |D_X\phi|
        \ \lesssim \ \sum_{ \substack{1\leqslant |I|\leqslant 3\\
        X\in\mathbb{L}}}\ \lp{\tp^{s-\frac{3}{2} -\epsilon}D_X^I\phi}
        {L^2_x(r < \frac{3}{2}t)} \ , \notag
\end{equation}
we easily have fixed time bound:
\begin{equation}
        \sup_{r < \frac{1}{2} t}\
        \tp^{s+1 - \epsilon}\, |D\phi |
        \ \lesssim \ \sum_{ \substack{1\leqslant |I|\leqslant 3\\
        X\in\mathbb{L}}}\ \lp{\tp^{s-\frac{3}{2} -\epsilon}D_X^I\phi}
        {L^2_x(r < \frac{3}{2}t)} \ . \notag
\end{equation}
By integrating the square of this last estimate over the time
interval $(0,t_0)$, and again applying the point-wise bound
\eqref{double_exp_bound}, we arrive at the proof of
\eqref{far_int_L2Linfty} above. This completes our discussion of the
$L^2(L^\infty)$ estimate in \eqref{phi_Linfty_est}.\\

\subsection*{The pure $L^\infty$ for $\frac{1}{r}D_L(r\phi)$ in the
region $t<2r$} We begin here by applying estimate
\eqref{cov_char_slice_sob} to the term $\frac{1}{r}D_L(r\phi)$ with
$\delta_+=s$. This yields the bound:
\begin{multline}
        \sup_{(t,x)\in \overline{C}(u)} \ \tp^{2s+3}
    |\frac{1}{r}D_L(r\phi)|^2 w_\gamma \ \lesssim \\
    \sum_{\substack{|I|\leqslant 1 , |J|\leqslant 2\\ |I|+|J|
    \leqslant 2\\ X\in\{S - \Omega_{0r}\}
    , Y\in \mathbb{O} }}\
    \lp{\tp^s (w)^\frac{1}{2}_\gamma \, D^I_X D^J_Y \frac{1}{r} D_L(r\phi)
    }{L^2\big(\td{C}(u)\big)} \ . \label{char_DL_phi_L2Linfty_prelim}
\end{multline}
To deal with the commutator, we apply the identities
\eqref{rot_alpha_comm_comp}--\eqref{alpha_S0r_comm_comp} which provide
the following refinement of \eqref{DL_comm_bound}:
\begin{multline}
        \sum_{\substack{|I| \leqslant 1 , |J|\leqslant
    2\\ |I|+|J| \leqslant 2\\
    X\in \{S-\Omega_{0r}\} , Y\in \mathbb{O} }}
    \ |D_X^I D_Y^J \frac{1}{r} D_{L}(r\phi)| \ \lesssim \
    \sum_{\substack{|I| \leqslant 1 , |J|\leqslant
    2\\ |I|+|J| \leqslant 2\\
    X\in \{S,\Omega_{0i}\} , Y\in \mathbb{O} }}\
    |\frac{1}{r} D_{L} (r D_X^I D_Y^J \phi)  | \ \\
    + \ \ \sum_{\substack{|I| + |J| \leqslant 1\\
    X\in \mathbb{L}\\ Y\in \{S-\Omega_{0r}\}\cup\mathbb{O}}}\
    \tp |\alpha(\mathcal{L}^I_X \tF)|\cdot |D^J_Y \phi|
    \ + \ \sum_{\substack{|I|\leqslant 1\\ X\in \mathbb{O}}}\
    \left( |u|\cdot |\rho | \cdot |D_X^I \phi | \
    + \ |\frac{u D_X^I\phi}{r^2}|\right) \ . \label{char_DL_comm_bound}
\end{multline}
Notice that the weight $\tm$ appearing in \eqref{DL_comm_bound} is
replaced by $|u|$ in the last term on the right hand side of
\eqref{char_DL_comm_bound} above. This is because we do not have any
term involving the radial derivative $\partial_r$ to contend with on
the left hand side. It is necessary to have this if we are to use
the characteristic energy bound contained in \eqref{k_phi_eng}. We
now insert the terms on the right hand side of
\eqref{char_DL_comm_bound} into the right hand side of
\eqref{char_DL_phi_L2Linfty_prelim}. Notice that the first such
expression is automatically bounded by the energy \eqref{k_phi_eng}.
Therefore, it suffices to treat the last two groups of terms which
we
do in reverse order.\\

To bound the last group of terms on the right hand side of
\eqref{char_DL_comm_bound}, we expand $\rho = \td{\rho} +
\overline{\rho}$ and group the term $|u|\cdot |\overline{\rho}|\cdot
D_X^I\phi$ with the term $r^{-2} |u|D_X^I\phi$ using the identity
\eqref{null_charge_tensor}. The resulting term is contained in
\eqref{k_phi_eng}, so we are left with bounding:
\begin{align}
        &\sum_{\substack{|I| \leqslant 1\\ X \in \mathbb{L}}}\
        \lp{\tp^s \tm (w)^\frac{1}{2}_\gamma
        \cdot |\td{\rho}|\cdot D_X^I\phi }{L^2\big(\td{C}(u)\big)}^2
        \ , \notag\\
        \lesssim \
        &\lp{\tm^s (w)^\frac{1}{2}_\gamma
        \, \td{\rho}}{L^2\big(\td{C}(u)\big)}^2\cdot
        \sum_{\substack{|I| \leqslant 1\\ X \in \mathbb{L}}}\
        \lp{\tp^s\tm^{1-s}\, D_X^I \phi}{L^\infty[0,t_0]}^2 \ , \notag\\
        \lesssim  \ &E^{(s,\gamma,\epsilon)}[F]\cdot
        \llp{\phi}{L^\infty[0,t_0](s,\gamma,\epsilon)}^2 \ . \notag
\end{align}\ret

To wrap things up here, we need to bound the middle set of terms on
the right hand side of \eqref{char_DL_comm_bound} substituted into the
right hand side of \eqref{char_DL_phi_L2Linfty_prelim}. Dealing first
with the term where the derivatives fall on $\alpha$, we estimate:
\begin{align}
        &\sum_{\substack{|I|\leqslant 1\\ X\in \mathbb{L}}}\
        \lp{\tp^{s+1} (w)^\frac{1}{2}_\gamma \cdot |\alpha
        (\mathcal{L}_X^I \tF)|\cdot
        \phi }{L^2\big(\td{C}(u) \big)}^2 \ , \notag\\
        \lesssim \ &\sum_{\substack{|I|\leqslant 1\\ X\in \mathbb{L}}}\
        \lp{\tp^{s} (w)^\frac{1}{2}_\gamma \ \alpha
        (\mathcal{L}_X^I \tF)}{L^2\big(\td{C}(u) \big)  }^2\cdot
        \lp{\tp \phi }{L^\infty(L^\infty)[0,t_0]}^2 \ , \notag \\
        \lesssim \ &E_1^{(s,\gamma,\epsilon)}(0,t_0)[F]\cdot
        \llp{\phi}{L^\infty [0,t_0](s,\gamma\epsilon)}^2
        \ . \notag
\end{align}
Therefore, it remains to deal with the case where the derivatives fall
on $\phi$. Here we have to split cases depending on whether we have
$D_Y\phi$ with $Y = S - \Omega_{0r}$, or $Y\in\mathbb{O}$. In the
first case, we write $Y = \bu L$ and we expand:
\begin{equation}
        |\bu L \phi | \ \leqslant \ \tp |\frac{1}{r}D_L(r\phi)| +
    \frac{\tp}{r} |\phi| \ . \notag
\end{equation}
Using the fact that $t<2r$ this allows us to estimate:
\begin{align}
        &\lp{\tp^{s+1} (w)^\frac{1}{2}_\gamma\cdot
        |\alpha|\cdot D_{S- \Omega_{0r}} \phi}{L^2\big(\td{C}(u)\big)}^2
        \ , \notag\\
        \lesssim \ &\lp{\tp^{s} (w)^\frac{1}{2}_\gamma\
        \alpha}{L^2\big(\td{C}(u)\big)}^2\cdot
        \left( \lp{\tp^2 \frac{1}{r}D_L(r\phi) }{L^\infty[0,t_0]}^2
        + \lp{\tp \phi }{L^\infty[0,t_0]}^2\right) \ , \notag\\
        \lesssim \ &E_1^{(s,\gamma,\epsilon)}(0,t_0)[F]\cdot
        \llp{\phi}{L^\infty[0,t_0](s,\gamma,\epsilon)}^2 \ . \notag
\end{align}
We are now left with estimating the expression $\sum_{Y\in
\mathbb{O}}\, \tp^{s+1} (w)^\frac{1}{2}_\gamma\cdot |\alpha| \cdot
|D_Y\phi|$ in the space $L^2(\td{C}(u))$. This turns out to require
an argument which is a bit more involved than what we have been
doing so far. Notice that we cannot directly use either $L^\infty$
\emph{or} $L^2$ bounds on the $D_{\Omega_{ij}}\phi$ term.
Furthermore, a Poincare estimate similar to \eqref{radial_poincare}
for light-cones will not work in this situation because it
introduces boundary terms which cannot be controlled. Therefore,
being unable to deal with things on the light-cone $\td{C}(u_0)$, we
extend the integral to the exterior truncated time slab:
\begin{equation}
        \mathcal{R}(u_0,t_0) = \left([\max\{0,u_0\},t_0]\times\RR\right)
        \cap \{u \leqslant u_0\} \ . \notag
\end{equation}
Doing this via Stokes theorem yields the estimate:
\begin{multline}
        \sum_{Y\in\mathbb{O}}\
        \lp{\tp^{s+1} (w)^\frac{1}{2}_\gamma\cdot |\alpha| \cdot D_Y\phi
        }{L^2\big(\td{C}(u_0)\big)}^2 \
        \leqslant \\
        \sum_{Y\in\mathbb{O}}\ \Big(\ \Big|
        \int\int_{\mathcal{R}(u_0,t_0)}\
        \tp^{2s+2} \ \bL \left( | \alpha |^2\cdot |D_Y\phi
        |^2\right) \ dxdt \ \Big|\cdot w_\gamma(u_0)\\
        + \ \lp{\tp^{s+1} (w)^\frac{1}{2}_\gamma\cdot |\alpha|
        \cdot
        D_Y\phi}{L^2(\{t=u_0\}\times\RR)}^2 \Big)
        \ . \label{cone_to_slab}
\end{multline}
In the first term on the right hand side of the above expression,
the weight $w_\gamma(u_0)$ denotes the constant value of the
$w_\gamma$ weight along the $\td{C}(u_0)$ cone. In the second term
on the right hand side above we have switched to the variable
$w_\gamma$ weight, which is possible due to the fact that this is a
non-decreasing function for $r\to\infty$. To bound this second term
we estimate:
\begin{multline}
        \sum_{Y\in\mathbb{O}}\ \lp{\tp^{s+1} (w)^\frac{1}{2}_\gamma\cdot
        |\alpha| \cdot
        D_Y\phi}{L^2(\{t=u_0\}\times\RR)}^2 \ \lesssim \\
        \sum_{Y\in\mathbb{O}}\
        \llp{F}{L^\infty[0,t_0](s,\gamma,\epsilon)}^2\cdot
        \lp{\tp^{-\frac{1}{2}} \,
        D_Y\phi}{L^2(\{t=u_0\}\times\RR)}^2 \ . \notag
\end{multline}
We are now in the familiar territory where estimate
\eqref{Dr_poincare_applied} can be applied to the second factor on
the right hand side above.\\

Moving on to the first term on the right hand side of
\eqref{cone_to_slab}, we expand the $\bL$ derivative to bound:
\begin{multline}
        \sum_{Y\in\mathbb{O}}\ \Big| \
    \int\int_{\mathcal{R}(u_0,t_0)}\
    \tp^{2s+2} \ \bL \left( | \alpha|^2 \cdot |D_Y\phi
    |^2\right) \ dxdt \ \Big|\cdot w_\gamma(u_0)
    \ \leqslant \\
    \sum_{\substack{ Y\in\mathbb{O}\\
    X\in\{\partial_r , \partial_t , S , \Omega_{0r}\} }}\
    \int_0^{t_0}\int_{\RR}\
    \tp^{2s+2} \tm^{-1} \ \big| X \left( | \alpha|^2 \cdot |D_Y\phi
    |^2\right) \big| \cdot w_\gamma \ dxdt \ . \label{hard_int}
\end{multline}
Using the formulas \eqref{alpha_lie_component} and the fact that
$\partial_t(\alpha_A) = \alpha_A(\mathcal{L}_{\partial_t} \tF)$,
we expand the sum:
\begin{multline}
        \sum_{\substack{ Y\in\mathbb{O}\\
    X\in\{\partial_r , \partial_t , S , \Omega_{0r}\} }}\
    \big| X \left( | \alpha|^2 \cdot  |D_Y\phi
    |^2\right) \big| \ \lesssim \\
    \sum_{\substack{ |I|\leqslant 2\\
    X\in\mathbb{L} }}
    \ |\alpha |^2 \cdot |D_X^I \phi|^2 \ + \
    \sum_{\substack{ |I|\leqslant 1 \\
    X\in \mathbb{O} , Y\in \mathbb{L}
    }} \ |\alpha|\cdot|\alpha(\mathcal{L}_Y \tF)|\cdot
    |D_Y \phi|^2 \ .
        \label{alpha_phi_der_exp}
\end{multline}
Substituting the first term on the right hand side of
\eqref{alpha_phi_der_exp} in the right hand integral in
\eqref{hard_int}, and using the mixed norm estimate for $\alpha$
 we bound:
\begin{multline}
        \sum_{\substack{ |I|\leqslant 2\\
    X\in\mathbb{L} }} \ \lp{
    \tp^{s+1} \tm^{-\frac{1}{2}} (w)^\frac{1}{2}_\gamma\cdot
    |\alpha| \cdot D^I_X\phi}{L^2([0,t_0]\times\RR)}^2 \ \lesssim \\
    \llp{F}{L^\infty[0,t_0](s,\gamma,\epsilon)}^2\cdot
    \sum_{\substack{ |I|\leqslant 2\\
    X\in\mathbb{L} }} \ \lp{
    \tm^{-\frac{1}{2}+\epsilon}\
    D^I_X\phi}{L^\infty(L^2)[0,t_0]}^2 \ . \notag
\end{multline}
Using now the condition that $\epsilon \leqslant s-\frac{1}{2}$, we
are in a position to again apply \eqref{Dr_poincare_applied} noting
that one may add an extra derivative to this estimate and still remain
bounded in terms of $E_2^{(s,\gamma,\epsilon)}(0,t_0)[\phi]$.\\

We have now reduced things to estimating the integral:
\begin{equation}
        I \ = \ \sum_{X\in\mathbb{L} , Y\in\mathbb{O}}\
        \int_0^{t_0}\int_{\RR}\
        \tp^{2s+2} \tm^{-1} \ | \alpha|\cdot|\alpha(\mathcal{L}_X
        \tF)|\cdot |D_Y\phi|^2 \cdot w_\gamma \ dxdt \ . \notag
\end{equation}
Since this cannot be done directly, we use the identities:
\begin{align}
        |D_Y\phi|^2 \ &= \ 2 Y(\Re \langle\phi , D_Y \phi \rangle
        - 2 \Re \langle\phi , D^2_Y\phi \rangle\ , \notag \\
        2 \Re \langle\phi , D_Y \phi \rangle \ &= \
        Y(|\phi|^2) \ . \notag
\end{align}
to integrate by parts several times with respect to the $D_{\Omega_{ij}}$
derivatives to achieve the bound:
\begin{align}
        I \ &\lesssim \ \sum_{ \substack{|I| , |J| \leqslant 2\\
        X\in\mathbb{L} , Y\in\mathbb{O}}}\
        \int_0^{t_0}\int_{\RR}\
        \tp^{2s+2} \tm^{-1} \ | \alpha|\cdot|\alpha(\mathcal{L}_X^I
        \tF)|\cdot |\phi|\cdot |D^J_Y\phi| \cdot w_\gamma \ dxdt \notag \\
        & \ \ \ \ \ \ + \ \sum_{ \substack{
        |I|\leqslant 2\\  X\in \mathbb{L} }}\ \int_0^{t_0}\int_{\RR}\
        \tp^{2s+2} \tm^{-1} \ |\alpha(\mathcal{L}^I_X
        \tF)|^2 \cdot |\phi|^2 \cdot w_\gamma \
        dxdt  \ , \notag \\
        &= \ I_1 + I_2 \ . \notag
\end{align}
To estimate the first integral on the right hand side above, we use
a multiple H\"older inequality and the simple estimate
$\tm^{-1}w_\gamma \leqslant \tm^{2\epsilon} w'_{\gamma,\epsilon}$ to
conclude:
\begin{align}
        I_1 \ &\lesssim \ \sum_{ \substack{|I| , |J| \leqslant 2\\
        X\in\mathbb{L} , Y\in\mathbb{O}}}\
        \lp{\tp^{s+1} \tm^\frac{1}{2}
        (w')^\frac{1}{2}_{\gamma,\epsilon}\,
        \alpha}{L^2(L^\infty)[0,t_0]}
        \cdot \lp{\tp^s
        (w')^\frac{1}{2}_{\gamma,\epsilon}\,
        \alpha( \mathcal{L}_X^I \tF) }{L^2(L^2)[0,t_0]} \notag \\
        &\ \ \ \ \ \ \ \ \ \ \ \ \ \ \ \ \ \ \ \ \ \ \ \ \ \ \ \ \ \
        \ \ \ \
        \cdot \lp{\tp\tm^\epsilon \phi}{L^\infty[0,t_0]}\cdot
        \lp{\tm^{-\frac{1}{2} + \epsilon} D_Y^J
        \phi}{L^\infty(L^2)[0,t_0]}
        \ , \notag\\ \notag \\
        &\lesssim  \
        \left(\llp{F}{L^\infty[0,t_0](s,\gamma,\epsilon)}^2
        + E_2^{(s,\gamma,\epsilon)}[F]\right)
        \cdot\left(\llp{\phi}{L^\infty[0,t_0](s,\gamma,\epsilon)}^2
        + E_2^{(s,\gamma,\epsilon)}[\phi]
        \right) \ . \notag
\end{align}
Notice that the passage to the last line above is guaranteed by the
condition $\epsilon \leqslant s - \frac{1}{2}$ and the estimate
\eqref{Dr_poincare_applied}. To estimate the  second integral $I_2$
above, we proceed  as follows:
\begin{align}
        I_2 \ &\lesssim \ \sum_{ \substack{
        |I|\leqslant 2\\  X\in \mathbb{L} }}\
        \lp{\tp^s (w')^\frac{1}{2}_{\gamma,\epsilon}\ \alpha(
        \mathcal{L}_X^I \tF )}{L^2(L^2)[0,t_0]}^2
        \cdot\lp{\tp \tm^\epsilon \phi }{L^\infty [0,t_0]}^2 \ ,
        \notag\\
        &\lesssim \  E_2^{(s,\gamma,\epsilon)}[F]\cdot
        \llp{\phi}{L^\infty[0,t_0](s,\gamma,\epsilon)}^2 \ . \notag
\end{align}
This completes the proof of the pure $L^\infty$ bound for
$\frac{1}{r} D_L(r\phi)$ in the region $t<2r$.\\ \\

Combining all of the estimates proved in the various subsections
above, we have proved the estimate:
\begin{multline}
        \llp{\phi}{L^\infty[0,t_0](s,\gamma,\epsilon)}^2 \ \lesssim \
    E_2^{(s,\gamma,\epsilon)}[0,t_0][\phi]\cdot\left(
    1 + \llp{F }{L^\infty[0,t_0](s,\gamma,\epsilon)}^2
    + E_2^{(s,\gamma,\epsilon)}[0,t_0][F] \right) \\
    + \
    \left( \llp{F }{L^\infty[0,t_0](s,\gamma,\epsilon)}^2
    + E_2^{(s,\gamma,\epsilon)}[0,t_0][F] \right) \cdot
    \llp{\phi}{L^\infty[0,t_0](s,\gamma,\epsilon)}^2 \ . \notag
\end{multline}
Using now the $L^\infty$ estimate \eqref{main_F_linfty_est} for the curvature, we
arrive at estimate \eqref{phi_Linfty_est}.
\end{proof}

\ret
%-------------------------------------------------------------------------
%%%%%%%%%%%%%%%%%%%%%%%%%%%%%%%%%%%%%%%%%%%%%%%%%%%%%%%%%%%%%%%%%%%%%%%%%%
%-------------------------------------------------------------------------

\section{Abstract Weight Notation
and the Bilinear Space-Time Estimates in General Form}

In this section we will organize and consolidate the $L^2$ and
$L^\infty$ type estimates which have been proved in the preceding
sections. This will help to streamline notation for the remainder of
the paper, and will ultimately reduce a lot of overlapping which
occurs in the commutator estimates needed to deal with the nonlinear
problem \eqref{basic_MKG}. This will also give a chance for the
reader to review the various notations, parameters, weights, etc.
which have been introduces thus far. Our first order of business
here is to set up some generic markers for the quantities which
arise in the estimates we consider. These are the null components of
the curvature $F$ given by \eqref{null_decomp}, the corresponding
null decomposition of the gradient $D\phi$, and the weighted scalar
field $\frac{\phi}{\tp}$. As we have already mentioned in Remark
\ref{analogy_rem} above, there is a straightforward analogy between
the different members of these sets. To make this analogy precise we
form the following sets of field quantities:
\begin{subequations}\label{abs_def}
\begin{align}
        \Psi_{(1)} \ &= \ \{\alpha \ ,\  \frac{1}{r} D_L
        (r\phi)\, \chi_{t<2r} + D_L\phi\, \chi_{r<\frac{1}{2}t}   \} \ ,
        &\Psi_{(-1)} \ &= \ \{\balpha \ ,\  D_{\bL} \phi   \} \ ,
        \label{abs_def1}\\
        \Psi_{(0)} \ &= \ \{\td{\rho} \ ,\  \sigma \ , \
        \sD \phi \ ,\  \sum_{|J|\leqslant 1 , Y\in\mathbb{L}}
        \frac{D_Y^J \phi}{\tp} \}\ ,
        &Q_{(0)} \ &= \  \{ \overline{\rho}  \} \ . \label{abs_def2}
\end{align}
\end{subequations}
We will also use the following designation for sets containing the
component decompositions of Lie derivatives of the various field
quantities listed above:
\begin{subequations}\label{abs_def_lie}
\begin{align}
        \mathcal{L}_{\mathbb{L}}^k \Psi_{(1)} \ &= \
     \bigcup_{  \substack{ |I|\leqslant k\\ X\in\mathbb{L} }}\
     \{\alpha(\mathcal{L}_X^I \tF) \ ,\  \frac{1}{r} D_L
        (rD_X^I \phi)\, \chi_{t<2r} + D_L(D_X^I\phi)\,
    \chi_{r<\frac{1}{2}t}   \} \ , \label{abs_def_lie1}\\
    \mathcal{L}_{\mathbb{L}}^k
    \Psi_{(0)} \ &= \
    \bigcup_{  \substack{ |I|\leqslant k\\ X\in\mathbb{L} }}\
    \{\rho(\mathcal{L}_X^I \tF) \ ,\
    \sigma(\mathcal{L}_X^I \tF) \ , \
    \sD (D_X^I\phi) \ ,\  \sum_{|J|\leqslant 1 , Y\in\mathbb{L}}
    \frac{D_X^I D_Y^J \phi}{\tp} \}\ , \label{abs_def_lie2} \\
    \mathcal{L}_{\mathbb{L}}^k
    \Psi_{(-1)} \ &= \
    \bigcup_{  \substack{ |I|\leqslant k\\ X\in\mathbb{L} }}\
    \{\balpha( \mathcal{L}_X^I \tF ) \ ,
    \  D_{\bL} (D_X^I \phi)   \} \ , \label{abs_def_lie3}\\
    \mathcal{L}_{\mathbb{L}}^k Q_{(1)} \
    &= \
    \bigcup_{  \substack{ |I|\leqslant k\\ X\in\mathbb{L} }}\
    \{ \alpha(\mathcal{L}^I_X \overline{F} )  \} \ , \label{abs_def_lie4}\\
    \mathcal{L}_{\mathbb{L}}^k Q_{(0)} \
    &= \
    \bigcup_{  \substack{ |I|\leqslant k\\ X\in\mathbb{L} }}\
    \{ \rho(\mathcal{L}^I_X \overline{F} ) \ , \
    \sigma(\mathcal{L}_X^I \overline{F} )\} \ , \label{abs_def_lie5} \\
    \mathcal{L}_{\mathbb{L}}^k Q_{(-1)} \
    &= \
    \bigcup_{  \substack{ |I|\leqslant k\\ X\in\mathbb{L} }}\
    \{ \balpha(\mathcal{L}^I_X \overline{F} )  \} \ . \label{abs_def_lie6}
\end{align}
\end{subequations}
While it is true that the above objects are properly sets, they can
be manipulated as if they were functions by making each designator
stand for the absolute sum of each element in the set, and letting
products of the various sets denote the absolute sum of products of
each element in the corresponding sets. As a first application of
this principle we use the above notation to recast the norms
\eqref{k_F_eng}, \eqref{k_phi_eng}, \eqref{F_Linfty_norm}, and
\eqref{phi_Linfty_norm}:
\begin{align}
        E_k^{(s,\gamma,\epsilon)}(0,t_0)[\Psi] \ &= \
    E_k^{(s,\gamma,\epsilon)}(0,t_0)[F ] \ + \
    E_k^{(s,\gamma,\epsilon)}(0,t_0)[\phi] \ , \label{abs_k_eng}\\
    \llp{\Psi}{L^\infty[0,t_0](s,\gamma,\epsilon)}^2 \ &= \
    |q(F)|^2 \ + \ \lp{\tp^{s+1}\tm^\frac{1}{2}
    (w')^\frac{1}{2}_{\gamma,\epsilon}
    \Psi_{(1)}}{L^2(L^\infty)[0,t_0]}^2 \ \label{abs_Linfty_norm}\\
    & \hspace{-1in}
    \sup_{0\leqslant t \leqslant t_0}\ \left(\tp^2\tm^{2s-1} |\phi|^2 +
    \tp^{2s+3} |\Psi_{(1)}|^2 + \tp^{2}\tm^{2s+1}
    |\Psi_{(-1)}|^2    +   \tp^{2s + 2}\tm |\Psi_{(0)}|^2
    \right)\cdot w_\gamma \ . \notag
\end{align}
Notice that we have included the extra $L^\infty$ norm for $\phi$ on
the right hand side of \eqref{abs_Linfty_norm}. This extra estimate
turns out to be very important for us, but will be used in a somewhat
auxiliary manner. \\

We can now write out the content of line \eqref{abs_k_eng} with the
symbolic estimate (neglecting the characteristic energy terms which
are only used to prove $L^\infty$ estimates and which don't quite
fit into our scheme):
\begin{multline}
        E_0^{(s,\gamma,\epsilon)}(0,t_0)[\Psi] \ \geqslant \
        \sup_{0\leqslant t \leqslant t_0}\
    \int_{\{t\}\times\RR}\ \left(\tp^{2s} |\Psi_{(1)}|^2 +
    \tm^{2s}|\Psi_{(-1)}|^2 +
    \tp^{2s}|\Psi_{(0)}|^2  \right)\cdot w_\gamma \ dx \\
    + \ \int_0^{t_0}\int_{\RR}\
    \Big( \tp^{2s} |\Psi_{(1)}|^2 + \tz^{1 + 2\epsilon}\big(
    \tm^{2s}|\Psi_{(-1)}|^2 + \tp^{2s}|\Psi_{(0)}|^2\big)
    \Big)\cdot w'_{\gamma,\epsilon} \ dxdt  \ . \label{abs_eng_content}
\end{multline}
We can also write the $L^\infty$ estimates
\eqref{alpha_bF_peel}--\eqref{other_bF_peel},
\eqref{main_F_linfty_est}, and \eqref{phi_Linfty_est} in
the following consolidated symbolic form:\\

\begin{prop}[$L^\infty$ estimates for field quantities in abstract from]
Let $0 \leqslant t_0$ be a given fixed time, and let $0 < \gamma
,\epsilon , s$ be parameters chosen so that $\epsilon \leqslant
s-\frac{1}{2}$, then one has the following abstract non-linear $L^\infty$
estimate for the field quantities
$\mathcal{L}_{\mathbb{L}}^k \Psi$ and $\mathcal{L}_\mathbb{L}^k Q$:
\begin{multline}
        \llp{\mathcal{L}_{\mathbb{L}}^k\Psi}
    {L^\infty[0,t_0](s,\gamma,\epsilon)}^2 \ \lesssim \
    E_{k+2}^{(s,\gamma,\epsilon)}(0,t_0)[\Psi] \cdot\left( 1 +
    E_{2}^{(s,\gamma,\epsilon)}(0,t_0)[\Psi]\right) \\ + \
    E_{2}^{(s,\gamma,\epsilon)}(0,t_0)[\Psi]\cdot
    \llp{\mathcal{L}_{\mathbb{L}}^k \Psi}
    {L^\infty[0,t_0](s,\gamma,\epsilon)}^2
    \ , \label{abs_Linfty_est}
\end{multline}
and:
\begin{subequations}\label{abs_charge_Linfty}
\begin{align}
        \mathcal{L}_\mathbb{L}^k Q_{(1)} \ &\leqslant \ C_k \
        |q|\cdot\tz\tp^{-2} \cdot \chi_{t < r+1} \ ,
        \label{abs_charge_Linfty1}\\
        \mathcal{L}_\mathbb{L}^k Q_{(0)} \ , \ \mathcal{L}_\mathbb{L}^k Q_{(-1)}
        \ &\leqslant \ C_k \ |q|\cdot\tp^{-2} \cdot \chi_{t < r+1}
        \ . \label{abs_charge_Linfty2}
\end{align}
\end{subequations}
\end{prop}\ret

\noindent We now conclude this section by proving all of the
estimates we will encounter in the rest of the paper in a simple
abstract form. It turns out that everything we need can be cast in
the language of $L^2$ bilinear space-time estimates. These in turn
follow directly from H\"olders inequality and the $L^2$ and
$L^\infty$ type estimates contained in the right hand sides of
\eqref{abs_eng_content} and
\eqref{abs_Linfty_norm}--\eqref{abs_charge_Linfty2}.

\begin{prop}[Weighted $L^2$ bilinear estimates for field quantities in
abstract form] Let $0 \leqslant t_0$ be a given fixed time, and let
$0 < \gamma ,\epsilon , s$ be parameters chosen so that $\epsilon
\leqslant s-\frac{1}{2}$, and define the auxiliary weight:
\begin{equation}
        \td{w}_{\gamma,\epsilon} \ = \
    \tm^{2\gamma} \chi_{t<r} + \tm^{-\epsilon}\chi_{r<t} \ ,
    \label{bilin_weight}
\end{equation}
then one has the following abstract estimates for
field quantities $\Psi$ and $\Phi$:
\begin{align}
        \lp{\tp^{2s+1} \,
    \td{w}_{\gamma,\epsilon} \
    \Psi_{(1)}\cdot\Phi_{(0)}}{L^2[0,t_0]}^2 \ &\lesssim \
    \label{abs_L21}\\
    &\hspace{-.5in} \
    E_0^{(s,\gamma,\epsilon)}(0,t_0)[\Psi_{(1)}]\cdot
    \llp{\Phi_{(0)}}{L^\infty[0,t_0](s,\gamma,\epsilon)}^2 \ ,
    \notag \\
    \lp{\tp^{2s+1} \, \td{w}_{\gamma,\epsilon}
    \ \Psi_{(1)}\cdot\Phi_{(0)}}{L^2[0,t_0]}^2 \ &\lesssim \
    \label{abs_L22}\\
    &\hspace{-.5in} \
    E_0^{(s,\gamma,\epsilon)}(0,t_0)[\Phi_{(0)}]\cdot
    \llp{\Psi_{(1)}}{L^\infty[0,t_0](s,\gamma,\epsilon)}^2 \ ,
    \notag\\
    \lp{\tp^{2s + \frac{1}{2} - \epsilon } \tm^{\frac{1}{2} + \epsilon}
      \td{w}_{\gamma,\epsilon} \
      \Psi_{(0)}\cdot\Phi_{(0)}}{L^2[0,t_0]}^2 \
    &\lesssim \
    \label{abs_L23}\\
    &\hspace{-.5in} \
    E_0^{(s,\gamma,\epsilon)}(0,t_0)[\Psi_{(0)}]\cdot
    \llp{\Phi_{(0)}}{L^\infty[0,t_0](s,\gamma,\epsilon)}^2 \ ,
    \notag \\
    \lp{\tp^{s+1} \tm^{s}
    \td{w}_{\gamma,\epsilon}  \
    \Psi_{(1)}\cdot\Phi_{(-1)}}{L^2[0,t_0]}^2 \
    &\lesssim \
    \label{abs_L24}\\
    &\hspace{-.5in} \
    E_0^{(s,\gamma,\epsilon)}(0,t_0)[\Psi_{(1)}]\cdot
    \llp{\Phi_{(-1)}}{L^\infty[0,t_0](s,\gamma,\epsilon)}^2 \ ,
    \notag\\
    \lp{\tp^{s+1} \tm^s
     \td{w}_{\gamma,\epsilon} \
      \Psi_{(1)}\cdot\Phi_{(-1)}}{L^2[0,t_0]}^2 \
    &\lesssim \
    \label{abs_L25}\\
    &\hspace{-.5in} \
    E_0^{(s,\gamma,\epsilon)}(0,t_0)[\Phi_{(-1)}]\cdot
    \llp{\Psi_{(1)}}{L^\infty[0,t_0](s,\gamma,\epsilon)}^2 \ ,
    \notag\\
    \lp{\tp^{s + \frac{1}{2}  - \epsilon } \tm^{s + \frac{1}{2} + \epsilon}
      \td{w}_{\gamma,\epsilon} \
      \Psi_{(0)}\cdot\Phi_{(-1)}}{L^2[0,t_0]}^2 \
    &\lesssim \
    \label{abs_L26}\\
    &\hspace{-.5in} \
    E_0^{(s,\gamma,\epsilon)}(0,t_0)[\Psi_{(0)}]\cdot
    \llp{\Phi_{(-1)}}{L^\infty[0,t_0](s,\gamma,\epsilon)}^2 \ ,
    \notag\\
    \lp{\tp^{s + \frac{1}{2}  - \epsilon } \tm^{s + \frac{1}{2} + \epsilon}
      \td{w}_{\gamma,\epsilon} \
      \Psi_{(0)}\cdot\Phi_{(-1)}}{L^2[0,t_0]}^2 \
    &\lesssim \
    \label{abs_L27}\\
    &\hspace{-.5in} \
    E_0^{(s,\gamma,\epsilon)}(0,t_0)[\Phi_{(-1)}]\cdot
    \llp{\Psi_{(0)}}{L^\infty[0,t_0](s,\gamma,\epsilon)}^2 \ ,
    \notag
\end{align}
Recalling the definition of the weight function $w_{\gamma}$ from
line \eqref{w_gamma_def1}, we also have the following space-time
estimates for the interaction of the charge $Q$ and the field
quantity $\Psi$:
\begin{align}
        \lp{\tp^{s +\frac{5}{2} - \epsilon} \tm^{-1+\epsilon}
        (w_\gamma)^\frac{1}{2} Q_{(1)}\cdot \Psi_{(0)}}{L^2[0,t_0]} \ &\lesssim \ |q|^2\cdot
        E_0^{(s,\gamma,\epsilon)}[\Phi] \ , \label{cg_abs_L21}\\
    \begin{split}
        \lp{\tp^{s+\frac{3}{2} -\epsilon}\tm^\epsilon
        (w_\gamma)^\frac{1}{2}\ \big( Q_{(-1)}\cdot \Psi_{(0)} +\
        Q_{(0)}\cdot &\Psi_{(0)} +\
        Q_{(1)}\cdot \Psi_{(-1)}
        \big)}{L^2[0,t_0]} \\
         &\lesssim \ |q|^2\cdot
        E_0^{(s,\gamma,\epsilon)}[\Phi] \ ,
    \end{split}\label{cg_abs_L22}\\
        \lp{\tp^{\frac{3}{2} -\epsilon} \tm^{s+\epsilon}
        (w_\gamma)^\frac{1}{2} Q_{(0)}\cdot \Psi_{(-1)} }{L^2[0,t_0]} \ &\lesssim \ |q|^2\cdot
        E_0^{(s,\gamma,\epsilon)}[\Phi] \ . \label{cg_abs_L23}
\end{align}
\end{prop}\ret

\begin{rem}\label{norm_estimate_remarks}
Since the proofs of these are a direct consequence of the estimates
\eqref{abs_eng_content} and
\eqref{abs_Linfty_norm}--\eqref{abs_charge_Linfty2} we simply make a
few comments here. First of all, notice that the product estimates
\eqref{abs_L22} and \eqref{abs_L25} must use the extra mixed space
$L^2(L^\infty)$ norm on the right hand side of
\eqref{abs_Linfty_norm}. For example, if in the latter  one were
instead to use the pure $L^\infty$ estimate for $\Psi_{(1)}$ in
conjunction with the space-time $L^2$ bound for $\Phi_{(-1)}$, there
would be an additional need for a factor of $\tz^\epsilon$. As we
shall see in Section \ref{phi_error_sect}, such an extra convergence
factor is not available in our application where estimate
\eqref{abs_L25} comes up (although it is merely a notational
convenience where \eqref{abs_L22} arises). This is precisely the
reason we have included the extra $L^2(L^\infty)$ norms. Thus, as
far as our proof of Theorem \ref{main_th} is concerned we never need
to make use of the pure $L^\infty$ estimate for the quantity
$\Psi_{(1)}$.\\

The second thing we call the readers attention to here is the fact
that the estimates \eqref{cg_abs_L21}--\eqref{cg_abs_L23} will be
used in their precise form for $2t < r$. Notice that there are no
extra factors of $\tm^\gamma$ which can be put in these estimates as
was the case for \eqref{abs_L21}--\eqref{abs_L27} above. From a
technical point of view (aside from convenience), this is the reason
it is necessary to have the exact space-time energy norms on the
right hand side of \eqref{abs_eng_content}. Without these, estimates
of the form \eqref{cg_abs_L21}--\eqref{cg_abs_L23} would lead to
logarithmic divergences which we see no other way of controlling.
\end{rem}\ret

Before we finish, we take one last look at the estimates
\eqref{abs_L21}--\eqref{cg_abs_L23}. Our purpose is to recast these
in such a way that they conform more closely to how they will be
applied to estimate the system \eqref{basic_MKG}. In doing this, a
surprisingly simple and elegant picture emerges of the underlying
structure of the system \eqref{basic_MKG}. This reenforces our point
of view that the complex scalar field can be treated as if it were a
tensorial quantity with roughly the same properties as the curvature
$F_{\alpha\beta}$. We will put this structure to use through the so
called ``parity condition''. This is a numerical device used to keep
track of weights and components in contractions such as the right
hand side of \eqref{cbox_comm_formula}, which is typical of the kind
of error terms we treat in the next two sections. It turns out that
all of the estimates we will need can be reduced to a single
streamlined form, the motivation for which will become a bit more
clear through its use in the sequel.\\

\begin{prop}[Abstract parity form of the weighted bilinear $L^2$
estimates.] Let $0 \leqslant t_0$ be a given fixed time, and let $0
< \gamma ,\epsilon , s$ be parameters chosen so that $2\epsilon
\leqslant s-\frac{1}{2}$. Recall the definition of the auxiliary
weight $w_{\gamma,\epsilon}$ on line \eqref{w_gamma_def1.5}, and
define the parity weights $w_{\gamma,\epsilon}(a)$ via the formulas:
\begin{align}
        w_{\gamma,\epsilon}(1) \ &= \ \tz^{- \epsilon} \cdot
        (w)^\frac{1}{2}_{\gamma,\epsilon} \ ,
        &w_{\gamma,\epsilon}(0) \ &= \ \tz^{\frac{1}{2} } \cdot
        (w)^\frac{1}{2}_{\gamma,\epsilon} \
        &w_{\gamma,\epsilon}(-1) \ &= \ \tz^{2s - \epsilon} \cdot
        (w)^\frac{1}{2}_{\gamma,\epsilon} \ . \label{parity_weight}
\end{align}
Also define the parity optical weights:
\begin{align}
        \tau_{(1)} \ &= \ \tm \ , &\tau_{(0)} \ &= \ \tau_{(-1)} \ =
        \ \tp \ . \label{parity_optical}
\end{align}
Then the following abstract bilinear $L^2$ estimate holds for the
field quantities $\Psi ,\Phi$:
\begin{multline}
        \lp{\tp^{s-\frac{1}{2}} \ w_{\gamma,\epsilon}(a+b+c+d)\
        \tau_{(a)}\tau_{(b)}\ \Psi_{(c)}\cdot \Phi_{(d)}
        }{L^2[0,t_0]}^2 \ \lesssim \\
        E_0^{(s,\gamma,\epsilon)}(0,t_0)[\Psi]
        \cdot\llp{\Phi}{L^\infty[0,t_0](s,\gamma,\epsilon)}^2 \ .
        \label{abs_parity_est}
\end{multline}
whenever the condition $-1\leqslant a+b+c+d \leqslant 1$ holds.
One also has the analogous estimate for the interaction of $Q$ and
$\Phi$ under the same condition on $a,b,c,d$:
\begin{multline}
        \lp{\tp^{s-\frac{1}{2}} \ w_{\gamma,\epsilon}(a+b+c+d)\
        \tau_{(a)}\tau_{(b)}\ Q_{(c)}\cdot \Phi_{(d)}
        }{L^2[0,t_0]}^2 \ \lesssim \\
        |q|^2 \cdot
        E_0^{(s,\gamma,\epsilon)}(0,t_0)[\Phi]
        \ . \label{charge_abs_parity_est}
\end{multline}
\end{prop}\ret

\begin{proof}[Proof of the estimates
\eqref{abs_parity_est}--\eqref{charge_abs_parity_est}] We
concentrate on the first estimate \eqref{abs_parity_est}, the second
being similar and much easier. First of all, we claim that this can
be reduced to the more restricted special case:
\begin{multline}
        \lp{\tp^{s+\frac{1}{2}} \ w_{\gamma,\epsilon}(b+c+d)\
        \tau_{(b)}\ \Psi_{(c)}\cdot \Phi_{(d)}
        }{L^2[0,t_0]}^2 \ \lesssim \\
        E_0^{(s,\gamma,\epsilon)}(0,t_0)[\Psi]
        \cdot\llp{\Phi}{L^\infty[0,t_0](s,\gamma,\epsilon)}^2 \ ,
        \label{abs_parity_est_reduced}
\end{multline}
where now we impose the two conditions:
\begin{align}
        -1 \ &\leqslant c+d \ \leqslant \ 1 \ ,
        &-1 \ &\leqslant  \ b + c +d \ \leqslant \ 1 \ . \label{bcd_range_res}
\end{align}
Notice that \eqref{abs_parity_est_reduced} implies
\eqref{abs_parity_est} with these extra conditions enforced because
in that case one simply has the bound:
\begin{equation}
        w_{\gamma,\epsilon}(a+ b+c+d)\cdot\tp^{-1}\cdot\tau_{(a)} \ \leqslant \
        w_{\gamma,\epsilon}(b+c+d) \ , \notag
\end{equation}
as long as we are not in the case where both $a=1$ and $b+c+d=-1$.
If this happens, we see from the previous conditions that $b=-1,0$ so
we can use the bound:
\begin{equation}
        w_{\gamma,\epsilon}(a+ b+c+d)\cdot\tp^{-1}\cdot\tau_{(a)}
        \tau_{(b)}\ \leqslant \
        w_{\gamma,\epsilon}(b+c+d+1)\cdot \tau_{(b+1)}
        \ , \notag
\end{equation}
which holds in this case. Replacing now $b'=b+1$ we are back to the
case of \eqref{abs_parity_est_reduced}--\eqref{bcd_range_res}.\\

To derive the estimate \eqref{abs_parity_est} in the case where
$c+d=-1$ and $b=-1$, notice that one must then have $a=1$. In this
case the bound $\tp^{-1} \tau_{(1)} \leqslant 1$ reduces things to
\eqref{abs_parity_est_reduced} with $c+d=-1$ and $b=0$. The case
$c+d=1$, $b=1$, $a=-1$ can be treated similarly.\\

To derive  \eqref{abs_parity_est} in the case where  $c=d=-1$, we
use the symbolic inequality:
\begin{equation}
        \tp^{-1}\tau_{(1)} \Psi_{(-1)} \ \leqslant \ \Psi_{(0)} \ ,
        \notag
\end{equation}
which puts us in a position where we may again apply
\eqref{abs_parity_est_reduced} under the conditions
\eqref{bcd_range_res}.\\

To deal with \eqref{abs_parity_est} in the extreme case where
$c=d=1$, we simply use the fact that while it is not listed in the
estimates \eqref{abs_L21}--\eqref{abs_L27}, the product
$\Psi_{(1)}\cdot \Phi_{(1)}$ is so favorable that it satisfies the
space-time $L^2$ estimate \eqref{abs_parity_est} for any of the
weights $w_{\gamma,\epsilon}(a)$.\\

It now remains to prove \eqref{abs_parity_est_reduced} under the
conditions \eqref{bcd_range_res}. This involves a simple case by
case analysis split along the value of $c+d$:

\subsection*{Case: $c+d = 1$} In this case we can either have $b=-1$ or
$b=0$. In either case the $\tau_{(b)}$ weight is the same, so it
suffices to consider the case which maximizes the weight
$w_{\gamma,\epsilon}(b+c+d)$. In this case we are dealing with
$w_{\gamma,\epsilon}(1)$. Substituting this into estimate
\eqref{abs_parity_est_reduced}, and using the bound:
\begin{align}
        \tp^{s+\frac{1}{2}}\,
        w_{\gamma,\epsilon}(1)\, \tau_{(0)} \ &= \
    \tp^{s + \frac{3}{2} + \epsilon}
        \tm^{-\epsilon} (w)^\frac{1}{2}_{\gamma,\epsilon} \ , \notag\\
    &\lesssim \
        \tp^{2s+1} \td{w}_{\gamma,\epsilon} \ , \notag
\end{align}
which holds due to the condition $2\epsilon \leqslant
s-\frac{1}{2}$, we see that the desired result follows from
estimates \eqref{abs_L21}--\eqref{abs_L22} above.\\

\subsection*{Case: $c+d=0$} In this case, we can have all choices
$-1\leqslant b \leqslant 1$. It suffices to consider the one which
maximizes the product $w_{\gamma,\epsilon}(b)\cdot \tau_{(b)}$ which
is easily seen to be $b=0$. In this case we compute the total weight
in estimate \eqref{abs_parity_est_reduced} to be:
\begin{align}
        \tp^{s+\frac{1}{2}} w_{\gamma,\epsilon}(0)\, \tau_{(0)}
      \ &= \ \tp^{s+1}\tm^\frac{1}{2}
      (w)^\frac{1}{2}_{\gamma,\epsilon} \ , \notag\\
      &\lesssim \ \tp^{s+1} \tm^s \td{w}_{\gamma,\epsilon} \ , \notag\\
      &\lesssim \  \tp^{2s + \frac{1}{2} - \epsilon}
      \tm^{\frac{1}{2} + \epsilon} \td{w}_{\gamma,\epsilon} \ ,
      \notag
\end{align}
where the last two inequalities follow from the condition
$2\epsilon\leqslant s - \frac{1}{2}$. We are now in a position to
directly apply estimates \eqref{abs_L23}--\eqref{abs_L25}.\\

\subsection*{Case: c+d=-1} In this case we can have either $b=1$ or
$b=0$ so it suffices to consider the one which maximizes the product
$w_{\gamma,\epsilon}(b-1)\cdot \tau_{(b)}$. This turns out to depend
on the value of $\frac{1}{2} < s \leqslant 1$. If it is the case
that $\frac{3}{4} + \frac{\epsilon}{2} < s$, then the choice $b=1$
maximizes the product in which case we are dealing with the total
weight:
\begin{align}
        \tp^{s+\frac{1}{2}} w_{\gamma,\epsilon}(0)\, \tau_{(1)}
        \ &= \ \tp^{s} \tm^{\frac{3}{2}}
        (w)^\frac{1}{2}_{\gamma,\epsilon} \ , \notag\\
        &\lesssim \ \tp^{s+\frac{1}{2} - \epsilon}\tm^{s+\frac{1}{2} +
        \epsilon} \td{w}_{\gamma,\epsilon} \ . \notag
\end{align}
Notice that the inequality on the last line follows simply from the
condition $\frac{1}{2} < s$ and the fact that we may assume
$\epsilon \leqslant \frac{1}{2}$.\\

If on the other hand it is the case that $s\leqslant \frac{3}{4} +
\frac{\epsilon}{2}$, then the choice $b=0$ maximizes the product
$w_{\gamma,\epsilon}(b-1)\cdot \tau_{(b)}$ in which case we are
dealing with the total weight:
\begin{align}
        \tp^{s+\frac{1}{2}} w_{\gamma,\epsilon}(-1)\, \tau_{(0)}
        \ &= \ \tp^{\frac{3}{2} + \epsilon -s} \tm^{2s- \epsilon}
        (w)^\frac{1}{2}_{\gamma,\epsilon} \ , \notag\\
        &\lesssim \ \tp^{s+\frac{1}{2} - \epsilon}\tm^{s+\frac{1}{2} +
        \epsilon} \td{w}_{\gamma,\epsilon} \ , \notag
\end{align}
where, as before, the inequality is guaranteed by the condition
$2\epsilon\leqslant s - \frac{1}{2}$. Substituting this last line
into estimate \eqref{abs_parity_est_reduced} we can then directly
apply estimates \eqref{abs_L26}--\eqref{abs_L27}. This completes our
proof of \eqref{abs_parity_est}.
\end{proof}

\ret
%-------------------------------------------------------------------------
%%%%%%%%%%%%%%%%%%%%%%%%%%%%%%%%%%%%%%%%%%%%%%%%%%%%%%%%%%%%%%%%%%%%%%%%%%
%-------------------------------------------------------------------------

\section{Differentiating the Field Equations I: Error Estimates for the
curvature $F_{\alpha\beta}$}\label{F_error_sect}

We are now ready to begin in earnest our proof of Theorem \ref{main_th}. We
assume that we have fixed some level of regularity for the problem,
say $k$ derivatives with $2\leqslant k$. The result will be
demonstrated through a bootstrapping argument on the energy
\eqref{abs_k_eng}. Recalling the definition of the norms for the
initial data on the left hand side of \eqref{initial_smallness}, we define:
\begin{multline}
        E_k^{(s,\gamma)}(0)[F,\phi] \ = \
        \lp{E^{df}}{H^{k,s + \gamma}(\RR)}^2 \ + \ \lp{H}{H^{k,s  +
        \gamma}(\RR)}^2 \ +\\
        \lp{\underline{D}\phi_0 }{H^{k,s + \gamma}(\RR)}^2 \ + \
        \lp{\dot{\phi}_0}{H^{k,s+ \gamma}(\RR)}^2 \ . \label{initial_energy}
\end{multline}
We will now show that:\\

\begin{thm}[Bootstrapping form of  Theorem
\ref{main_th}]\label{main_th_bot_form} Let $k$ be a given level of
regularity, and assume that we are given parameters $0 <
s,\gamma,\epsilon$ with the properties that $s\leqslant 1$, and
$4\epsilon \leqslant s-\frac{1}{2}$, and $s+\gamma<\frac{3}{2}$. Let
$0 <T$ be a given time parameter, and let $\Psi$ denote the totality
of components of the system \eqref{basic_MKG} as defined by
\eqref{abs_def_lie} with the associated $k^{th}$ level energy
content \eqref{abs_k_eng}. Let $E^{(s,\gamma)}_k(0)[F,\phi]$ denote
the initial $k^{th}$ level energy content in the initial data
\eqref{initial_data}. Then there exists a constant $1 \leqslant
C_{k,s,\gamma,\epsilon}$ which depend only on $k,s,\epsilon,\gamma$,
but \emph{not} on $T$, such that if one first assumes that both:
\begin{equation}
        E_2^{(s,\gamma)}(0)[F,\phi] \ , \
    E_2^{(s,\gamma,\epsilon)}(0,T)[\Psi] \ \leqslant \
    C^{-1}_{k,s,\gamma,\epsilon} \ , \label{main_boot_assumption}
\end{equation}
then the following nonlinear estimate also holds:
\begin{equation}
        E_k^{(s,\gamma,\epsilon)}(0,T)[\Psi] \ \leqslant \
    C_{k,s,\gamma,\epsilon}\left( E_k^{(s,\gamma)}(0)[F,\phi]
    + \sum_{l=2}^7\ \left[E_k^{(s,\gamma,\epsilon)}(0,T)[\Psi]\right]^l
    \right) \ . \label{main_boost_est}
\end{equation}
\end{thm}\ret

\noindent In particular, we see from Theorem \ref{main_th_bot_form}
that the assumptions:
\begin{align}
        E_k^{(s,\gamma)}(0)[F,\phi] \ &\leqslant \
    \frac{1}{8} C^{-2}_{k,s,\gamma,\epsilon}  \ ,
    &E_k^{(s,\gamma,\epsilon)}(0,T)[\Psi] \ \leqslant \
    \frac{1}{2} C^{-1}_{k,s,\gamma,\epsilon} \ , \notag
\end{align}
together imply that:
\begin{equation}
        E_k^{(s,\gamma,\epsilon)}(0,T)[\Psi]\leqslant
    \frac{1}{4} C^{-1}_{k,s,\gamma,\epsilon} \ . \notag
\end{equation}
This, combined with the usual local existence theorem for the system
\eqref{basic_MKG}, and keeping in mind the main $L^\infty$ estimate
\eqref{abs_Linfty_est} (assuming that $C_{k,s,\gamma,\epsilon}$ is
chosen so large that the condition \eqref{main_boot_assumption}
allows one to absorb the extra $L^\infty$ terms on the right hand
side of this estimate) implies the claim of Theorem \ref{main_th}.
Therefore, for the remainder of this paper we will assume that the
solution exists up to time $T$ and we will concentrate on proving
the non-linear a-priori estimate \eqref{main_boost_est}.\\

In this section we concentrate on proving \eqref{main_boost_est} for
the portion of $\Psi$ which contains the curvature.  By the energy
estimate \eqref{main_F_L2_est2}, this boils down to estimating the
differentiated current vector $\mathcal{L}_X^I J$ at time $t=0$, as
well as over the space--time slab $[0,T]$. We do this separately.
For the latter it suffices to be able to prove the following
bounds:\\

\begin{prop}[Space-time error bounds for the curvature $F_{\alpha\beta}$]
Let $k$ be a given level of regularity, and assume that we are given
parameters $0 < s,\gamma,\epsilon$ with the properties that
$s\leqslant 1$,  and $2\epsilon \leqslant s-\frac{1}{2}$, and
$s+\gamma<\frac{3}{2}$. Let $0 <T$ be a given time parameter, and
let:
\begin{equation}
        J \ = \ \Im ( \phi \overline{D\phi}) \ , \notag
\end{equation}
be the current vector for the system \eqref{basic_MKG}. Recall
the current vector norm \eqref{J_norm}. Then there exists
a constant $1 \leqslant C_{k,s,\gamma,\epsilon}$ depending only on
$k,s,\epsilon,\gamma$, such that if:
\begin{equation}
        E_2^{(s,\gamma,\epsilon)}(0,T)[\Psi] \ \leqslant \
    C^{-1}_{k,s,\gamma,\epsilon} \ , \label{J_boot_assump}
\end{equation}
then one has the following weighted space-time estimate:
\begin{equation}
        \sum_{\substack{|I|\leqslant k\\
    X\in\mathbb{L}} }\
    \llp{\mathcal{L}_X^I J}{L^2[0,T](s,\gamma,\epsilon)}^2 \
    \lesssim \ \sum_{l=2}^5\
    \left[E_k^{(s,\gamma,\epsilon)}(0,T)[\Psi]\right]^l
    \ . \label{J_norm_est}
\end{equation}
\end{prop}\ret

\begin{proof}[Proof of estimate \eqref{J_norm_est}]
In light of the abstract $L^\infty$ estimate \eqref{abs_Linfty_est}
and the assumption \eqref{J_boot_assump} with the constant
$C^{-1}_{k,s,\gamma,\epsilon}$ chosen small enough that the right
hand side of \eqref{abs_Linfty_est} containing the extra $L^\infty$
norm can be absorbed into the left hand side, it suffices to be able
to show that:
\begin{multline}
        \sum_{\substack{|I|\leqslant k\\
        X\in\mathbb{L}} }\
        \llp{\mathcal{L}_X^I J}{L^2[0,T](s,\gamma,\epsilon)}^2 \
        \lesssim \\
        E_k^{(s,\gamma,\epsilon)}(0,T)[\Psi]\cdot
        \llp{\mathcal{L}_\mathbb{L}^{k-2} \Psi}
        {L^\infty[0,T](s,\gamma,\epsilon)}^2\left( 1 +
        \llp{\mathcal{L}_\mathbb{L}^{k-2} \Psi}
        {L^\infty[0,T](s,\gamma,\epsilon)}^2
        \right) \ , \label{J_norm_est_Linfty}
\end{multline}
where $\llp{\mathcal{L}_X^I J}{L^2[0,T](s,\gamma,\epsilon)}$ is the
notation from line \eqref{J_norm} above. Expanding out the norm on
the left hand side of this estimate, and using the weight notation
\eqref{parity_weight} we see that:
\begin{multline}
        \sum_{\substack{|I|\leqslant k\\
        X\in\mathbb{L}} } \ \llp{\mathcal{L}_X^I
        J}{L^2[0,T](s,\gamma,\epsilon)}^2
        \ = \  \sum_{\substack{|I|\leqslant k\\
        X\in\mathbb{L}} }\ \Big(
        \lp{\tp^{s+\frac{1}{2}}
        w_{\gamma,\epsilon}(1)\  (\mathcal{L}_X^I J)_L }{L^2[0,T]}^2 \\
        + \ \lp{\tp^{s+\frac{1}{2}} w_{\gamma,\epsilon}(0)\ \mathcal{L}_X^I
        \sJ }{L^2[0,T]}^2 \ + \
        \lp{\tp^{s+ \frac{1}{2}} w_{\gamma,\epsilon}(-1)\
        (\mathcal{L}_X^I J)_{\bL} }{L^2[0,T]}
        \Big) \ , \notag
\end{multline}
where  we are using the notation $|\mathcal{L}_X^I \sJ|^2 =
\delta^{AB} (\mathcal{L}_X^I J)_A\cdot (\mathcal{L}_X^I J)_B$.
Therefore, using the abstract parity estimates
\eqref{abs_parity_est} in the form of estimate
\eqref{abs_parity_est_reduced} as well as estimate
\eqref{charge_abs_parity_est}, it suffices to prove the following
symbolic bounds:
\begin{align}
        &\sum_{\substack{|I|\leqslant k\\
        X\in\mathbb{L}} }\ |(\mathcal{L}_X^I J)_L| \ \label{JL_parity_est}\\
        \lesssim \
        &\sum_{\substack{l+m = k-1\\ \hbox{or}\ l=0 , m=k \\ \hbox{and}\
        a+b+c = 1 }} \ \tau_{(a)} \ \left(
        \mathcal{L}^l_\mathbb{L} \Psi_{(b)} +
        \mathcal{L}^l_\mathbb{L} Q_{(b)}
        \right)\cdot \mathcal{L}^m_\mathbb{L} \Psi_{(c)}\cdot\left( 1 +
        \llp{\mathcal{L}_\mathbb{L}^{k-2}\Psi}{L^\infty[0,T](s,\gamma,\epsilon)}
        \right) \ , \notag\\
        &\sum_{\substack{|I|\leqslant k\\
        X\in\mathbb{L}} }\ |(\mathcal{L}_X^I J)_{\bL}| \ \label{JbL_parity_est}\\
        \lesssim \
        &\sum_{\substack{l+m =k-1\\ \hbox{or}\ l=0 , m=k \\
        \hbox{and}\
        a+b+c = -1 }} \ \tau_{(a)} \ \left(
        \mathcal{L}^l_\mathbb{L} \Psi_{(b)} +
        \mathcal{L}^l_\mathbb{L} Q_{(b)}
        \right)\cdot \mathcal{L}^m_\mathbb{L} \Psi_{(c)}\cdot\left( 1 +
        \llp{\mathcal{L}_\mathbb{L}^{k-2}\Psi}{L^\infty[0,T](s,\gamma,\epsilon)}
        \right) \ , \notag\\
        &\sum_{\substack{|I|\leqslant k\\
        X\in\mathbb{L}} }\ |\mathcal{L}_X^I \sJ| \ \label{sJ_parity_est}\\
        \lesssim \
        &\sum_{\substack{l+m =k-1\\ \hbox{or} \ l=0 , m=k \\ \hbox{and}\
        a+b+c = 0 }} \ \tau_{(a)} \ \left(
        \mathcal{L}^l_\mathbb{L}\Psi_{(b)}+
        \mathcal{L}^l_\mathbb{L} Q_{(b)}
        \right)\cdot \mathcal{L}^m_\mathbb{L} \Psi_{(c)}\cdot\left( 1 +
        \llp{\mathcal{L}_\mathbb{L}^{k-2}\Psi}{L^\infty[0,T](s,\gamma,\epsilon)}
        \right) \ . \notag
\end{align}
We will only prove the bounds
\eqref{JL_parity_est}--\eqref{sJ_parity_est} in the extended
exterior region $t<2r$. These bounds in the complimentary region
$r<\frac{1}{2}t$ follow from similar reasoning and are in fact much
simpler because the weights $\tau_{(a)}$ are all identical there.\\

Our first step is the following simple inductive calculation of the
Lie derivative $\mathcal{L}_X^I J$ which is based on repeated use of
the formula \eqref{basic_J_Lie}:
\begin{align}
        \mathcal{L}_X^I J \ &= \ \frac{1}{2}\cdot \sum_{\substack{K_1 + K_2 = I\\X_1 ,
        X_2 \in
        \mathbb{L}}}\ \left(D_{X_1}^{K_1}\phi \cdot \overline{ {\LC_{X_2}}^{K_2} D\phi}
        \right) + \left(D_{X_2}^{K_2}\phi \cdot \overline{ {\LC_{X_1}}^{K_1} D\phi}
        \right) \ , \label{J_lie_AB_calc}\\
        &= \ \frac{1}{2}\cdot  \sum_{\substack{K_1 + K_2 = I\\X,Y\in
        \mathbb{L}}}\ \left(D_{X_1}^{K_1}\phi \cdot \overline{ D( D_{X_2}^{K_2} \phi)}
        \right) + \left(D_{X_2}^{K_2}\phi \cdot \overline{ D( D_{X_1}^{K_1} \phi)}
        \right) \  \notag\\
        &\ \ \ \ \ \ \ \ \ \ \ + \ \frac{1}{2}\cdot \sum_{\substack{K_1 + K_2 = I\\X_1,X_2\in
        \mathbb{L}}}\ \left( D_{X_1}^{K_1}\phi \cdot \overline{
        [{\LC_{X_2}}^{K_2}  , D]\phi
        }\right) + \left(D_{X_2}^{K_2}\phi \cdot \overline{ [{\LC_{X_1}}^{K_1}  , D]\phi
        }\right) \ , \notag\\
        &= \ A + B \ . \notag
\end{align}
We now compute each of the $A$ and $B$ terms on the right hand side
of the above identities separately. Each of these terms can be seen
as real valued two forms, and we denote their components by
$A_\alpha$
and $B_\alpha$ respectively.\\

To compute the $A$ term, notice that since the sum is symmetric in
the $K_1$ and $K_2$ multiindices, we may replace the $D_L$
derivative by $\frac{1}{r} D_L(r\cdot)$ in the $A_L$ component.
Doing this, and putting absolute values around the different
components of $A$ while using the condition $|I|\leqslant k$ we have
the bounds:
\begin{align}
        |A_L| \ &\lesssim \  \sum_{\substack{|K| = k \\
        X\in\mathbb{L}}}\
    \tp\ |\frac{\phi}{\tp}|\cdot
    |\frac{1}{r}D_L(r D_X^K \phi)| \ \label{AF_bound1}\\
    &\ \ \ \ \ \ \ \ + \ \sum_{\substack{|K_1|+ |K_2| \leqslant k-1
        \\X_1 , X_2 ,Y\in \mathbb{L}}} \
    |D_Y D_{X_1}^{K_1}\phi |\cdot
        |\frac{1}{r}D_L(r D_{X_2}^{K_2} \phi )|
    \ , \notag\\
    |A_{\bL}| \ &\lesssim \ \sum_{\substack{|K| = k \\
        X\in\mathbb{L}}}\
    \tp\ |\frac{\phi}{\tp}|\cdot
    |D_{\bL}(D_X^K \phi)| \ \label{AF_bound2}\\
    &\ \ \ \ \ \ \ \ + \
    \sum_{\substack{|K_1|+ |K_2| \leqslant k-1 \\X_1 , X_2 ,Y\in
        \mathbb{L}}} \ |D_Y ( D_{X_1}^{K_1}\phi) |
    \cdot |D_{\bL} (D_{X_2}^{K_2}\phi )|
    \ , \notag\\
    |\sA | \ &\lesssim \
    \sum_{\substack{|K| = k \\
        X\in\mathbb{L}}}\
    \tp\ |\frac{\phi}{\tp}|\cdot
    |\sD (D_X^K \phi)| \ \label{AF_bound3}\\
    &\ \ \ \ \ \ \ \ + \
    \sum_{\substack{|K_1|+ |K_2| \leqslant k-1 \\X_1 , X_2 ,Y\in
        \mathbb{L}}} \ |D_Y( D_{X_1}^{K_1}\phi) |\cdot
        |\sD(D_{X_2}^{K_2}\phi )|
    \ . \notag
\end{align}
For the first term in each of the above sums, it suffices to merely
recall the designations \eqref{abs_def}--\eqref{abs_def_lie} and the
definition of the parity weights $\tau_{(a)}$ to achieve the bounds
\eqref{JL_parity_est}--\eqref{sJ_parity_est}. To achieve these
bounds for the second term on each line
\eqref{AF_bound1}--\eqref{AF_bound3} above, it suffices to show the
bounds:
\begin{equation}
        \sum_{\substack{|K_1| \leqslant l \\
        \\X_1 , Y\in \mathbb{L}}}\ |D_Y(D_{X_1}^{K_1}\phi )| \
        \lesssim \ \sum_{ a + b = 0 }\\
        \tau_{(a)}\ \mathcal{L}_\mathbb{L}^l \Psi_{(b)} \ . \notag
\end{equation}
This last line follows from a straight forward application of the
bounds:
\begin{align}
        |X^L| \ , \ |X^A| \ &\lesssim \ \tp \ ,
    &|X^{\bL}| \ \lesssim \ \tm \ , & &X\in\mathbb{L}
    \ , \label{lor_comp_bound}
\end{align}
which follows from an inspection of the identities \eqref{lor_null_decomp},
together with the designations
\eqref{abs_def}--\eqref{abs_def_lie} and the weight
definition \eqref{parity_optical}.\\

We now move on to dealing with the $B$ portion of identity
\eqref{J_lie_AB_calc}. This boils down to computing the commutator
actions $[{\LC_{X}}^K  , D]\phi$. By a repeated use of the formula
\eqref{D_complex_lie} and the fact that $\mathbb{L}$ is a Lie
algebra (so all its commutators involve constant coefficient sums),
we easily have the bounds:
\begin{align}
        |B_\alpha| \ \lesssim  \
    \sum_{\substack{|K_1| +|K_2| + |K_3| \leqslant |I|-1\\
    X_1,X_2,X_3 , Y \in \mathbb{L}}}\
    \ |(i_Y \mathcal{L}_{X_1}^{K_1}
    F)_\alpha|\cdot
    |D_{X_2}^{K_2}\phi|\cdot|D_{X_3}^{K_3}\phi| \ . \label{BF_bound}
\end{align}
We now use the fact that either $|K_2| \leqslant k-2$ or $|K_3|
\leqslant k-2$ which comes from the restriction given above on their
sum, to employ the $L^\infty$ estimate for $D_X^K\phi$ contained in
\eqref{abs_Linfty_est} to bound the product of the second and third
factor in this last sum as follows:
\begin{equation}
    |D_{X_2}^{K_2}\phi|\cdot|D_{X_3}^{K_3}\phi| \ \lesssim \
    \mathcal{L}_\mathbb{L}^{\max\{|K_2|,|K_3|\}}
    \Psi_{(0)}\cdot \lp{\mathcal{L}_\mathbb{L}^{k-2} \Psi}{
    L^\infty[0,T](s,\gamma,\epsilon)} \ . \notag
\end{equation}
Therefore, to achieve the right hand side of
\eqref{JL_parity_est}--\eqref{sJ_parity_est} for the $B$ term it
suffices to show the bounds:
\begin{align}
        \sum_{\substack{|K| \leqslant l \\
    \\X , Y \in \mathbb{L}}}\
    \ |(i_Y \mathcal{L}_{X}^{K} F)_L | \ &\lesssim \
    \sum_{ a + b = 1}\ \tau_{(a)}\left(
    \mathcal{L}_\mathbb{L}^{l}\Psi_{(b)} +
    \mathcal{L}_\mathbb{L}^{l}Q_{(b)}\right)
    \ , \label{F_i_parity_bound1}\\
    \sum_{\substack{|K| \leqslant l \\
    \\X , Y \in \mathbb{L}}}\
    \ |(i_Y \mathcal{L}_{X}^{K} F)_{\bL} | \ &\lesssim \
    \sum_{ a + b = -1}\ \tau_{(a)}\left(
    \mathcal{L}_\mathbb{L}^{l}\Psi_{(b)} +
    \mathcal{L}_\mathbb{L}^{l}Q_{(b)}\right) \ , \label{F_i_parity_bound2}\\
    \sum_{\substack{|K| \leqslant l \\
    \\X , Y \in \mathbb{L}}}\
    \ \sum_A \ |(i_Y \mathcal{L}_{X}^{K} F)_A | \ &\lesssim \
    \sum_{ a + b = 0}\ \tau_{(a)}\left(
    \mathcal{L}_\mathbb{L}^{l}\Psi_{(b)} +
    \mathcal{L}_\mathbb{L}^{l}Q_{(b)}\right)\ . \label{F_i_parity_bound3}
\end{align}
Splitting $\mathcal{L}_X^K F = \mathcal{L}_X^K \tF + \mathcal{L}_X^K
\overline{F}$ we see that it suffices to do this calculation for the
$\tF$ portion of things because the same computation for the charge
field $\overline{F}$ is identical. Since this is an abstract
counting argument for an arbitrary two form we can drop the Lie
derivatives and the tilde notation.  Making now the identifications:
\begin{align}
        L \ &\Leftrightarrow \ (1) \ ,
    &\bL \ &\Leftrightarrow \ (-1) \ ,
    &A,B \ &\Leftrightarrow \ (0) \ , \label{null_parity_numbers}
\end{align}
we see that in the groupings \eqref{abs_def}--\eqref{abs_def_lie}
the component $F_{\alpha\beta}$ is put in the set that corresponds
to the sum of $\alpha$ and $\beta$ (according to the numerical
designations \eqref{null_parity_numbers}) as they range over the
null frame $\{L,\bL,e_A\}$. Using this observation in conjunction
with the bound \eqref{lor_comp_bound} and the weight definitions
\eqref{parity_optical}, it is seen that the estimates
\eqref{F_i_parity_bound1}--\eqref{F_i_parity_bound3} are an
immediate consequence of translating the contractions $X^\beta
F_{\alpha\beta}$, computed in a null frame, into the more simple
numerical parity sum form. This completes the proof of the estimate
\eqref{J_norm_est}.
\end{proof}\ret

To complete the proof of the bootstrapping estimate
\eqref{main_boost_est} for the
electro-magnetic field $F_{\alpha\beta}$, it suffices to be able to
bound the initial data type norms  for the current vector $J$ given on
the right hand side of \eqref{main_F_L2_est2} in terms of the energy
\eqref{initial_energy}. To do this it is enough to show:\\

\begin{prop}[Initial data bounds for the system \eqref{basic_MKG} ]
Let $2\leqslant k$ be a given level of regularity, and assume that we are given
parameters $0 < \gamma$, and $\frac{1}{2} < s$, and $s+\gamma<\frac{3}{2}$.
Let $(F,\phi)$ be a
solution to the system \eqref{basic_MKG} with initial data
$(E,H,\phi_0,\dot{\phi}_0)$. Then there exists a constant $C_{k,s_0}$
such that if one first assumes the smallness condition:
\begin{equation}
        E_k^{(s,\gamma)}(0)[F,\phi] \ \leqslant \
        C^{-1}_{k,s,\gamma,\epsilon} \ , \label{data_smallness_cond}
\end{equation}
then one has the initial bounds:
\begin{multline}
         \sum_{|I|\leqslant k}
     \lp{(1+r)^{s+\gamma +|I|} \nabla_x^I
     J_0(0) }{L_x^\frac{6}{5}}^2 \  \
     + \ \ \sum_{|I|\leqslant k-1}\
     \lp{(1+r)^{s+\gamma +|I|+1} \nabla_{t,x}^I
     J(0) }{L_x^2}^2  \\ \lesssim \ \Big[ \
       E_k^{(s,\gamma)}(0)[F,\phi]\ \Big]^2 \ . \label{initial_J_bound}
\end{multline}
\end{prop}\ret

\begin{proof}[Proof of the estimate \eqref{initial_J_bound}]
We begin by bounding the first term on the left hand side of
\eqref{initial_J_bound}. At time $t=0$ we directly compute that in
terms of the initial data:
\begin{equation}
        \sum_{|I| = l} \ |\nabla_x^I J_0| \ \lesssim \
    \sum_{|K_1| + |K_2| = l}\ |\underline{D}^{K_1}\phi_0|\cdot
    |\underline{D}^{K_2} \dot{\phi}_0 | \ . \label{expand_der_J}
\end{equation}
This, combined with the two $L^\infty$ estimates:
\begin{align}
        (1+r)^{s+\gamma+\frac{1}{2}}\, |\phi_0| \ &\lesssim \
        \sum_{|I|\leqslant 1 }\ \lp{(1+r)^{s+\gamma +|I|} \underline{D}^I
        \underline{D} \phi_0 }{L^2(\RR)} \ , \label{uD_phi_Linfty1}\\
        (1+r)^{s+\gamma+\frac{3}{2}}\, |\dot{\phi}_0| \ &\lesssim \
        \sum_{|I|\leqslant 2 }\ \lp{(1+r)^{s+\gamma +|I|} \underline{D}^I
        \dot{\phi}_0 }{L^2(\RR)} \ , \label{uD_phi_Linfty2}
\end{align}
which in turn follow for the estimates
\eqref{phi_0_Linfty}--\eqref{dotphi_0_Linfty}, achieves the bound
\eqref{initial_J_bound} for this first term through a use of the
$L^2\cdot L^3$ H\"older inequality and directly integrating the
resulting $L^3$ estimate using
\eqref{uD_phi_Linfty1}--\eqref{uD_phi_Linfty2} and the fact that
$\frac{1}{2} < s+\gamma$. \\

To estimate the second term on the left hand side of
\eqref{initial_J_bound}, we first expand the derivatives
$\nabla_{t,x} J$ as we did in line \eqref{expand_der_J} above. Doing
this and using the $L^\infty$ estimate:
\begin{equation}
        \sum_{\substack{|K|\leqslant k-1 \\ X\in\{\partial_\alpha\}}}
        (1+r)^{1  + |K| } |D^K_X \phi| \ \lesssim \
        \sum_{\substack{1\leqslant |I_1| + |I_2| \leqslant k+1\\
        Y\in\{\partial_i\} ,
        X\in\{\partial_\alpha\}}}
        \ \lp{(1+r)^{s+\gamma -1 +|I_1| +|I_2| } D_Y^{I_1}
        D_X^{I_2} \phi }{L^2(\RR)}
        \ , \label{special_data_Linfty}
\end{equation}
which follows from the estimate \eqref{phi_0_Linfty} applied to
functions of the form $(1+r)^{ |K| }\, D^K_X \phi$ and the fact that
$1 < \frac{1}{2} + s + \gamma$ as well as $\frac{1}{2} < s +
\gamma$, we are reduced to having to demonstrate the general
estimate:
\begin{equation}
        \sum_{\substack{1\leqslant |I|\leqslant k+1\\ X\in\{\partial_\alpha\}}}
        \ \lp{(1+r)^{s+\gamma -1 +|I| } D_X^I \phi }{L^2(\RR)}^2 \
        \lesssim \ E_k^{(s,\gamma)}(0)[F] \ , \label{data_to_energy_bound}
\end{equation}
holds, provided that one first assumes the smallness condition
\eqref{data_smallness_cond} as well as the condition $2\leqslant
k$.\\

It is clear that the work in showing \eqref{data_to_energy_bound} is
simply a matter of eliminating the $D_t$ derivative in favor of
$D_i$ derivatives. This in turn boils down to controlling the
commutators of the form $[D_t , \underline{D}^I]$ because our goal
it to move the $D_t$ operator as far to the right as possible and
then to use the field equation \eqref{complex_field}. While this
overall strategy is simple to describe, the implementation becomes a
bit tedious because the commutator formula \eqref{curvature_def}
introduces the field electric $E_i$ which must then be
differentiated many times.\\

The simplest wa to deal with all of this seems to be to fix a value
for $k$ and then to induct on the number of $D_t$ derivatives on the
left hand side of this estimate. If this number is zero then the
claim is obvious, so we may now fix the number at some value
$1\leqslant l \leqslant k+1$. Our goal is to prove the bound:
\begin{equation}
        \lp{(1+r)^{s+\gamma -1 +|I| }
        D_X^{K_1} D_t D_Y^{K_2}\phi }{L^2(\RR)}^2 \ \lesssim \
        E_k^{(s,\gamma)}(0)[F] \ , \label{induct_bound_to_prove}
\end{equation}
where $|K_1| + |K_2| \leqslant |I|-1 $ and $X\in\{\partial_\alpha\}$
while $Y\in\{\partial_i\}$. By induction, we may assume that the
same estimate holds with the $D_t$ derivative  replaced by one of
$D_i$. In this regard it will be useful for us to employ the
following notation:\ If $K$ is a multiindex, then we denote by
$|K^0|$ the number of time derivatives in the operator $D_X^K$ (note
that we are only working with $X\in\{\partial_\alpha\}$ now). With
this convention, we have that $|K^0_1|\leqslant l-1$ in
\eqref{induct_bound_to_prove} above. Now, computing some multiple
commutators we  have the point-wise bound:
\begin{multline}
        |D_X^{K_1} D_t D_Y^{K_2}\phi | \ \ \lesssim  \ \
        |D_X^{K_1} D_Y^{K_2} D_t \phi | \\ + \
        \sum_{\substack{|I_1| + |I_2| = |K_2| -1\\
        X\in\{\partial_\alpha\}\\
        Y_1, Y_2 , Z \in \{\partial_i\}}}\
        \big|D_X^{K_1} D^{I_1}_{Y_1} [D_t, D_Z] D^{I_2}_{Y_2}\phi \big|
        \ \ . \label{commutator_line}
\end{multline}
We estimate the second term in the above expression first.
Substituting this into the left hand side of
\eqref{induct_bound_to_prove} and expanding out the single
commutator in terms of the electric field $E_i$, distributing the
left hand derivatives according to the covariant Leibnitz rule,
using the field equations \eqref{EH_Maxwell} to remove time
derivatives from $E$, and splitting result along the Hodge
decomposition $E = E^{df} + E^{cf}$ we have the bound:
\begin{align}
        &\sum_{\substack{|K_1^0| \leqslant l-1\\
        |K_2| + |I_1| + |I_2| \leqslant k-1\\
        X\in\{\partial_\alpha\}\\
        Y_1, Y_2 , Z \in \{\partial_i\}}}\
        \lp{(1+r)^{s+\gamma +1 + |I_1| + |I_2|}
        D_X^{K_1} D^{I_1}_{Y_1} [D_t, D_Z]
        D^{I_2}_{Y_2}\phi}{L^2(\RR)} \ , \label{big_initial_bound}\\
        \lesssim \
        &\sum_{\substack{|I_1| + |I_2| \leqslant k-1\\
        |I_2^0| \leqslant l-1\\
        X\in \{\partial_\alpha\} }} \
        \Big( \big( \ \lp{(1+r)^{s+\gamma +1 + |I_1| + |I_2|}
        (\nabla_{x}^{I_1} E^{df} )\cdot (D_X^{I_2} \phi )
        }{L^2(\RR)} \notag\\
        &\ \ \ \ \ \ \ \ \ \ \ \ \ \ \
        \ \ \ \ \ \ \ \ \ \ \ \ \ \ \
        + \
        \lp{(1+r)^{s+\gamma +1 + |I_1| + |I_2|}
        (\nabla_{x}^{I_1} H )\cdot (D_X^{I_2} \phi )
        }{L^2(\RR)}\ \big) \notag \\ \notag \\
        &\ \ \ \ \ \ \ \ \ \ \ \ \ \ \ + \
        \lp{(1+r)^{s+\gamma +1 + |I_1| + |I_2|}
        (\nabla_{x}^{I_1} E^{cf} )\cdot (D_X^{I_2} \phi )
        }{L^2(\RR)}\ \Big) \notag \\ \notag \\
        &\ \ \ \ \ \ \ \ \ \ \ \ \ \ \ \ \ \ \ \ \ \
        + \ \sum_{\substack{|I_1| + |I_2| \leqslant k-2 \\
        |I_1^0| + |I_2^0| \leqslant l-1\\
        X\in\{\partial_i\}}}
        \lp{(1+r)^{s+\gamma  + 2 + |I_1| + |I_2|}
        (\nabla_{t,x}^{I_1} J)\cdot (D_X^{I_2} \phi ) }{L^2(\RR)} \
        , \notag\\
        = \ &A + B + C \ . \notag
\end{align}
We now bound each of the three terms on the right hand side above
separately. For the  $A$ term we are done by induction after a
simple use of the $L^\infty$ estimate \eqref{special_data_Linfty}.\\

To bound the $B$ term, we add and subtract off the term $\nabla_x
\frac{1}{2\pi r}\cdot \chi^+(r-2)$ from the curl-free component
$E^{cf}$. Then using the estimate \eqref{gen_weighted_grad_sob} in
conjunction with the bounds we have already established for the
first term on the left hand side of \eqref{initial_J_bound}, and by
making another use of the estimate \eqref{special_data_Linfty} this
we are reduced to showing that:
\begin{equation}
        \sum_{\substack{ |I_2| \leqslant k-1\\
        |I_2^0| \leqslant l-1\\
        X\in \{\partial_\alpha\} }} \
        \lp{(1+r)^{s+\gamma - 1 + |I_2|}
        D_X^{I_2} \phi }{L^2(\RR)}^2 \ \lesssim \
        E_k(0)^{(s,\gamma)}[\phi,F] \ . \notag
\end{equation}
But this last estimate follows from our inductive hypothesis unless
$I_2=0$ in the sum. In this latter case we can use the Poincare
estimate \eqref{first_poincare} to make our claim.\\

Our next task here is to prove the (inductive) bound
\eqref{initial_J_bound} for the $C$ term on the right hand side of
line  \eqref{big_initial_bound} above. As before, we use the
estimate \eqref{special_data_Linfty} to deal with the factor
involving $(1+r)^{1 + |I_2|} D_X^{I_2} \phi$. Doing this and making
use of the continuity equation \eqref{J_cont} and the smallness
condition \eqref{data_smallness_cond}, we arrive at the bound (as
long as $0<l-1$):
\begin{equation}
        C \ \lesssim \ \sum_{\substack{|I_1| \leqslant k-2 \\
        |I_1^0| \leqslant l-2 }}
        \lp{(1+r)^{s+\gamma  + 1 + |I_1| }
        \nabla_{t,x}^{I_1} J\, (0) }{L^2(\RR)} \ . \notag
\end{equation}
We can now reduced things to the point where we can go back to the
paragraph containing line \eqref{special_data_Linfty} above and
repeat the entire argument up to this point which allows one to
bound the right hand side of this last expression via induction.\\

To finish here, we only need to deal with  the first term on the
left hand side of \eqref{commutator_line}. First of all, if it is
the case that $|K^0_1| = 0$ then we are automatically done. If
however $1\leqslant |K^0_1|$ then we simply repeat the steps
outlined above for the second term on the right hand side
\eqref{commutator_line}, but this time applied to the second term on
the right hand side of the bound:
\begin{multline}
            |D_X^{\td{K}_1} D_t D_Y^{K_2} D_t \phi | \ \ \lesssim  \ \
    |D_X^{\td{K}_1} D_Y^{K_2} D^2_t \phi | \\ + \
    \sum_{\substack{|I_1| + |I_2| = |K_2| -1\\
    X\in\{\partial_\alpha\}\\
    Y_1, Y_2 , Z \in \{\partial_i\}}}\
    \big|D_X^{\td{K}_1} D^{I_1}_{Y_1} [D_t, D_Z] D^{I_2}_{Y_2}D_t \phi \big|
    \ \ , \label{comm_line2}
\end{multline}
where this time $|\td{K}_1| = |K_1|-1$ and $\td{K}^0_1 = l-2$.
Notice that these numerical conditions guarantee that no more than a
total of $k-1$ derivatives altogether, and no more than $l-1$  of
the $D_t$ derivatives in particular, can fall on $\phi$ in this
second term.  Therefore, we are again left with estimating terms
which are similar to those appearing on the right hand side of
\eqref{big_initial_bound} above.\\

Finally, to deal with the first term on the right hand side of
\eqref{comm_line2} we simply use the wave equation
\eqref{complex_field} which trades the derivatives $D^2_t$ for the
covariant Laplacian $D^i D_i$. Having reduced the total number of
time derivatives in this last expression we are done by induction.
This completes the proof of the estimate \eqref{initial_J_bound}.
\end{proof}\ret

\ret
%-------------------------------------------------------------------------
%%%%%%%%%%%%%%%%%%%%%%%%%%%%%%%%%%%%%%%%%%%%%%%%%%%%%%%%%%%%%%%%%%%%%%%%%%
%-------------------------------------------------------------------------

\section{Differentiating the Equations II: Error Estimates for the
complex scalar field $\phi$}\label{phi_error_sect}

To wrap things up, we need to prove the boot-strapping estimate
\eqref{main_boost_est} for the $\phi$ portion of things. Using the
energy estimate \eqref{main_phi_L2_est} for complex scalar fields,
and using the assumption that $C^{-1}_{k,s,\gamma,\epsilon}$ in
\eqref{J_boot_assump} is chosen so small that the $L^\infty$ norm
term on the right hand side of  \eqref{main_phi_L2_est} can be
absorbed into the right hand side of that inequality, and using the
estimates \eqref{first_poincare} and \eqref{data_to_energy_bound}
above to deal with the terms involving the initial data we see that
it is enough to prove the following $L^2$ estimate for the
commutator of covariant derivatives
and the complex D'Lambertian:\\

\begin{prop}[Error bounds for the complex scalar field $\phi$]
Let $k$ be a given level of regularity, and assume that we are given
parameters $0 < s,\gamma,\epsilon$ with the properties that
$s\leqslant 1$, and $2\epsilon \leqslant s-\frac{1}{2}$, and
$s+\gamma < \frac{3}{2}$. Let $0 <T$ be a given time value. Then
there exists a constant $1 \leqslant C_{k,s,\gamma,\epsilon}$
depending only on $k,s,\epsilon,\gamma$, such that if:
\begin{equation}
        E_2^{(s,\gamma,\epsilon)}(0,T)[\Psi] \ \leqslant \
    C^{-1}_{k,s,\gamma,\epsilon} \ , \label{phi_boot_assump}
\end{equation}
then one has the following weighted space-time estimate for commutators:
\begin{equation}
        \sum_{\substack{|I|\leqslant k\\ X\in\mathbb{L}}}
    \lp{\tp^s \tm^\frac{1}{2} (w)^\frac{1}{2}_{\gamma,\epsilon}\
      \big[\cBox , D_X^I\big] \phi }{L^2[0,T]} \ \lesssim \
    \sum_{l=2}^7\
    \left[E_k^{(s,\gamma,\epsilon)}(0,T)[\Psi]\right]^l
    \ , \label{main_phi_comm_est}
\end{equation}
where $w_{\gamma,\epsilon}$ is the weight function from line
\eqref{w_gamma_def1.5} above.
\end{prop}\ret

\begin{proof}[Proof of estimate \eqref{main_phi_comm_est}]
As in the previous section, we may assume that
$C^{-1}_{k,s,\gamma,\epsilon}$ is chosen small enough that the right
hand side of \eqref{abs_Linfty_est} containing the extra $L^\infty$
norm can be absorbed into the left hand side. Taking this into
account, it suffices to be able to show:
\begin{multline}
        \sum_{\substack{|I|\leqslant k\\
        X\in\mathbb{L}} }\
        \lp{\tp^s\tm^\frac{1}{2} (w)^\frac{1}{2}_{\gamma,\epsilon}\
        \big[\cBox, D^I_X \big]\phi}{L^2[0,T]}\
        \lesssim \\
        E_k^{(s,\gamma,\epsilon)}(0,T)[\Psi]\cdot
        \llp{\mathcal{L}_\mathbb{L}^{k-2} \Psi}
        {L^\infty[0,T](s,\gamma,\epsilon)}^2\left( 1 +
        \llp{\mathcal{L}_\mathbb{L}^{k-2} \Psi}
        {L^\infty[0,T](s,\gamma,\epsilon)}^2
        \right)^2 \ . \label{phi_norm_est_Linfty}
\end{multline}
Our first step in this process is to expand the commutator
$\big[\cBox, D^I_X \big]$. This is done through the following
multiindex identity:
\begin{equation}
        [\cBox, D^I_X ] \phi \ = \ \sum_{\substack{K_1 + K_2 + K_3 = I\\ |K_2|=1\\
        X_1,X_2,X_3 \in \mathbb{L} }}\
        D_{X_1}^{K_1} [\cBox , D^{K_2}_{X_2} ]D^{K_3}_{X_3} \phi \ .
        \label{multicom_form}
\end{equation}
Now, using the single derivative commutator formula
\eqref{cbox_comm_formula} in conjunction with a repeated application
the two identities:
\begin{multline}
        D_X\left(Y^\alpha F_{\alpha\beta} D^\beta\phi \right) \ = \
        [X,Y]^\alpha F_{\alpha\beta}D^\beta\phi \\
        + \ Y^\alpha \mathcal{L}_X F_{\alpha\beta}D^\beta\phi \ + \
        Y^\alpha F_{\alpha\beta}D^\beta(D_X\phi) \\
        - \ Y^\alpha F_{\alpha\beta}\ {}^{(X)}\pi^{\beta\gamma}D_\gamma
        \phi \ + \ \sqrt{-1} Y^\alpha F_{\alpha\beta}
        F_\gamma^{\ \beta}
        X^\gamma \cdot \phi \ , \notag
\end{multline}
and:
\begin{multline}
        D_X\left( \nabla^\beta (Y^\alpha F_{\alpha\beta})\cdot \phi
        \right) \ = \ \nabla^\beta ([X,Y]^\alpha F_{\alpha\beta})\cdot
        \phi\
        + \ \nabla^\beta ( Y^\alpha \mathcal{L}_X F_{\alpha\beta}
        )\cdot \phi \\ + \ \nabla^\beta (Y^\alpha F_{\alpha\beta})\cdot
        D_X \phi\ - \ {}^{(X)}\pi^{\beta\gamma} \nabla_\gamma (Y^\alpha
        F_{\alpha\beta})\cdot \phi \ , \notag
\end{multline}
which in turn follow from the geometric formulas
\eqref{dagger_form}, \eqref{F_div_iden},
\eqref{complex_lie_cont_id}, \eqref{D_complex_lie}, and
\eqref{div_omega_calc} (respectively), we arrive at the following
point-wise bound for the left hand side of \eqref{multicom_form}
above:
\begin{align}
        &\sum_{\substack{|I|\leqslant k\\ X\in\mathbb{L}} } \ | [\cBox, D^I_X
        ]\phi| \ , \label{phi_ABC_bound}\\
        \lesssim \ &\sum_{\substack{ |K_1| + |K_2| \leqslant k-1\\
        X_1 , X_2 , Y \in \mathbb{L}}} \ \Big(
        \big| Y^\beta ( \mathcal{L}_{X_1}^{K_1} J_\beta) \big|\cdot
        \big|D_{X_2}^{K_2} \phi \big| \notag\\
        &\ \ \ \ \ \ \ + \ \big| 2 Y^\alpha(
        \mathcal{L}_{X_1}^{K_1} F_{\alpha\beta})\cdot D^\beta(D_{X_2}^{K_2}\phi)
        - \nabla^\alpha(Y^\beta)\cdot( \mathcal{L}_{X_1}^{K_1} F_{\alpha\beta})\cdot
        (D_{X_2}^{K_2}\phi) \big| \Big) \notag \\
        & \ \ \ \ \ \ \ \ \ \ \ + \ \sum_{\substack{ |K_1| +
    |K_2| +|K_3| \leqslant k-2\\
        X_1 , X_2 , X_3 , Y_1 , Y_2 \in \mathbb{L}}}\
        \big|Y_1^\alpha (\mathcal{L}_{X_1}^{K_1} F_{\alpha\beta})\big|\cdot
         \big|Y_2^\gamma (\mathcal{L}_{X_2}^{K_2} F^{\beta}_{\ \gamma})\big|\cdot
        \big| D_{X_3}^{K_3}\phi\big| \ , \notag\\
        = \ &A + B + C \ . \notag
\end{align}
Our goal in now to prove the bound \eqref{phi_norm_est_Linfty} for
each of these three terms. We do this separately and in order. As
with the proof of \eqref{JL_parity_est}--\eqref{sJ_parity_est}
above, we will only carry out our calculations in the region $t<2r$
as the calculations necessary to handle the region $r <
\frac{1}{2}t$ are essentially trivial because all components of
$(F,D\phi)$ can be treated on an equal footing there.\\

Our first step is to deal with  the $A$ term above. Taking in to
account the abstract parity estimates
\eqref{abs_parity_est}--\eqref{charge_abs_parity_est}, our task is
reduced to showing the symbolic bound:
\begin{multline}
        A \ \lesssim \\
    \sum_{\substack{l+m =k-1\\ a+b+c +d = 0 }} \ \tp^{-1}\ \tau_{(a)}
    \tau_{(b)} \ \left(
        \mathcal{L}^l_\mathbb{L}\Psi_{(c)}+
        \mathcal{L}^l_\mathbb{L} Q_{(c)}
        \right)\cdot \mathcal{L}^m_\mathbb{L} \Psi_{(d)}\cdot\left( 1 +
    \llp{\mathcal{L}_\mathbb{L}^{k-2}\Psi}{L^\infty[0,T](s,\gamma,\epsilon)}
        \right)^2 \ . \label{A_abs_parity_est}
\end{multline}
Expanding out the contraction with respect to $\beta$ in the first
factor of the $A$ term, and using the identifications
\eqref{null_parity_numbers} we have the straight forward bound:
\begin{equation}
        \big| Y^\beta \mathcal{L}_{X_1}^{K_1} J_\beta \big|
        \ \lesssim \ \sum_{a + b =0} \ \tau_{(a)}\
        \big| \mathcal{L}_{X_1}^{K_1} J_{(b)} \big| \
        . \label{J_contr_parity_bound}
\end{equation}
To bound the $\mathcal{L}_{X_1}^{K_1} J_{(b)}$ term itself, we can
simply use lines \eqref{JL_parity_est}--\eqref{sJ_parity_est} above.
However, it will be necessary for us to use the following refinement
of those bounds, which can easily be checked by recalling lines
\eqref{AF_bound1}--\eqref{BF_bound}:
\begin{multline}
        \sum_{\substack{ |K_1| \leqslant l \\ X_1 \in \mathbb{L}}}\
        \big| \mathcal{L}_{X_1}^{K_1} J_{(b)} \big| \ \lesssim \\
        \sum_{\substack{ c + d = b\\ l_1 + l_2 = l}}\ \tp^{-1}\
        \tau_{(c)} \ \left(
        \mathcal{L}^{l_1}_\mathbb{L}\Psi_{(d)}+
        \mathcal{L}^{l_1}_\mathbb{L} Q_{(d)}
        \right)\cdot\sum_{\substack{|K_1| \leqslant l_2 \\ X_1 \in \mathbb{L} }}\
        | D^{K_1}_{X_1} \phi |\
        \cdot \left( 1 +
        \llp{\mathcal{L}_\mathbb{L}^{l-2}\Psi}{L^\infty[0,T](s,\gamma,\epsilon)}
        \right) \ . \notag
\end{multline}
Expanding this into \eqref{J_contr_parity_bound} above, and
tacking on the extra factor of $D_{X_2}^{K_2}\phi$ while using the bound:
\begin{equation}
        \sum_{\substack{ |K_1| \leqslant l_2 , |K_2| \leqslant m \\
        l_2 + m \leqslant k-1 \\
        X_1 , X_2 \in
        \mathbb{L}}}\ |D_{X_1}^{K_1}\phi | \cdot  |D_{X_2}^{K_2}\phi |
    \ \lesssim \  \mathcal{L}_\mathbb{L}^{\max\{l_2 , m
    \}}\Psi_{(0)}\cdot \llp{\mathcal{L}_\mathbb{L}^{k-2}\Psi}
    {L^\infty[0,T](s,\gamma,\epsilon)} \ , \notag
\end{equation}
noting that by design we have $l_1 + \max \{l_2 , m\} \leqslant
k-1$, we see that we have achieved \eqref{A_abs_parity_est}. Notice
that here we are making (crucial) use of the special bound
$|D_X^I\phi|\lesssim \tp^{-1}$, when $|I|\leqslant k-1$, contained
in the norm \eqref{abs_Linfty_norm}.\\

We now move on to bounding the $B$ term in line \eqref{phi_ABC_bound} above.
Here we will prove the symbolic bound:
\begin{equation}
        B \ \lesssim \ \sum_{\substack{ a + b + c = 0 \\ l + m
        \leqslant k-1}} \ \tau_{(a)}\ \left(
        \mathcal{L}_\mathbb{L}^l \Psi_{(b)} +
        \mathcal{L}_\mathbb{L}^l Q_{(b)}
        \right)\cdot \mathcal{L}_\mathbb{L}^m \Psi_{(c)}
    \ . \label{Bphi_parity_bound}
\end{equation}
Through an application of estimate \eqref{abs_parity_est}, this will
prove the bound \eqref{phi_norm_est_Linfty} for the $B$  portion of
things. To simplify matters, it suffices to prove the above bound
with $l=m=0$. The same estimate containing derivatives is simply a
matter of notation. With this in mind, it suffices to show that
both:
\begin{align}
        \big| Y^\alpha ( F_{\alpha\beta})\cdot  \frac{1}{r}D^\beta(r\phi)
    \big| \ &\lesssim \ \sum_{a+b + c= 0} \ \tau_{(a)} \
    \left(Q_{(b)} + \Psi_{(b)} \right)\cdot \Psi_{(c)}
    \ , \label{main_psiB_bound1}\\
    \big| \big(\frac{2}{r} \nabla^\alpha(r)\cdot X^\beta - \nabla^\alpha(X^\beta)
    \big)\cdot F_{\alpha\beta}\big|\cdot |\phi| \ &\lesssim \
    \sum_{a+b + c= 0} \ \tau_{(a)} \
    \left(Q_{(b)} + \Psi_{(b)} \right)\cdot \Psi_{(c)}
    \ . \label{main_psiB_bound2}
\end{align}
Notice that the  proof \eqref{main_psiB_bound1} is a simple matter
of applying the definitions \eqref{abs_def}, \eqref{parity_optical},
and the weight bounds \eqref{lor_comp_bound}. The key observation
here is  that since this is a full contraction, its parity weight
must be zero. Also, note that the term $\frac{1}{r} D_{\bL} (r\phi)$
can be treated as a $\Psi_{(-1)}$ term on account of the symbolic
bounds (i.e. any bound for the right hand side is satisfied by the
left hand side):
\begin{equation}
        \Psi_{(0)} \ \lesssim \ \Psi_{(-1)} \ . \notag
\end{equation}
Therefore, to prove \eqref{Bphi_parity_bound} it now remains to show
the bounds \eqref{main_psiB_bound2}. For the most part this is a
simple matter of noticing that from the formulas
\eqref{cov_ders_vect_fields} (or by homogeneity!), each component of
the covariant derivatives on the left hand side  of
\eqref{main_psiB_bound2} satisfies the bounds:
\begin{equation}
        | \frac{2}{r} \nabla^\alpha(r) \cdot X^\beta | +
    | \nabla^\alpha X^\beta | \ \lesssim \ 1 \ . \notag
\end{equation}
Thus, if the parity weight of $F_{\alpha\beta}$ is $(0)$ or $(1)$ we
can pass to the right hand side of \eqref{main_psiB_bound2} after
multiplication by the bounded factor $\tp\cdot r^{-1}$. The only
case where this general procedure does not work is when
$F_{\alpha\beta}$ has weight $(-1)$. In this case, to pass to the
right hand side of \eqref{main_psiB_bound2} we must pick up an extra
factor of $\tz$. This indeed turns out to be the case, and is
perhaps the most striking structural property of the commutator
\eqref{cbox_comm_formula} as was pointed out in \cite{Shu_MKG} (see
p. 226 of that work) in a somewhat different form. What we need to
show is that:
\begin{equation}
        \big|\frac{2}{r} \nabla^{\bL}(r)\cdot X^A - \nabla^{\bL}(X^A)
        + \nabla^A(X^{\bL}) \big| \lesssim \ \tz \ , \label{special_canc}
\end{equation}
for each $X\in\mathbb{L}$ in the region $t<2r$. This follows from a
direct use of the identities \eqref{cov_ders_vect_fields}. Notice
that for the case $X\in\{\partial_\alpha , S\}$, the above sum
either vanishes or is of the order $O(r^{-1})$. Thus, the main thing
to check is that \eqref{special_canc} holds when $X\in\{\Omega_{ij}
, \Omega_{i0}\}$. This in turn follows easily from the formulas:
\begin{align}
        \nabla^{\bL}(\Omega_{ij}^A)\ &= \ - \frac{1}{2r}\
        \Omega_{ij}^A \ ,
        &\nabla^{A}(\Omega_{ij}^{\bL}) \ &= \ \frac{1}{2r}\ \Omega_{ij}^A \
        , \notag\\
        \nabla^{\bL}(\Omega_{i0}^A)\ &= \  - \ \frac{1}{2}\
        \omega_i^A \ ,
        &\nabla^{A}(\Omega_{i0}^{\bL}) \ &= \ \frac{1}{2}\ \omega_{i}^A \
        . \notag
\end{align}\\

To finish our proof of estimate \eqref{phi_norm_est_Linfty}, we only
need to show this bound for the $C$ term on line \eqref{phi_ABC_bound}
above. Writing $F = \tF + \overline{F}$ and
expanding out this term into parity notation, we have the
symbolic bound:
\begin{align}
        &C \ \lesssim \notag \\
        &\sum_{\substack{l+m =k-2\\ a+b+c + d= 0 }} \ \tp^{-1}\
        \tau_{(a)} \tau_{(b)}\
        \left(
        \mathcal{L}^l_\mathbb{L} \Psi_{(c)} +
        \mathcal{L}^l_\mathbb{L} Q_{(c)}
        \right)\cdot \mathcal{L}^m_\mathbb{L} \Psi_{(d)}\cdot
        \llp{\mathcal{L}_\mathbb{L}^{k-2}\Psi}{L^\infty[0,T](s,\gamma,\epsilon)}
        \  \notag\\
        &\ \ \ \ \ \ \ \ \ \ + \
        \sum_{\substack{l+m =k-2\\ a+b+c + d= 0 }} \ \tp \
        \tau_{(a)} \tau_{(b)}\
        \mathcal{L}^{l}_\mathbb{L} Q_{(c)}\cdot
        \mathcal{L}^{m}_\mathbb{L} Q_{(d)}\cdot
        \mathcal{L}^{k-2}_\mathbb{L} \Psi_{(0)} \ , \notag \\
        & \ \ \ \ = \ C_1 + C_2 \ . \notag
\end{align}
To prove the bound \eqref{phi_norm_est_Linfty} for the term $C_1$ is
a simple matter of referring to the abstract parity estimate
\eqref{abs_parity_est}. Thus, we are left with bounding the term
$C_2$. This reduces directly to a special case of the estimate
\eqref{charge_abs_parity_est} after an application of the $L^\infty$
bound \eqref{abs_charge_Linfty} which easily implies that their
exists a cutoff function $\td{\chi}^+$ such that:
\begin{multline}
        \sum_{\substack{l+m =k-2\\ a+b+c + d= 0 }} \ \tp \
        \tau_{(a)} \tau_{(b)}\
        \mathcal{L}^{l}_\mathbb{L} Q_{(c)}\cdot
        \mathcal{L}^{m}_\mathbb{L} Q_{(d)} \ \lesssim
        \tp \cdot |q|\frac{\td{\chi}^+(r-t-2)}{r^2}\cdot
        |q| \ . \notag
\end{multline}
Notice that the right hand side of this last line behaves the same
in estimate \eqref{charge_abs_parity_est} as the expression:
\begin{equation}
        \tp^{-1}\ \tau_{(0)} \tau_{(0)}\
        \mathcal{L}^{k-2}_\mathbb{L} Q_{(0)} \cdot
        \llp{\mathcal{L}_\mathbb{L}^{k-2}\Psi}
        {L^\infty[0,T](s,\gamma,\epsilon)}
        \ . \notag
\end{equation}
Multiplying this by the factor $\mathcal{L}^{k-2}_\mathbb{L}
\Psi_{(0)}$ and then substituting  the result into left hand side of
\eqref{phi_norm_est_Linfty}, we are done by an application of the
parity estimate \eqref{abs_parity_est}. This completes the proof of
\eqref{main_phi_comm_est}.
\end{proof}

\ret
%-------------------------------------------------------------------------
%%%%%%%%%%%%%%%%%%%%%%%%%%%%%%%%%%%%%%%%%%%%%%%%%%%%%%%%%%%%%%%%%%%%%%%%%%
%-------------------------------------------------------------------------

%\section{$L^2$ Scattering: The use of Modified Coulomb Gauges}

\ret
%-------------------------------------------------------------------------
%%%%%%%%%%%%%%%%%%%%%%%%%%%%%%%%%%%%%%%%%%%%%%%%%%%%%%%%%%%%%%%%%%%%%%%%%%
%-------------------------------------------------------------------------

\section{Appendix}\label{appendix}

This appendix contains the proofs of several more or less standard
weighted Sobolev type estimates which are used at various places in
the main paper. In particular in Sections \ref{F_L2_section} and
\ref{F_Linfty_sect}--\ref{phi_Linfty_sect}. We make no claim of
originality here, rather our purpose is to have things stated in the
precise form in which we find them convenient to use in our work.
The first such estimate is a simple weighted version of the usual
$L^2
\hookrightarrow L^6$ embedding:\\

\begin{lem}\label{weighted_grad_sob_lem}
Let $\varphi$ be a real valued test function on $\RR$, and set:
\begin{equation}
        q \ = \ \int_{\RR} \ \varphi \ dx \ , \label{charge_varphi_def}
\end{equation}
the average of $\varphi$. Let $\delta$ be given such that
$\frac{1}{2} < \delta < \frac{3}{2}$. Then the following weighted
integral inequality holds:
\begin{equation}
        \int_{\RR} \ r^{2\delta}\, \Big| \nabla
        \left(\frac{1}{\Delta} \varphi + \frac{q}{4\pi
        r}\right)\Big|^2 \ dx \ \lesssim \ \lp{r^\delta\,
        \varphi}{L^\frac{6}{5}}^2 \ . \label{weighted_grad_sob}
\end{equation}
In the above estimate, the implicit constants depend on $\delta$.
\end{lem}\ret

\begin{rem}
In effect, what estimate \eqref{weighted_grad_sob} shows is that is
possible to commute the weight $r^\delta$ past the Riesz operator
$\nabla\Delta^{-1}$ so long as one first subtracts off the leading
order term in the asymptotic expansion of $\Delta^{-1}\varphi$. This
is:
\begin{equation}
        \nabla_i \frac{1}{\Delta}\varphi \ \sim
    \  \frac{q\cdot\omega_i}{4\pi r^2} \ . \label{leading_order}
\end{equation}
That such a correction is necessary is clear from the range of
weights we are considering because for $\frac{1}{2} < \delta$ a
function with this asymptotic behavior at infinity cannot be in
$L^2$ weighted by $r^\delta$. The idea behind
\eqref{weighted_grad_sob} is, of course, that once the right hand
side of \eqref{leading_order} is subtracted off from
$\nabla\Delta^{-1}\varphi$ the resulting function will decay enough
to be in $r^\delta$ weighted $L^2$ as long as the appropriately
weighted $L^\frac{6}{5}$ norm of $\varphi$ is bounded. However,
there is a limit to how much weight one can expect to apply this
way. This is because one can further expand \eqref{leading_order} to
include higher order terms:
\begin{equation}
        \nabla_i \frac{1}{\Delta}\varphi \ \sim
        \  -\ \nabla_i\left( \frac{q}{4\pi r} +
        \sum_k \frac{q_k\cdot \omega_k}{4\pi r^2} +
        \{higher\}\right) \ , \label{multipole}
\end{equation}
where:
\begin{equation}
        q_k \ = \ \int_{\RR}\ y_k\, \varphi(y)\ dy \ , \notag
\end{equation}
are  the higher moments of $\varphi$. Without the additional
vanishing of these other quantities, it is clear that one cannot put
$r^\delta\, \nabla\Delta^{-1}\varphi$ in $L^2$ for
$\frac{3}{2}\leqslant \delta$. \eqref{multipole} is sometimes
referred to as the \emph{multipole} expansion of $\nabla
\frac{1}{\Delta}\varphi$. Since we are only interested in the decay
of the initial data \eqref{initial_data} modulo charge of the order no
greater than $O(r^{-3})$, the first term in \eqref{multipole} will
be the only one which concerns us in this work.
\end{rem}\ret

\begin{proof}[Proof of \eqref{weighted_grad_sob}]
Our first step here is to simply integrate by parts several times on
the left hand side of \eqref{weighted_grad_sob}. This yields the
identity:
\begin{multline}
        \hbox{(L.H.S.)}\eqref{weighted_grad_sob} \ = \
    \delta(2\delta+1)\ \int_{\RR}\ r^{2\delta-2}\
    \left(\frac{1}{\Delta} \varphi + \frac{q}{4\pi
        r}\right)^2\ dx \\
    - \ \int_{\RR} \ r^{2\delta} \ \varphi\cdot
    \left(\frac{1}{\Delta} \varphi + \frac{q}{4\pi
        r}\right)\ dx \ . \label{lhs_wsob_int}
\end{multline}
To bound the integral in the first term on the right hand side of
\eqref{lhs_wsob_int} above, we use the definition
\eqref{charge_varphi_def} to compute:
\begin{equation}
        \Big| \int_{\RR}\ r^{2\delta-2}\
        \left(\frac{1}{\Delta} \varphi + \frac{q}{4\pi
        r}\right)^2\ dx \Big| \ \lesssim  \
        \int_{\RR}\ r^{2\delta-2}\ \left(\int_{\RR}\ \big|\frac{1}{|x-y|}-
        \frac{1}{|y|}\big|\ |\varphi(y)|\ dy \right)^2 dx \ . \notag\\
\end{equation}
We now split cases depending on whether $\frac{1}{2} < \delta
\leqslant 1$ or $1 < \delta < \frac{3}{2}$. In the first case we use
the bound:
\begin{equation}
        \big|\frac{1}{|x-y|}- \frac{1}{|y|}\big| \ \lesssim \
        |y|^\delta \left(\frac{1}{|x|^\delta\, |x-y|} +
        \frac{1}{|x|\, |x-y|^\delta}\right) \ . \notag
\end{equation}
This reduces our work to estimating the two integrals:
\begin{align}
        A \ &= \ \int_{\RR}\ \frac{1}{|x|^2}\left(\int_{\RR}\
        \frac{1}{|x-y|}\ |y|^\delta\, |\varphi(y)| \ dy\right)^2\
        dx \ , &\frac{1}{2} \ < \ \delta \ \leqslant \ 1 \ ,
        \notag\\
        B \ &= \ \int_{\RR}\ \frac{1}{|x|^{4-2\delta}}\left(\int_{\RR}\
        \frac{1}{|x-y|^\delta}\ |y|^\delta\, |\varphi(y)| \ dy\right)^2\
        dx \ , &\frac{1}{2} \ < \ \delta \ \leqslant \ 1 \ .
        \notag
\end{align}
In the second case, we simply use the bound:
\begin{equation}
        \big|\frac{1}{|x-y|}- \frac{1}{|y|}\big| \ \leqslant \
        \frac{|y|}{|x|\, |x-y|} \ . \notag
\end{equation}
This reduces us to bounding the integral:
\begin{align}
        C \ &= \ \int_{\RR}\ \frac{1}{|x|^{4-2\delta}}\left(\int_{\RR}\
        \frac{1}{|x-y|\, |y|^{\delta-1}}\ |y|^\delta\, |\varphi(y)| \ dy\right)^2\
        dx \ , &1 \ < \ \delta \ < \ \frac{3}{2} \ . \notag
\end{align}
In all cases, the bound we wish to prove is:
\begin{equation}
        A ,  B , C , \ \lesssim \ \lp{r^\delta\, \varphi
        }{L^\frac{6}{5}}\ . \label{ABC_bound}
\end{equation}
This last estimate follows from several different instances of the
generalized fractional integration Lemma \ref{frac_int_lem}
below. Notice that in the case of integral $A$ we have (notation of
Lemma \ref{frac_int_lem})
$\alpha=\beta=1$, and $\gamma=0$. In case of integral $B$ we have
$\alpha=2-\delta$, $\beta = \delta$, and  $\gamma = 0$. In the case of
integral $C$ above we have $\alpha=2-\delta$, $\beta=1$, and
$\gamma=\delta-1$. In all cases we have that $p=p'=2$ and
$q=\frac{6}{5}$ and $\alpha+\beta+\gamma =2$ so the scaling condition
\eqref{frac_int_lem_sc_cond} is satisfied. Also, note that in each
case we have that $\alpha < \frac{3}{2}$ and $\gamma < \frac{1}{2}$
so the ``gap'' condition \eqref{frac_int_lem_gap_cond} is
satisfied. This completes the proof of the bound:
\begin{equation}
        \Big| \int_{\RR}\ r^{2\delta-2}\
    \left(\frac{1}{\Delta} \varphi + \frac{q}{4\pi
        r}\right)^2\ dx \Big| \ \lesssim  \ \lp{r^\delta\, \varphi
        }{L^\frac{6}{5}}\ . \label{wsob_first_term}
\end{equation}\ret

It remains to bound the second term on the right hand side expression
\eqref{lhs_wsob_int} above. By H\"olders inequality we have that:
\begin{equation}
        \Big| \int_{\RR} \ r^{2\delta} \ \varphi\cdot
    \left(\frac{1}{\Delta} \varphi + \frac{q}{4\pi
        r}\right)\ dx\ \Big| \ \leqslant \
    \lp{r^\delta\, \varphi}{L^\frac{6}{5}}\cdot
    \lp{r^\delta \left(\frac{1}{\Delta} \varphi + \frac{q}{4\pi
        r}\right)}{L^6} \ . \notag
\end{equation}
By the usual $L^2 \hookrightarrow L^6$ Sobolev embedding, the Leibnitz
rule, and the triangle inequality, the second factor on the right hand
side above can be bounded by:
\begin{equation}
        \lp{r^\delta \left(\frac{1}{\Delta} \varphi + \frac{q}{4\pi
        r}\right)}{L^6} \ \lesssim \
    \lp{r^{\delta-1} \left(\frac{1}{\Delta} \varphi + \frac{q}{4\pi
        r}\right)}{L^2} + \lp{r^\delta \, \nabla\left(\frac{1}{\Delta}
    \varphi + \frac{q}{4\pi r}\right)}{L^2} \ . \notag
\end{equation}
Therefore, using the bound \eqref{wsob_first_term} above, we can
estimate:
\begin{multline}
        \Big| \int_{\RR} \ r^{2\delta} \ \varphi\cdot
    \left(\frac{1}{\Delta} \varphi + \frac{q}{4\pi
        r}\right)\ dx\ \Big| \\
    \lesssim \ \lp{r^\delta\,
    \varphi}{L^\frac{6}{5}}\cdot\left(
    \lp{r^\delta\, \varphi}{L^\frac{6}{5}} +
    \lp{r^\delta \, \nabla\big(\frac{1}{\Delta}
    \varphi + \frac{q}{4\pi r}\big)}{L^2}
    \right) \ . \label{wsob_second_term}
\end{multline}
Adding estimates \eqref{wsob_first_term} and
\eqref{wsob_second_term} into the identity \eqref{lhs_wsob_int},
dividing through by the quantity $\lp{r^\delta \, \nabla\big(\frac{1}{\Delta}
\varphi + \frac{q}{4\pi r}\big)}{L^2}$ and resquaring, we have
achieved the desired estimate \eqref{weighted_grad_sob}.
\end{proof}\ret

\noindent To complete the proof of \eqref{weighted_grad_sob}, we need
to show the bound \eqref{ABC_bound} for the integrals $A$,$B$, and $C$
above. This will be a consequence of the following generalization of
the classical Hardy-Littlewood-Sobolev fractional integration theorem
(see \cite{LL}):\\

\begin{lem}[General fractional integration lemma]\label{frac_int_lem}
Let $1 < q < p < \infty$ and $0\leqslant \alpha,\beta,\gamma$ be given
parameters. Then the following integral estimate holds for (positive)
test functions $F$ and $G$ on $\RR$:
\begin{equation}
        \int_{\RR}\int_{\RR} \ G(x) \ \frac{1}{|x|^\alpha\,
    |x-y|^\beta\, |y|^\gamma }\ F(y)\ dxdy \ \lesssim \
    \lp{G}{L^{p'}}\cdot\lp{F}{L^q} \ , \label{frac_int_th}
\end{equation}
where the following conditions hold on the various indices:
\begin{align}
        3-(\alpha+\beta+\gamma) \ &= \  3(\frac{1}{q} - \frac{1}{p}) \ ,
    &\hbox{(scaling)} \ , \label{frac_int_lem_sc_cond}\\
    \alpha \ < \ 3(1-\frac{1}{p}) \ , \
    \gamma \ &< \ 3(1-\frac{1}{q}) \ , &\hbox{(``gap'')} \ .
    \label{frac_int_lem_gap_cond}
\end{align}
The implicit constant in the above estimate depends on $p,q,\alpha,\beta,\gamma$.
\end{lem}\ret

\begin{proof}[Proof of \eqref{frac_int_th}]
By the restricted weak type version of the Marcinkiewicz interpolation
theorem, it suffices to prove that:
\begin{equation}
        \int_{\RR}\int_{\RR} \ \chi_E(x) \ \frac{1}{|x|^\alpha\,
    |x-y|^\beta\, |y|^\gamma }\ \chi_F(y)\ dxdy \ \lesssim \
    |E|^\frac{1}{p'}\cdot|F|^\frac{1}{q} \ , \notag
\end{equation}
for measurable sets $E$ and $F$ where $1 < q < p < \infty$ and the
conditions
\eqref{frac_int_lem_sc_cond}--\eqref{frac_int_lem_gap_cond} hold. By
the Riesz rearrangement inequality (see \cite{LL}) this is further
reduced to showing that:
\begin{equation}
        \int_{\RR}\int_{\RR} \ (\chi_E(x)\cdot |x|^{-\alpha})^*
        \ \frac{1}{
        |x-y|^\beta}\ (\chi_F(y)\cdot |y|^{-\gamma})^* \ dxdy \ \lesssim \
        |E|^\frac{1}{p'}\cdot|F|^\frac{1}{q} \ , \label{frac_int_red}
\end{equation}
where $f^*$ denotes the symmetric decreasing
rearangement\footnote{Here we use a definition of $f^*$ which does
not include the usual normalization factor of $4\pi$ to avoid
additional typesetting. This is $\chi^*_E(r) =
\chi_{[0,|E|^\frac{1}{3})}$. Clearly this convention does not effect
the use of rearrangements in the Riesz inequality.} of the function
$f$. We now compute that:
\begin{align}
        (\chi_E(x)\cdot |x|^{-\alpha})^*(r) \ &= \ \int_0^\infty\
        \chi^*_{ \{ \chi_E\cdot |x|^{-\alpha} > s\} }(r) \ ds \ ,
        \notag\\
        &= \ \int_0^\infty \
        \chi_{  [0 , |E\cap  B( s^{-\frac{1}{\alpha}} )|^\frac{1}{3} )    }    (r)
        \ ds \ , \notag\\
        &\leqslant \ \chi_{[0,|E|^\frac{1}{3})}(r)
        \cdot \int_0^\infty\ \chi_{[0,4\pi s^{-\frac{1}{\alpha}})}(r)\ ds \
        , \notag\\
        &= \ (4\pi)^{-1} (\chi_E)^*(r)\cdot r^{-\frac{1}{\alpha}}
        \ , \notag
\end{align}
and similarly for  $(\chi_F(x)\cdot |x|^{-\gamma})^*$. We now plug
these formulas into the right hand side of \eqref{frac_int_red} and
apply the usual fractional integration theorem to yield:
\begin{multline}
        \int_{\RR}\int_{\RR} \ (\chi_E(x))^*\cdot |x|^{-\alpha}
    \ \frac{1}{
    |x-y|^\beta}\ (\chi_F(y))^*\cdot |y|^{-\gamma} \ dxdy \\
    \lesssim \ \lp{(\chi_E(x))^*\cdot |x|^{-\alpha}}{L^{\theta_1}}
    \cdot\lp{  (\chi_F(x))^*\cdot |x|^{-\gamma}  }{L^{\theta_2}} ,
    \label{usual_frac_int}
\end{multline}
where we choose $\theta_1$ and $\theta_2$ to be given by the
expressions $\frac{1}{\theta_1} = \frac{1}{p'} + \frac{\alpha}{3}$
and $\frac{1}{\theta_2} = \frac{1}{q} + \frac{\gamma}{3}$
respectively. Notice that from the scaling condition
\eqref{frac_int_lem_sc_cond} and some quick algebra we have that:
$6-\beta \ = \ 3(\frac{1}{\theta_1} + \frac{1}{\theta_2})$, which is
the scaling condition needed for \eqref{usual_frac_int} to hold.
Furthermore, notice that by the gap condition
\eqref{frac_int_lem_gap_cond} we have that both $1 <
\theta_1,\theta_2 < \infty$ which is also needed for
\eqref{usual_frac_int} to hold. The proof of \eqref{frac_int_red} is
now accomplished by simply computing:
\begin{align}
        \lp{(\chi_E(x))^*\cdot |x|^{-\alpha}}{L^{\theta_1}}^{p'} \ &=
        \ \Big(\int_0^{|E|^\frac{1}{3}} \ r^{2-\alpha\theta_1}\ dr
        \Big)^\frac{p'}{\theta_1} \ , \notag\\
        &\lesssim \ |E|^{\frac{p'}{3\theta_1}(3-\alpha\theta_1)} \ ,
        \notag\\
        &= \ |E| \ . \notag
\end{align}
Notice that the inequality follows because we automatically have that
$\alpha\theta_1 < 3$. An identical computation shows that:
\begin{equation}
        \lp{(\chi_F(x))^*\cdot |x|^{-\gamma}}{L^{\theta_2}}^{q} \
        \lesssim \ |F| \ . \notag
\end{equation}
This completes the proof of \eqref{frac_int_red} and therefore the proof of
estimate \eqref{frac_int_th}.
\end{proof}\ret

The second main set of estimates we prove here are a localized
version of the so called global Sobolev inequalities. These were
first utilized in \cite{K_semi} to prove global existence and
asymptotic behavior for general non-linear wave equations. The
versions which we state here are sufficiently ``atomic'' to provide
all of the $L^\infty$ we will need in this paper. To discuss these,
we shall employ the same notation as introduced in the beginning of Section
\ref{F_Linfty_sect}. In particular, the notion of a dyadic cone
distance (CD) shell and the
associated $\tp(\mathcal{J}),\tm(\mathcal{J})$ notation.\\

\begin{lem}[Exterior global Sobolev estimate]\label{basic_decay_lem1}
In the exterior region $1\leqslant t< 2r$ let $\mathcal{J}$ be a
given dyadic CD shell, and let $q$ be given such that $2 < q <
\infty$. Then for test functions $f$ the following estimates hold:
\begin{equation}
        \sup_\omega\ |\chi f(t,r,\omega)| \ \lesssim \
        \tp^{-\frac{2}{q}}(\mathcal{J})\ \sum_{ i < j}\
        \lp{\Omega_{ij}  \ \chi f(t,r,\cdot)}{L^q(\mathbb{S}_r^2)}
        \ , \label{ext_sob1}
\end{equation}
where $\chi$ is the cutoff on the dyadic CD shell $\mathcal{J}$.
Furthermore, in the case where $2\leqslant q < 4$ we also have the
estimate:
\begin{multline}
        \lp{ \chi f }
        {L^\infty_r(L^q(\mathbb{S}_r^2))} \ \lesssim \\
        \tp^{-2(\frac{1}{2} - \frac{1}{q})}(\mathcal{J})\,
        \tm^{-\frac{1}{2}}(\mathcal{J})
        \Big( \lp{\chi f}{L^2} \ +  \
        \lp{ \tm \, \partial_r \ \chi f }{L^2}
        \ +  \  \sum_{ i < j }\ \lp{ \Omega_{ij} \chi f }{L^2}
        \Big) \ . \label{ext_sob2}
\end{multline}
\end{lem}\ret

\begin{proof}[Proof of estimates \eqref{ext_sob1}--\eqref{ext_sob2}]
The proof of both of these is entirely standard as they are
essentially just rescaled versions of the usual translation
invariant Sobolev estimates. For example, \eqref{ext_sob1} is just
the Sobolev theorem on spheres $\mathbb{S}_r$. For the sake of
completeness we give here a proof of \eqref{ext_sob2}.\\

The first step is to introduce a new set of variables which we
denote by $(\td{r},\omega)$ and  are supported in the unit dyadic
annular region $\mathcal{A} = [\frac{1}{8},8]\times\mathbb{S}^2$
(where $\mathbb{S}^2$ without a subscript denotes the unit sphere),
and which satisfy the change of variable formula:
\begin{equation}
          (\tm(\mathcal{J}) \cdot \td{r} + t,r\cdot\omega) \ = \
          (r,\omega_r) \in \mathbb{R}\times \mathbb{S}^2_r
          \ . \notag
\end{equation}
That is, we rescale the standard radial variable around the point
$r=t$ by the factor $\tm{\mathcal{J}}$ \emph{and} we rescale the
angular variable by the factor of $r$. This procedure gives the
change of volume formula:
\begin{align}
        dr\, d\omega_r \ &= \ v(r)\cdot d\td{r}\, d\omega \ ,
        &v(r) \ &\sim \ \tp^2(\mathcal{J})\cdot \tm(\mathcal{J}) \ ,
        \label{volume_change}
\end{align}
where the comparison on the right hand side above holds in the
region $x\in\mathcal{J}$. Using the analog of \eqref{volume_change}
for integration on a fixed sphere we have the equivalence:
\begin{equation}
        \lp{ \chi f }{L^\infty_r(L^q(\mathbb{S}_r^2))} \ \sim \
    \tp^\frac{2}{q}(\mathcal{J})\cdot \lp{\chi
    f}{ L^\infty_{\td{r}} (L^q(\mathbb{S}^2))} \ . \notag
\end{equation}
Furthermore, a direct calculation using \eqref{volume_change} gives:
\begin{multline}
    \tp(\mathcal{J})\tm^\frac{1}{2}(\mathcal{J})\cdot\Big(
    \lp{ \chi f}{L^2} \ +  \
    \sum_{ i < j }
    \lp{ \Omega_{ij} \ \chi f}{L^2} \ +  \
    \lp{ \tm \, \partial_r \ \chi f }{L^2}\Big) \ \sim \\
    \lp{\chi f}{  L^2(\mathcal{A})} \  + \
    \lp{\partial_{\td{r}} \chi f}{L^2(\mathcal{A})} \ + \
    \lp{\snabla \chi f}{L^2(\mathcal{A})} \ . \notag
\end{multline}
Using these last two expressions in tandem we see that the proof of
\eqref{ext_sob2} reduces to being able to show that:
\begin{equation}
        \lp{\chi f}{ L^\infty_{\td{r}} (L^q(\mathbb{S}^2))} \ \lesssim \
    \lp{\chi f}{H^1(\mathcal{A})} \ . \label{a_sob}
\end{equation}
Notice that by design, the rescaled $\chi f$ does not intersect the
boundary $\partial\mathcal{A}$.
Using two charts on the interior of $\mathcal{A}$,
estimate \eqref{a_sob} follows from the general mixed
Sobolev embedding on $\RR$ for test functions $\varphi$:
\begin{equation}
        \lp{\varphi}{L_x^p(L^q_y)} \ \lesssim \ \lp{\varphi}{H^1} \ ,
    \label{general_sobolev}
\end{equation}
where $(x,y) \in \mathbb{R}\times\mathbb{R}^2$ and
$\frac{3}{2} - \frac{1}{p} - \frac{2}{q} < 1$.
This last estimate in turn follows from running
the Sobolev lemma in the $x$ and $y$ variables separately to achieve:
\begin{equation}
        \lp{\varphi}{L_x^p(L^q_y)} \ \lesssim \
    \lp{\langle D_x \rangle^{(\frac{1}{2} - \frac{1}{p}) +}
      \langle D_y\rangle^{2(\frac{1}{2} - \frac{1}{q}) }
      \varphi}{L^2} \ . \notag
\end{equation}
where the ``$+$'' notation denotes an arbitrarily small positive
constant which gets us around the case $p=\infty$.
\eqref{general_sobolev} now follows from this last estimate and the
symbol bounds:
\begin{equation}
        (1 + |\xi_x|^2)^{\frac{1}{2}(\frac{1}{2} - \frac{1}{p}) +}
        (1 + |\xi_y|^2)^{ (\frac{1}{2} - \frac{1}{q})  } \ \lesssim \
        (1 + |\xi_x|^2 + |\xi_y|^2 )^\frac{1}{2} \ , \notag
\end{equation}
whenever $\frac{3}{2} - \frac{1}{p} - \frac{2}{q} < 1$. We have now shown
\eqref{ext_sob2}.
\end{proof}\ret

The second global Sobolev estimate we prove here is the analog of
Lemma \ref{basic_decay_lem1} for the interior region $r < \frac{1}{2}
t$. This is:\\

\begin{lem}[Interior global Sobolev estimate]\label{basic_decay_lem2}
Fix $1\leqslant t$, then in the interior region $r < \frac{3}{4}t$
one has the following weighted estimate for compactly supported (in
that region) test functions $f$:
\begin{equation}
        \lp{f}{L^q} \ \lesssim \
        t^{-3(\frac{1}{p}-\frac{1}{q})}
        \Big( \lp{f}{L^p} \ + \
        \sum_{ X\in \{S,\Omega_{0i} \}}\
        \lp{X(f) }{L^p} \Big) \ , \label{int_sob}
\end{equation}
whenever $\frac{1}{p} - \frac{1}{q} < \frac{1}{3}$.
\end{lem}\ret

\begin{proof}[Proof of estimate \eqref{int_sob}]
Note that this estimate is taking place in a ball of radius $t$. By
simply rescaling the usual Sobolev estimate we have that:
\begin{equation}
        \lp{f}{L^q} \ \lesssim \
        t^{-3(\frac{1}{p}-\frac{1}{q})}
        \big( \lp{f}{L^p} \ + \ t \ \lp{\nabla f}{L^p}
        \big) \ . \notag
\end{equation}
Therefore, the desired result follows directly from the bound:
\begin{equation}
        t \ \lp{\nabla f}{L^p} \ \lesssim \ \sum _{ X\in \{S,\Omega_{0i} \}}\
        \lp{X(f) }{L^p} \ , \notag
\end{equation}
which in turn follows from the identities:
\begin{align}
        t\, \partial_i \ &= \ \Omega_{i0} - x_i\partial_t \ , \notag\\
        \partial_t \ &= \ (t^2-r^2)^{-1} \left(
        tS - x^i\Omega_{i0}\right) \ . \notag
\end{align}
\end{proof}\ret

Finally, we end this appendix with a characteristic version of the
exterior estimates \eqref{ext_sob1}--\eqref{ext_sob2}. This will be
used  to prove the peeling properties of the best decaying
components of $F$ and $\phi$ in the main part of the paper. Although
strictly speaking it is not necessary for our proof, we would like
to know that these good decay estimates hold a-priori at each fixed
time, without having to assume that the solution to
\eqref{basic_MKG} satisfies appropriate global decay estimates.
Because of this we will need to cut the characteristic estimate off
sharply inside the time slabs $1\leqslant t \leqslant t_0$.\\

\begin{lem}[Characteristic (truncated) global Sobolev
estimate]\label{basic_decay_lem3}
Let $C(u)$ denote the cones $u=const.$ and let $t_0$ be a fixed
parameter. Now define the truncated cone:
\begin{equation}
        \overline{C}(u) \ = \ C(u)\cap\{1\leqslant t
        \leqslant t_0\}\cap\{t < 2r\} \ , \notag\\
\end{equation}
Let  $\mathcal{I}$ be a dyadic shell along the extended exterior
region $t <2r$ of the cone $C(u)$, and $\chi$ its smooth cutoff.
Then for any test function $f$ one has the following estimate for $2
< q < \infty$:
\begin{equation}
        \sup_\omega \ |\chi \, f(t,r,\omega)| \ \lesssim \ \tp^{-\frac{2}{q}}(\mathcal{I})\
        \sum_{i<j} \ \lp{\Omega_{ij} \chi\,
        f(t,r,\cdot)}{L^q(\mathbb{S}_{r}^2)} \ . \label{char_sob1}
\end{equation}
Furthermore, in the case where $2\leqslant q < 4$ we also have the
estimate:
\begin{multline}\label{char_sob2}
        \lp{\chi\,
    f}{L_u^\infty(L^q(\mathbb{S}^2_r)(\overline{C}(u))} \ \lesssim \  \\
    \tp^{-\frac{3}{2} - \frac{2}{q}}(\mathcal{I})
    \ \Big(\lp{\chi\, f}{L^2(\overline{C}(u))} \ + \
    \lp{\bu L (\chi\, f)}{L^2(\overline{C}(u))} \ + \
    \sum_{i<j}\ \lp{\Omega_{ij}(\chi\, f)}{L^2(\overline{C}(u))}
    \Big) \ .
\end{multline}\ret

\begin{proof}[Proof of estimates \eqref{char_sob1}--\eqref{char_sob2}]
The first estimate \eqref{char_sob1} is of course just a restatement
of \eqref{ext_sob1}. Therefore, we concentrate on \eqref{char_sob2}.
This estimate is proved using the same steps for \eqref{ext_sob2}
above, that is through a rescaling argument followed by the  mixed
norm Sobolev estimate \eqref{general_sobolev}. However, we need to
take a little care to deal with the sharp cutoff inside the slab
$\{1\leqslant t \leqslant t_0\}$. This is done by utilizing a
extension argument for $H^1$ functions (see
\cite{S_SI}).\\

The details are as follows. We first rescale all variables along the
cone $\chi \overline{C}(u)$ by the factor $\tp(\mathcal{I})$. Doing
this, we are again reduced to showing the annular estimate
\eqref{a_sob}. However, this time it is possible for $\chi f$ to
extend past the boundary $\partial \mathcal{A}$ even though we are
only allowed to estimate its behavior inside the region
$\mathcal{A}$. We get around this problem by using the fact that
$\mathcal{A}$ is a $C^\infty$ submanifold with boundary of $\RR$,
which allows us to extend $\chi_\mathcal{A} \,\chi f$ to a function
$\td{f}$ with the properties:
\begin{align}
        \td{ f}|_\mathcal{A} \ &= \ \chi f |_\mathcal{A} \ ,
        &\lp{\td{f} }{H^1(\mathbb{R}^3)} \ &\lesssim \ \lp{\chi\,
        f}{H^1(\mathcal{A})} \ . \notag
\end{align}
Applying  estimate \eqref{general_sobolev} to the function $\td{f}$
and using this last bound yields the desired result.
\end{proof}
\end{lem}\ret

\ret
%------------------------------------------------------------------------
%%%%%%%%%%%%%%%%%%%%%%%%%%%%%%%%%%%%%%%%%%%%%%%%%%%%%%%%%%%%%%%%%%%%%%%%%%
%%%%%%%%%%%%%%%%%%%%%%%%%%%%%%%%%%%%%%%%%%%%%%%%%%%%%%%%%%%%%%%%%%%%%%%%%%
%-------------------------------------------------------------------------

%-------------------------------------------------------------------------
%%%%%%%%%%%%%%%%%%%%%%%%%%%%%%%%%%%%%%%%%%%%%%%%%%%%%%%%%%%%%%%%%%%%%%%%%%
%%%%%%%%%%%%%%%%%%%%%%%%%%%%%%%%%%%%%%%%%%%%%%%%%%%%%%%%%%%%%%%%%%%%%%%%%%
%-------------------------------------------------------------------------

\end{document}